\documentclass[11pt,leqno,twoside]{amsart}

\usepackage{amsmath, amsfonts, amssymb, amsthm, bm}
\usepackage{iftex}
\ifPDFTeX
  \usepackage[utf8]{inputenc} 
  \usepackage[T1]{fontenc} 
  \usepackage{lmodern} 
\else
  \usepackage{fontspec} 
\fi
\usepackage{accents}

\newcommand{\thickbar}[1]{\accentset{\rule{0.45em}{0.8pt}}{#1}}

\usepackage{booktabs}
\newtheorem{thm}{Theorem}[section] 
\newtheorem{rem}[thm]{Remark} 

\usepackage[scale=0.8]{geometry}

\usepackage{xcolor}
\usepackage{morefloats}
\usepackage{mathtools}
\usepackage{enumitem}
\usepackage[english]{babel}

\usepackage{caption}
\usepackage{subcaption}
\usepackage{svg}
\usepackage{mathrsfs}

\numberwithin{equation}{section}

\usepackage{tikz}
\usepackage{multirow}

\usepackage{algorithmic}
\usepackage{algorithm}

\usepackage{todonotes}

\usepackage{amsmath, amsthm}
\usepackage{amsfonts, amssymb}
\usepackage{bm}
\usepackage[colorlinks=true, allcolors=blue]{hyperref}

	\newtheorem{theorem}{Theorem}[section]
	\newtheorem*{theorem*}{Theorem}
	\newtheorem{lemma}{Lemma}[section]
	
	\newtheorem*{proposition*}{Proposition}

\theoremstyle{remark}

\theoremstyle{definition}

\newcommand{\ve}{\varepsilon}

\newcommand{\RR}{\mathbb{R}}

\newcommand{\Q}{\mathcal{Q}}

\title[Kinetic--Fluid Coupling for Single-Temperature Gas Mixtures]{A Realizability-Based Kinetic-Fluid Domain Decomposition Method for Single-Temperature Gas Mixtures}

\author[D. Caparello]{Domenico Caparello}
\address{Domenico Caparello \\
Université Côte d’Azur, CNRS, LJAD, Parc Valrose, F-06108 Nice, France\\
Department of Mathematics and Computer Science, University of Ferrara \\
}
\email{Domenico.CAPARELLO@univ-cotedazur.fr}
\email{domenico.caparello@unife.it}

\author[L. Pareschi]{Lorenzo Pareschi}
\address{Lorenzo Pareschi \\
Maxwell Institute for Mathematical Sciences and Department of Mathematics,
Heriot-Watt University \\
Department of Mathematics and Computer Science, University of Ferrara \\
}
\email{L.Pareschi@hw.ac.uk}

\author[T. Rey]{Thomas Rey}
\address{Thomas Rey \\
Université Côte d’Azur, CNRS, LJAD, Parc Valrose, F-06108 Nice, France
}
\email{Thomas.REY@univ-cotedazur.fr}



\begin{document}
\begin{abstract}
    This paper introduces a domain-adaptive numerical method for multi-species kinetic models that relax toward a single-temperature and single-velocity hydrodynamic limit. To this end, we derive moment realizability matrices for this regime. In this way the ideas of Levermore, Morokoff, and Nadiga \cite{lev-mor-nad-1998} are extended to non-equilibrium mixtures in which inter-species and intra-species collisions occur at comparable rates, so that all species relax to a common velocity and a common temperature.
    The realizability conditions are used as an indicator to divide the computational domain, at each time step, into macroscopic and kinetic regions. In the kinetic regions, where the continuum description fails (e.g., across strong shocks, non-equilibrium interfaces, or kinetic boundary layers), the method solves the multi-species kinetic equations. Close to local thermodynamic equilibrium, the Euler equations of the mixture are solved instead. Confining the kinetic computations to the non-equilibrium regions reduces the simulation time, while the results remain close to those of the full kinetic solver.
    The method is discretized on Cartesian grids with high-order schemes in space and time. The kinetic equations are integrated with an Asymptotic-Preserving Implicit-Explicit (IMEX) scheme, so that a consistent discretization of the single-velocity, single-temperature Euler system of the mixture is recovered when the Knudsen number goes to zero. To the best of our knowledge, this is the first realizability-based adaptive spatial domain decomposition method for gas mixtures in this single-velocity, single-temperature setting.
\end{abstract}

\maketitle

 \section{Introduction}\label{sec1}
    \setcounter{equation}{0}
    Gas mixtures appear in many problems, from aerospace engineering to evaporation--condensation \cite{jin-qin-2013} and plasma physics \cite{bot-cow-sch-2024}. In transitional regimes, as in the re-entry of spacecraft at high altitude \cite{dill-cor-2008} or in micro-electro-mechanical systems (MEMS), the continuum hypothesis is no longer valid. In these conditions macroscopic models such as the Navier--Stokes equations do not correctly describe localized structures like shock profiles and kinetic boundary layers, and kinetic models, such as the multi-species Boltzmann equation, are needed.

    The theory of multi-species kinetic models is well developed \cite{mar-2016, bri-dau-2016}, and more recent works also consider granular mixtures and their long-time behaviour \cite{rey-ten-2024}. These equations couple a linear transport operator with nonlinear collision terms, which become stiff in the fluid regime. Their numerical solution is expensive, because the distribution functions depend on up to seven variables \cite{dim-par-2014}. When the mean free path goes to zero, the solution of the kinetic system approaches that of the Euler equations of the mixture, which are much cheaper to solve. For this reason, using a kinetic solver in the whole domain is wasteful when the flow is close to equilibrium in most of it.

    Hybrid methods, or multi-physics domain decomposition strategies, reduce this cost \cite{dim-par-2014, bur-boy-2009, fil-rey-2014, tiw-1998, cap-par-rey-2025, cap-par-rey-2026, cap-ten-2025}. They divide the computational domain into regions where different models are solved, according to the local distance from thermodynamic equilibrium.
    
    Kinetic--fluid coupling for gas mixtures has also been studied in other works. A particle-based micro--macro method was proposed in \cite{cre-kli-pir-2020}. Domain decomposition methods for multi-species flows were proposed in \cite{dat-shu-2023}, where a multi-species kinetic model is coupled with a multi-fluid model, and, closer to the present setting, in \cite{chu-sag-mur-haa-2024}, where a multi-species Vlasov--BGK model is coupled with its single-velocity, single-temperature hydrodynamic limit through buffer regions. In both domain decomposition methods the kinetic and fluid regions are fixed before the computation and do not change in time. Here, instead, the partition is dynamic: it is updated at every time step with a realizability-based indicator, so that the kinetic and fluid regions evolve in time with the flow.
    
    The stiffness of the collision terms is handled with Asymptotic-Preserving (AP) schemes \cite{cor-per-1991, jin-1995, caf-shi-rus-1995, jin-par-tos-1998, jin-1999, jin-par-2001, klar-1998, gos-tos-2002, diaz-chen-2015, fil-jin-2011, fil-jin-2010, jin-qin-2013, hu-wan-2015, bos-cho-rus-yun-2021}. AP schemes maintain numerical stability across highly collisional regimes without requiring spatial or temporal step sizes to resolve the microscopic scales. They guarantee that the appropriate macroscopic continuum equations are recovered at the discrete level, allowing the time step to be selected independently of the local Knudsen number \cite{shi-2010}.

    The main objective of this paper is the design of a spatially adaptive, automated numerical method for gas mixtures operating within a single-temperature, single-velocity hydrodynamic limit. At the microscopic level, this physical regime is characterized by intra- and inter-species collisions acting on the same asymptotic time scale. Thanks to these interactions, the distinct components of the mixture rapidly converge toward a common macroscopic velocity and temperature. Since all species of the mixture converge toward a single local equilibrium state, our domain decomposition algorithm treats the whole mixture as a single entity (i.e. all species switch together between the macroscopic solver and the kinetic solver).

    From a mathematical point of view, the main new contribution of this work is the derivation of a moment realizability matrix for this single-velocity and single-temperature limit. This matrix is used to detect the regions out of local equilibrium, and therefore to decompose the domain. In this paper the criterion is derived and tested for a multi-species BGK operator. The collision model enters the indicator only through the Chapman--Enskog corrections, so the same construction can be used for the multi-species Boltzmann operator once these corrections are known. The coupling between the two solvers does not depend on the collision operator, so the proposed strategy can be extended to the full multi-species Boltzmann collision integral.

    To reduce the cost, the kinetic solver is used only in the cells selected by the indicator. The coupling is \emph{non-intrusive}. The fluid and kinetic solvers only exchange cell values near the interfaces, and they only need to use the same Runge--Kutta stages.

    The opposite regime, in which collisions between particles of the same species dominate, is considered in the companion paper \cite{PaperA}. There the hydrodynamic limit has a velocity and a temperature for each species, the indicator is evaluated species by species, and the method is second order on unstructured triangular meshes. Here all collisions occur on the same time scale, so the indicator is derived for the mixture as a whole and the same partition of the domain is used for all species. The kinetic model, the implicit relaxation step and the coupling between the solvers follow \cite{PaperA}, with the present scaling, and are briefly recalled to keep the paper self-contained. Besides the indicator, the main difference is the order of accuracy: both solvers use a third-order CWENO reconstruction and a third-order IMEX scheme \cite{bos-par-rus-2013, bos-par-rus-2024}, with a penalization of the collision term, on Cartesian meshes.

    The article is organized as follows. Section \ref{secModel} presents the kinetic model and its hydrodynamic limit. In Section \ref{secIndicators} we derive the realizability conditions and the switching criterion. Section \ref{secNumSch} describes the numerical schemes, in particular the AP IMEX time integration and the coupling between the solvers. Section \ref{secSimulations} presents the numerical tests, and Section \ref{secNumEff} discusses the computational cost. Conclusions and future work are given in Section \ref{secCon}.

    \section{Kinetic model for single-temperature, single-velocity gas mixtures}\label{secModel}

    Following \cite{bis-gro-mar2021, bob-bis-gro-spi-pot-2018}, we consider a mixture of $N$ monatomic gases with elastic collisions. The distribution $f_i(t,x,v)$ of species $i$, with $x,v\in\RR^3$ and $t\in\RR_+$, satisfies
    \begin{equation}
        \begin{cases}\label{eq1}
            \displaystyle\frac{\partial f_i}{\partial t} + v\cdot \nabla_x f_i = \mathcal{Q}_i[f],\\            
            \displaystyle\mathcal{Q}_i[f]=\sum_{k=1}^N\mathcal{Q}_{i,k}[f_i, f_k],\\
            \displaystyle f_i(0, x, v) = f_{i,0}(x, v).\\          
        \end{cases}   
    \end{equation}
    Here, the term $\mathcal{Q}_i[f]$ represents the total collisional source for the $i$-th species. The specific binary interaction integral $\Q_{i,k} [f_i,f_k]$ between species $i$ and $k$ is expressed as
    \begin{align}
        \label{eqQBoltzMulti}
        \Q_{i,k} [f_i,f_k] & \,=\, \int_{\RR^{3}} \int_{\mathbb{S}^{2}}  g_{i,k}(|y|,\hat{y}\cdot\Omega) \, \left[ f_i(v') f_k(w') \,-\, f_i(v) f_k(w) \right] \, d\Omega \, dw.
    \end{align}
    In this integral, $y = v - w$ denotes the relative velocity vector of the interacting particle pair $v$ and $w$, while $v'$ and $w'$ refer to their respective velocities prior to the collision. The intermolecular forces enter through the collision kernel $g_{i,k}$, which depends on the relative speed $\vert{}y\vert{}$ and on $\mu = \hat{y} \cdot \Omega$, with the unit vector defined as $\hat{y} = y / \vert{}y\vert{}$. 

    The number density $n_i(t,x)$, the mean velocity $u_i(t,x)$ and the temperature $T_i(t,x)$ of each species are
    \begin{equation}
        \begin{cases}
            n_i \displaystyle = \int_{\RR^3} f_i(v) dv,\\[10pt]
            u_i \displaystyle = \frac{1}{n_i}\int_{\RR^3}v f_i (v) dv,\\[10pt]
            T_i \displaystyle = \frac{m_i}{3n_i}\int_{\RR^3}|v-u_i|^2 f_i(v)dv.
        \end{cases}
    \end{equation}
    The equilibrium states of \eqref{eq1} are local Maxwellians with a common velocity $u$ and a common temperature $T$:
    \begin{equation}
        f_i^{eq}=n_i\left(\frac{m_i}{2\pi T}\right)^{\frac{3}{2}} \exp{\left(-\frac{m_i}{2T}|v-u|^2\right)},
    \end{equation}
    where the total number density $n$, the total mass density $\rho$, the velocity $u$ and the temperature $T$ of the mixture are
    \begin{equation}
        \begin{cases}
            n \displaystyle = \sum_{i=1}^N n_i,\\
            \rho \displaystyle = \sum_{i=1}^N \rho_i = \sum_{i=1}^N m_i n_i,\\
            u \displaystyle = \frac{1}{\rho}\sum_{i=1}^N m_i n_i u_i,\\
            T \displaystyle = \frac{1}{3n}\sum_{i=1}^N m_i\int_{\RR^3}|v-u|^2 f_i(v) dv =
                              \frac{1}{n}\sum_{i=1}^N n_i T_i + \frac{1}{3n}\sum_{i=1}^N m_i n_i |u_i - u|^2.
        \end{cases}
    \end{equation}
    As in \cite{bis-gro-mar2021, bob-bis-gro-spi-pot-2018}, we replace the Boltzmann collision operator with a BGK-type model, in which the collision term is a sum of relaxation terms:
    \begin{equation}
        \begin{cases}\label{eq8}
            \displaystyle\mathcal{Q}_i[f]=\sum_{k=1}^N\tilde{\mathcal{Q}}_{i,k}[f_i, f_k],\\
            \tilde{\mathcal{Q}}_{i,k}[f_i, f_k] = \nu_{ik}(n_i M_{ik} - f_i).
        \end{cases}        
    \end{equation}
    Here $\nu_{ik}$ is the collision frequency, independent of $v$, between species $i$ and $k$, and $M_{ik}$ is a Maxwellian with the mass $m_i$ of species $i$ and with auxiliary velocity $u_{ik}$ and temperature $T_{ik}$:
    \begin{equation}
        M_{ik} = \left(\frac{m_i}{2\pi T_{ik}}\right)^{\frac{3}{2}} \exp{\left(-\frac{m_i}{2T_{ik}}|v-u_{ik}|^2\right)}.
    \end{equation}
    These parameters are chosen so that the exchange rates of momentum and energy of the BGK operator coincide with those of the Boltzmann operator for Maxwellian interactions. In particular, the total momentum and the total energy of the mixture are conserved \cite{bis-gro-mar2021}. This gives \cite{bis-gro-mar2021}
    \begin{equation}\label{eqUT}
        \begin{cases}
            u_{ik} = (1 - a_{ik})u_i + a_{ik}u_k,\\
            T_{ik} = (1 - b_{ik})T_i + b_{ik}T_k + \gamma_{ik}|u_i - u_k|^2.
        \end{cases}
    \end{equation}
    The coefficients $a_{ik}$, $b_{ik}$ and $\gamma_{ik}$ are
    \begin{equation}
    \label{eq_aik_bik_gammaik}
        \begin{cases}
            a_{ik} \displaystyle = \frac{\lambda_{ik}n_k m_k}{\nu_{ik}(m_i + m_k)},\\[10pt]
            b_{ik} \displaystyle = \frac{2 a_{ik}m_i}{m_i + m_k},\\[10pt]
            \gamma_{ik} \displaystyle = \frac{m_i a_{ik}}{3}\left(\frac{2m_k}{m_i + m_k}-a_{ik}\right),
        \end{cases}
    \end{equation}
    where $\lambda_{ik}\displaystyle =2\pi \int_{-1}^{1}(1-\mu)g_{i,k}(\mu) d\mu$; for Maxwellian interactions $g_{i,k}$ does not depend on $|y|$ and $\lambda_{ik}$ is a constant.
    
    For $k=i$ one has $u_{ii}=u_i$ and $T_{ii}=T_i$, so the term describing collisions within species $i$ is the usual single-species BGK operator.

    We assume $\lambda_{ik}=\lambda_{ki}$, as required by pairwise conservation of momentum and energy. We consider inert mixtures in which collisions between particles of the same species and of different species occur at comparable rates. Denoting by $\varepsilon$ the Knudsen number, that is the ratio between the mean free path and a macroscopic reference length, the scaled multi-species BGK equations read \cite{bis-gro-mar2021}
    \begin{equation}
        \frac{\partial f_i}{\partial t}+v\cdot \nabla_x f_i = \frac{1}{\varepsilon}\sum_{k=1}^N\nu_{ik}(n_i M_{ik}-f_i), \qquad i=1, ..., N.
    \end{equation}
   The model differs from the one considered in \cite{PaperA} in the scaling: there only the collisions between particles of the same species are stiff, while here the collisions between different species are stiff as well, and the whole mixture relaxes to a single velocity and a single temperature. For $\varepsilon \ll 1$ collisions are frequent and the mixture behaves as a \emph{fluid}. For $\varepsilon \gtrsim 1$ collisions are rare and the gas is rarefied.

\section{Computation of indicators}\label{secIndicators}
In this section we adapt the indicators based on the moment realizability matrix to a gas mixture whose hydrodynamic limit is the single-temperature, single-velocity Euler system.
Moment realizability was introduced in \cite{lev-mor-nad-1998} to detect local departures from equilibrium in single-species gases, and was then extended by several authors, see for instance \cite{fil-rey-2014, cap-ten-2025, fil-xio-2018, tiw-2000, xio-qiu-2017, cap-par-rey-2026}.
It uses macroscopic quantities to estimate departure from equilibrium by comparing successive approximations of a moment matrix. Loss of positive semi-definiteness is a sufficient warning that the truncated moment description is not realizable. The resulting indicators are cheap, because they do not require the distribution function \cite{lev-mor-nad-1998, fil-rey-2014}. To our knowledge, this is the first application of realizability-based domain decomposition to gas mixtures with a single-temperature, single-velocity hydrodynamic limit.
In \cite{PaperA} a realizability matrix is built for each species, on its own velocity and temperature, because in that regime each species has its own Euler equations. Here the fluid model is the Euler system of the whole mixture, and a single matrix is built on the mixture velocity $u$ and temperature $T$. The indicators are built from the Chapman--Enskog expansion, which we briefly recall below.

\subsection{The Chapman--Enskog expansions}
Following \cite{bis-gro-mar2021}, we outline the Chapman--Enskog expansion suited for mixtures in which intra- and inter-species collisions occur on the same asymptotic time scale.
First, recall the BGK system \cite{bis-gro-mar2021}
\begin{equation}\label{rescaledBGK}
    \frac{\partial f_i}{\partial t}+v\cdot \nabla_x f_i = \frac{1}{\varepsilon}\sum_{k=1}^N\nu_{ik}(n_i M_{ik}-f_i), \qquad i=1, ..., N.
\end{equation}
The goal is to determine the asymptotic expansion $f_i = f_i^{(0)}+\varepsilon f_i^{(1)}$. In parallel, expansions are postulated for the number densities, flow velocities, temperatures, and Maxwellian attractors,
\[
n_i = n_i^{(0)}+\ve n_i^{(1)}, \quad 
u_i = u_i^{(0)}+\ve u_i^{(1)},
\quad
T_i = T_i^{(0)}+\ve T_i^{(1)},
\quad
M_{ik} = M_{ik}^{(0)}+\ve M_{ik}^{(1)},
\]
and auxiliary parameters
\[
\nu_{ik} = \nu_{ik}^{(0)}+\ve \nu_{ik}^{(1)}, \quad 
u_{ik} = u_{ik}^{(0)}+\ve u_{ik}^{(1)},
\quad
T_{ik} = T_{ik}^{(0)}+\ve T_{ik}^{(1)}.
\]
Under this formulation, requiring the first-order variations of $n_i$, $u$, and $T$ to vanish leads to the constraints \cite{bis-gro-mar2021}
\begin{equation}
    n_i^{(1)}=0, \quad \sum_{i=1}^N m_i n_i u_i^{(1)}=0, \quad \sum_{i=1}^N n_i T_i^{(1)}=0.
\end{equation}
In addition we have that \cite{bis-gro-mar2021}
\begin{equation}\label{eq4}
    f_i^{(0)}=f_i^{eq}=n_i\left(\frac{m_i}{2\pi T}\right)^{\frac{3}{2}} \exp{\left(-\frac{m_i}{2T}|v-u|^2\right)},
\end{equation}
and concerning the mixed Maxwellian $M_{ik}$ \cite{bis-gro-mar2021}
\begin{equation}\label{eq6}
    \begin{cases}
        u_{ik}^{(0)} \displaystyle = (1-a_{ik}^{(0)})u_i^{(0)} + a_{ik}^{(0)}u_k^{(0)},\\[10pt]
        u_{ik}^{(1)} \displaystyle = (1-a_{ik}^{(0)})u_i^{(1)} + a_{ik}^{(0)}u_k^{(1)} + a_{ik}^{(1)}(u_k^{(0)}-u_i^{(0)}),\\[10pt]
        T_{ik}^{(0)} \displaystyle = (1 - b_{ik}^{(0)})T_i^{(0)} + b_{ik}^{(0)}T_k^{(0)} + \gamma_{ik}^{(0)}|u_i^{(0)} - u_k^{(0)}|^2,\\[10pt]
        \begin{aligned}
        T_{ik}^{(1)} &= (1-b_{ik}^{(0)})T_i^{(1)}+b_{ik}^{(0)}T_k^{(1)}\\
        &\quad+2\gamma_{ik}^{(0)}(u_i^{(1)}-u_k^{(1)})\cdot(u_i^{(0)}-u_k^{(0)})\\
        &\quad+b_{ik}^{(1)}(T_k^{(0)}-T_i^{(0)})+\gamma_{ik}^{(1)}|u_i^{(0)}-u_k^{(0)}|^2,
        \end{aligned}
    \end{cases}
\end{equation}
where regarding the coefficients $a_{ik}, b_{ik},  \gamma_{ik}$ one gets \cite{bis-gro-mar2021}
\begin{equation}
    \begin{cases}    a_{ik}^{(0)}=\displaystyle\frac{\lambda_{ik}^{(0)}n_k^{(0)}m_k}{\nu_{ik}^{(0)}(m_i+m_k)},\\[10pt]
    a_{ik}^{(1)}=\displaystyle\frac{m_k}{\nu_{ik}^{(0)}(m_i+m_k)}\Bigg(\lambda_{ik}^{(1)}n_k^{(0)}+\lambda_{ik}^{(0)}n_k^{(1)}-\frac{\lambda_{ik}^{(0)}n_k^{(0)}\nu_{ik}^{(1)}}{\nu_{ik}^{(0)}}\Bigg),\\[10pt]
    b_{ik}^{(0)}=\displaystyle\frac{2a_{ik}^{(0)}m_i}{m_i+m_k}, \qquad b_{ik}^{(1)}=\displaystyle\frac{2a_{ik}^{(1)}m_i}{m_i+m_k}, \\[10pt]
    \gamma_{ik}^{(0)}=\displaystyle\frac{m_i a_{ik}^{(0)}}{3}\Bigg(\frac{2m_k}{m_i+m_k}-a_{ik}^{(0)}\Bigg), \qquad \gamma_{ik}^{(1)}=\displaystyle\frac{2}{3}a_{ik}^{(1)}\Bigg(m_{ik}-m_i a_{ik}^{(0)}\Bigg),
    \end{cases}
\end{equation}
where $m_{ik}=m_im_k/(m_i+m_k)$ is the reduced mass.

To obtain macroscopic system laws for $n_i, u, T$ at either Euler or Navier--Stokes scale, we consider the weak formulation \cite{bis-gro-mar2021} of kinetic model \eqref{rescaledBGK}:
\begin{equation} \label{weakForm}
\begin{aligned}
    &\frac{\partial}{\partial t} \int_{\mathbb{R}^3} \phi_i(v) \left( f_i^{(0)} + \varepsilon f_i^{(1)} \right) dv 
    + \nabla_x \cdot \int_{\mathbb{R}^3} \phi_i(v) v \left( f_i^{(0)} + \varepsilon f_i^{(1)} \right) dv \\[8pt]
    &\quad = \sum_{k=1}^N \left( \nu_{ik}^{(0)} + \varepsilon \nu_{ik}^{(1)} \right) 
    \int_{\mathbb{R}^3} \phi_i(v) \left( n_i M_{ik}^{(1)} - f_i^{(1)} \right) dv,
\end{aligned}
\end{equation}
where $\phi_i(v)=1, m_i v, \frac{1}{2} m_i |v|^2$. Truncating at zeroth order and summing the momentum and energy equations over the species produces the single-temperature, single-velocity Euler system. On the other hand, incorporating the $\mathcal{O}(\varepsilon)$ perturbation $f_i^{(1)}$ introduces dissipative terms (species diffusion, viscosity and thermal conduction), yielding the single-temperature, single-velocity Navier--Stokes equations.

\subsection{Euler equations}
At leading order, the weak formulation \eqref{weakForm} yields the Euler system for a mixture with common velocity and temperature \cite{bis-gro-mar2021}:
\begin{equation}
    \begin{cases}\label{eq7}
        \displaystyle \frac{\partial n_i}{\partial t}+\nabla_x\cdot(n_i u) = 0,\\[10pt]
        \displaystyle \frac{\partial}{\partial t}(\rho u) + \nabla_x\cdot(\rho u\otimes u)+\nabla_x(n T) \displaystyle=0,\\[10pt]
        \displaystyle \frac{\partial}{\partial t}\left(\frac{1}{2}\rho |u|^2 + \displaystyle\frac{3}{2}n T\right) + \nabla_x\cdot \left[\left(\frac{1}{2}\rho |u|^2 + \frac{5}{2}n T\right)u\right] \displaystyle=0.
    \end{cases}
\end{equation}

\subsection{Navier--Stokes equations}
To reach the Navier--Stokes hydrodynamic level via weak form \eqref{weakForm}, the $\mathcal{O}(\varepsilon)$ non-equilibrium contribution $f_i^{(1)}$ must be retained. For single-temperature, single-velocity mixture flows, this first-order perturbation takes the form \cite{bis-gro-mar2021}
\begin{equation}\label{eq100}
\begin{aligned}
    f_i^{(1)}\displaystyle &=\frac{f_i^{(0)}}{\sum_{k=1}^N \nu_{ik}^{(0)}}\Bigg\{ \frac{m_i}{T}c\cdot \sum_{k=1}^N\nu_{ik}^{(0)}u_{ik}^{(1)}   + \Bigg(\frac{m_i}{2}\frac{|c|^2}{T}-\frac{3}{2}\Bigg)\frac{1}{T}\sum_{k=1}^N \nu_{ik}^{(0)}T_{ik}^{(1)}\\
            &-m_ic\cdot \Bigg(\frac{\nabla_x n_i}{\rho_i}-\frac{\nabla_x n}{\rho}\Bigg) - m_i\Bigg(\frac{c}{T^{\frac{1}{2}}}\otimes \frac{c}{T^{\frac{1}{2}}} - \frac{1}{3}\frac{|c|^2}{T} I\Bigg):\nabla_x u\\
            &- \Bigg(\frac{m_i}{2}\frac{|c|^2}{T} - \frac{3}{2} - \frac{m_i n}{\rho}\Bigg)\frac{c}{T^{\frac{1}{2}}} \cdot \frac{\nabla_x T}{T^{\frac{1}{2}}} \Bigg\},
\end{aligned}
\end{equation}
where $c = v - u$ and $I$ is the identity matrix. We refer to \cite{bis-gro-mar2021} for the derivation of this expression and for the computation of the transport coefficients.

The auxiliary corrections $u_{ik}^{(1)}$ and $T_{ik}^{(1)}$ follow from \eqref{eq6} once the species corrections are known. Taking the first velocity moment of \eqref{eq100} gives
\begin{equation}\label{eq:CE_diffusion}
\sum_{k\ne i}\lambda_{ik}^{(0)}m_{ik}n_in_k
 \big(u_i^{(1)}-u_k^{(1)}\big)
 =-\nabla_x(n_iT)+\frac{\rho_i}{\rho}\nabla_x(nT),
\qquad i=1,\dots,N.
\end{equation}
Together with $\sum_i\rho_i u_i^{(1)}=0$, this system determines the diffusion velocities whenever the inter-species collision graph is connected. The temperature moments give
$\sum_k\nu_{ik}^{(0)}b_{ik}^{(0)}(T_i^{(1)}-T_k^{(1)})=0$.
For a connected collision graph, the constraint $\sum_i n_iT_i^{(1)}=0$ therefore yields $T_i^{(1)}=0$ for all $i$, and hence $T_{ik}^{(1)}=0$ in this monatomic single-temperature scaling.

Inserting the first-order expansion into weak form \eqref{weakForm} and keeping terms up to $\mathcal{O}(\varepsilon)$ (neglecting $\mathcal{O}(\varepsilon^2)$ contributions) incorporates non-equilibrium fluxes into the moment integrals, leading to the single-temperature, single-velocity Navier--Stokes system \cite{bis-gro-mar2021}:
\begin{equation}
\begin{cases}
    \displaystyle\frac{\partial n_i}{\partial t} + \nabla_x \cdot (n_i u) + \ve \nabla_x \cdot \Big(n_i u_i^{(1)}\Big) = 0 \\[10pt]
    \displaystyle\frac{\partial}{\partial t}(\rho u) + \nabla_x \cdot (\rho u \otimes u) + \nabla_x(n T) + \varepsilon \nabla_x \cdot P^{(1)} =0 \\[10pt]
    \displaystyle\frac{\partial}{\partial t}\left(\frac{1}{2}\rho |u|^2 + \frac{3}{2}n T\right) + \nabla_x \cdot \left[\left(\frac{1}{2}\rho |u|^2 + \frac{5}{2}n T\right)u\right] \\[10pt]
    \qquad + \varepsilon \nabla_x \cdot \left(P^{(1)} \cdot u\right) + \varepsilon \nabla_x \cdot q^{(1)} = 0,
\end{cases}
\end{equation}
where \cite{bis-gro-mar2021}
\begin{equation}
  \begin{gathered}
      P^{(1)}=\sum_{i=1}^Nm_i\int_{\RR^3}c\otimes c f_i^{(1)}(v)dv, 
      \qquad
      q^{(1)} = \frac{1}{2}\sum_{i=1}^N m_i\int_{\RR^3}c|c|^2f_i^{(1)}(v)dv,\\
      u_i^{(1)}=\frac{1}{n_i}\int_{\RR^3}v f_i^{(1)}dv.
  \end{gathered}
\end{equation}

\subsection{Regime indicators}
We now derive the regime indicators, extending the moment realizability approach of \cite{lev-mor-nad-1998, fil-rey-2014, fil-xio-2018, tiw-2000, li-son-wan-2021, xio-qiu-2017, cap-par-rey-2025, cap-ten-2025, cap-par-rey-2026, PaperA} to the present multi-species kinetic model. We first define the average particle mass of the mixture, $m = \rho/n$, and then we introduce the following polynomial elements:
\begin{equation}
    \begin{cases}\label{soninePol_scaled}
        \displaystyle {A}_i = \frac{m_i}{n} \Bigg[\frac{(v-u)\otimes (v-u)}{T} - \frac{1}{3}\frac{|v-u|^2}{T}I\Bigg],\\[12pt]
        \displaystyle {B}_i = \frac{\sqrt{m}}{n}\Bigg[\frac{m_i}{2T}|v-u|^2-\frac{3}{2}-\frac{m_i n}{\rho}\Bigg]\frac{v-u}{T^{\frac{1}{2}}},\\[12pt]
        \displaystyle {C}_i = \frac{m\, m_i}{n}\Bigg[\frac{1}{2T}|v-u|^2 - \frac{3}{2m_i} - \frac{n}{\rho}\Bigg]^2.
    \end{cases}
\end{equation}
We then introduce their macroscopic averages over the velocity space: 
\begin{equation}
    \begin{cases}
        \displaystyle{\bar{A}} = \sum_i\int {A}_i f_i dv, \\
        \displaystyle{\bar{B}} = \sum_i\int {B}_i f_i dv, \\
        \displaystyle{\bar{C}} = \sum_i\int {C}_i f_i dv.
    \end{cases}
\end{equation}

With these scaled quantities, the first-order perturbation $f_i^{(1)}$ in \eqref{eq100} can be recast as:
\begin{equation}\label{eq3_scaled}
\begin{aligned}
    f_i^{(1)}\displaystyle &=\frac{f_i^{(0)}}{\sum_{k=1}^N \nu_{ik}^{(0)}}\Bigg\{ \frac{m_i}{T}c\cdot \sum_{k=1}^N\nu_{ik}^{(0)}u_{ik}^{(1)}   + \Bigg(\frac{m_i}{2}\frac{|c|^2}{T}-\frac{3}{2}\Bigg)\frac{1}{T}\sum_{k=1}^N \nu_{ik}^{(0)}T_{ik}^{(1)}\\
            &-m_ic\cdot \Bigg(\frac{\nabla_x n_i}{\rho_i}-\frac{\nabla_x n}{\rho}\Bigg) - n {A}_i:\nabla_x u\\
            &- \frac{n}{\sqrt{m}} {B}_i \cdot \frac{\nabla_x T}{T^{\frac{1}{2}}} \Bigg\},
\end{aligned}
\end{equation}
where $c=v-u$.

In \cite{lev-mor-nad-1998} the realizability criterion indicates the transition between the Euler and the Navier--Stokes regimes of a single gas. Here we use it for the mixtures \eqref{rescaledBGK}, in which all collisions occur on the same time scale and the hydrodynamic limit has a single velocity and a single temperature. As for a single species \cite{fil-rey-2014, cap-par-rey-2025, tiw-2000, cap-par-rey-2026, cap-ten-2025}, the indicator is used to switch between the fluid and the kinetic solvers.

Non-negativity of $f_i$ implies $\int \psi^2 f_i dv \ge 0$ for any choice of velocity test function $\psi(v)$. Choosing $\psi = a^\intercal K_i$ with $a\in \RR^p$ and $K_i(v)$ a $p$-dimensional vector of polynomial basis functions gives
\begin{equation}
    a^\intercal \left(\int_{\RR^3} m_i K_i K_i^\intercal f_i\, dv\right) a = \int_{\RR^3} m_i (a^\intercal K_i)^2 f_i\, dv\geq 0, \quad a\in \RR^p.
\end{equation}
Thus, if $f_i$ is non-negative, then for any vector $K_i$, the symmetric matrix defined by $\int m_i K_i K_i^\intercal f_i dv$ is positive semi-definite. Evaluating this moment matrix at successive orders of the Chapman--Enskog expansion gives a simple check: if the matrix is not positive semi-definite, the truncated expansion does not correspond to a non-negative distribution, and the fluid closure is not valid.

We choose
\begin{equation}
    K_i = \begin{bmatrix}
        1 \\[10pt]
        \displaystyle \sqrt{m}\frac{v-u}{\sqrt{T}} \\[10pt]
        \displaystyle m\left(\frac{1}{2T}|v-u|^2 - \frac{3}{2m_i} - \frac{n}{\rho}\right)
    \end{bmatrix}.
\end{equation}
With this scaling, the realizability matrix below takes a simple form.
We define the global mixture realizability matrix 
\begin{equation}\label{eq:Wmix}
    \mathcal{W} = \frac{1}{\rho}\sum_{i=1}^N \int_{\RR^3} m_i K_i K_i^\intercal f_i \, dv.
\end{equation}
This choice follows the structure of the single-temperature, single-velocity limit and of its Chapman--Enskog expansion. Since $\mathcal{W}$ is the sum of the species matrices, if $\mathcal{W}$ is not positive semi-definite, then at least one of them is not.

By computing the integrals one gets:
\begin{equation}
    \mathcal{W} = 
    \begin{pmatrix}
        1  & 0^\intercal & -1\\[2mm]
        0    & I + {\bar{A}} & {\bar{B}}\\[2mm]
        -1 & {\bar{B}}^\intercal & {\bar{C}}
    \end{pmatrix}.
\end{equation}

We first derive a lower bound for ${\bar{C}}$. Substituting the definition of ${C}_i$:
\begin{equation}\label{eqTempC_scaled}
\begin{aligned}
    {\bar{C}} &= \frac{m}{n} \sum_i \int_{\RR^3} m_i \Bigg[\frac{1}{2T}|v-u|^2 - \frac{3}{2m_i} - \frac{n}{\rho}\Bigg]^2 f_i \, dv \\[2mm]
    &= \frac{m}{n} \left[ \frac{1}{4T^2} \sum_i m_i \int_{\RR^3} |v-u|^4 f_i \, dv - \frac{3}{2T} \sum_i \int_{\RR^3}|v-u|^2 f_i \, dv + \frac{9}{4} \sum_i \frac{n_i}{m_i} + \frac{n^2}{\rho} \right].
\end{aligned}
\end{equation}
A lower bound for the fourth-order moment $\int_{\RR^3} |v-u|^4 f_i \, dv$ can be estimated using Hölder's inequality 
\begin{equation}
 \label{eq::Holder}
    \int\left|g(v)h(v)\right|dv \le \left(\int|g(v)|^{p} dv\right)^{1/p} \left(\int|h(v)|^{q}dv\right)^{1/q},  
 \end{equation}
 
 with $p=q=2$, $h(v)=\sqrt{f_i}|v-u|^2$, and $g(v)=\sqrt{f_i}$. By the equality condition in H\"older's inequality \cite{mit-1970}, equality holds only if $|g|^2$ and $|h|^2$ are proportional almost everywhere, that is, only if $|v-u|$ is constant on the support of $f_i$. Under the non-degeneracy assumption that $f_i$ is not concentrated on a sphere $|v-u|=\mathrm{const}$, in both the continuous and the discrete setting, the inequality is therefore strict. Since $\int g^2\,dv=n_i$ and $\int gh\,dv=\int_{\RR^3}|v-u|^2f_i\,dv$, this gives
\begin{equation}\label{eq::bound4th}
    \int_{\RR^3} |v-u|^4 f_i \, dv > \frac{1}{n_i} \left( \int_{\RR^3} |v-u|^2 f_i \, dv \right)^2.
\end{equation}
Substituting inequality \eqref{eq::bound4th} into Equation \eqref{eqTempC_scaled} and completing the square for each species $i$ yields:
\begin{equation}
    {\bar{C}} > \frac{m}{n} \left[ \sum_i \frac{\left( m_i \int_{\RR^3} |v-u|^2 f_i \, dv - 3 n_i T \right)^2}{4 m_i n_i T^2} + \frac{n^2}{\rho} \right] \ge \frac{m}{n} \left(\frac{n^2}{\rho}\right).
\end{equation}
Since $m = \rho/n$, we obtain
\begin{equation}
    \bar{C} > \frac{m}{n} \left(\frac{n^2}{\rho}\right) = \frac{\rho}{n^2} \left(\frac{n^2}{\rho}\right) = 1.
\end{equation}

Since $\bar{C} > 1 > 0$, $\mathcal{W}$ is positive semi-definite if and only if its Schur complement $\mathcal{R}$ with respect to $\bar{C}$ \cite{horn2012} is positive semi-definite. We have
\begin{equation}
    {\mathcal{R}}=
    \begin{pmatrix}
        \displaystyle 1 - \frac{1}{{\bar{C}}}    & \displaystyle \frac{{\bar{B}}^\intercal}{{\bar{C}}} \\[4mm]
        \displaystyle \frac{{\bar{B}}}{{\bar{C}}}  & \displaystyle I + {\bar{A}} - \frac{{\bar{B}} {\bar{B}}^\intercal}{{\bar{C}}}
    \end{pmatrix}.
\end{equation}
Since $1 - 1/\bar{C} > 0$, we take the Schur complement $\mathcal{S}$ of $\mathcal{R}$ with respect to its first entry and, after simplification, we obtain
\begin{equation}\label{eq5}
    {\mathcal{S}} = I + {\bar{A}} - \frac{{\bar{B}} {\bar{B}}^\intercal}{{\bar{C}}}\Bigg[\frac{1}{1 - \frac{1}{{\bar{C}}}}\Bigg] = I + {\bar{A}} - \frac{{\bar{B}} {\bar{B}}^\intercal}{{\bar{C}} - 1}.
\end{equation}
Hence $\mathcal{W}$ is positive semi-definite if and only if $\mathcal{S}$ is.

The tensors ${\bar{A}}$, ${\bar{B}}$, and ${\bar{C}}$ depend on the underlying distribution function $f_i$. Evaluating them at zeroth order ($f_i = f_i^{(0)}$) and first order ($f_i = f_i^{(0)} + \varepsilon f_i^{(1)}$) yields the following explicit expressions.
At zeroth order one gets:
\begin{equation}
    \begin{cases}
        {\bar{A}}^{(0)} \displaystyle= 0,\\[4pt]
        {\bar{B}}^{(0)} = 0,\\[4pt]
        {\bar{C}}^{(0)} =  \displaystyle 1 + \frac{3}{2}  \frac{m}{n} \sum_{i=1}^N \frac{n_i}{m_i} 
    \end{cases}
\end{equation}
and at the first order:
\begin{equation}
    \begin{cases}
        {\bar{A}}^{(1)} \displaystyle = - \frac{\varepsilon}{n} \left( \sum_{i=1}^N \frac{ n_i}{\sum_{k=1}^N \nu_{ik}^{(0)}} \right) D(u),\\[16pt]
        {\bar{B}}^{(1)} = \varepsilon \displaystyle \frac{\sqrt{m}}{n} \sum_{i=1}^N \frac{n_i}{\left({\sum_{k=1}^N \nu_{ik}^{(0)}}\right) T^{\frac{1}{2}}} \left( 1 - \frac{m_i n}{\rho} \right) \left[ \sum_{k=1}^N \nu_{ik}^{(0)} u_{ik}^{(1)} - T \left( \frac{\nabla_x n_i}{\rho_i} - \frac{\nabla_x n}{\rho} \right) \right] \\[12pt]
        \qquad \quad - \, \varepsilon \displaystyle \frac{\sqrt{m}}{n} \left[ \sum_{i=1}^N \frac{1}{\left({\sum_{k=1}^N \nu_{ik}^{(0)}}\right) T^{\frac{1}{2}}} \left( \frac{7 n_i}{2 m_i} - \frac{2 n_i n}{\rho} + \frac{\rho_i n^2}{\rho^2} \right) \right] \nabla_x T,\\[16pt]
        {\bar{C}}^{(1)} \displaystyle = {\bar{C}}^{(0)} + \varepsilon \frac{m}{n} \sum_{i=1}^N \frac{3 n_i}{\left({\sum_{k=1}^N \nu_{ik}^{(0)}}\right) T} \left( \frac{1}{m_i} - \frac{n}{\rho} \right) \sum_{k=1}^N \nu_{ik}^{(0)}T_{ik}^{(1)},
    \end{cases}
\end{equation}
where $D(u) = \nabla_x u + (\nabla_x u)^\intercal - \frac{2}{3} (\nabla_x \cdot u)I$. Substituting these into \eqref{eq5} determines ${\mathcal{S}}^{(0)}$ and ${\mathcal{S}}^{(1)}$. In particular, since $\bar A^{(0)}=0$ and $\bar B^{(0)}=0$, we have $\mathcal{S}^{(0)}=I$.

The aggregate matrix monitors the mixture moment corrections; it does not separately bound every species diffusion flux. In particular, compensating species contributions may be invisible to a mixture-level indicator. For the monatomic scaling considered here, $T_{ik}^{(1)}=0$ and therefore $\bar C^{(1)}=\bar C^{(0)}>1$, so the denominator in the Schur-complement formula remains positive.
If $\mathcal{S}^{(0)}$ and $\mathcal{S}^{(1)}$ are close, the zeroth-order approximation is accurate, while large differences between their eigenvalues indicate that higher-order non-equilibrium effects are important. We therefore consider the Euler description (Layer 0) valid as long as ${\mathcal{S}}^{(1)}$ stays positive definite and its eigenvalues remain close to those of ${\mathcal{S}}^{(0)}$. 

The switching rule from the Euler solver (Layer 0) to the kinetic BGK solver (Layer 1) is activated when:
\begin{equation}
  \label{crit:CompEuler}
    \exists \, j \in \{1, 2, 3\} \quad \text{such that} \quad \left| \theta_j^1 - \theta_j^0 \right| > \eta,
\end{equation}
where $\theta_j^0 \in {\rm Sp}({\mathcal{S}}^{(0)})$ and $\theta_j^1 \in {\rm Sp}({\mathcal{S}}^{(1)})$ denote the corresponding eigenvalues (ordered consistently), and $\eta > 0$ is a small, single threshold parameter for the mixture. When this condition holds, the first-order correction is too large with respect to the zeroth-order description, and the Euler equations (Layer 0) are replaced by the kinetic equations (Layer 1).

\begin{rem}
    The reverse transition from Layer 1 (kinetic solver) back to Layer 0 (Euler solver) is governed by the same spectral indicator. A grid cell reverts to Layer 0 once the local spectrum relaxes toward equilibrium:
    \begin{equation}
    |\theta_j^1(t, x) - \theta_j^0(t, x)| \le \eta \qquad \forall \, j \in \{1, 2, 3\}.
    \end{equation}
\end{rem}

\begin{rem}\label{rem:hierarchy}
The indicator is based on the Chapman--Enskog expansion, and the realizability matrices can be computed at any truncation order. The method is therefore not limited to two levels: it can be used to build a hierarchy of models including, for instance, Navier--Stokes or Burnett-type closures, as in \cite{cap-par-rey-2025, cap-par-rey-2026} for a single gas. In such a hierarchy, the passage from level $l$ to level $l+1$ is decided as in \eqref{crit:CompEuler}, by comparing the eigenvalues of the matrices $\mathcal S^{(l)}$ and $\mathcal S^{(l+1)}$ obtained at two consecutive orders of the expansion, with a threshold for each level.
\end{rem}

The spectral indicator measures the Chapman--Enskog corrections represented in these moments; it is not a bound on the distance of an arbitrary kinetic distribution from a Maxwellian. In particular, spatially homogeneous non-equilibrium modes with equilibrium hydrodynamic moments need not be detected. The switching strategy is therefore assessed here for the initial data and flow configurations of the numerical tests.

When the initial data are not at equilibrium, a simple remedy is to compute the first time steps with the kinetic model everywhere (see also \cite{PaperA}): if the distribution functions are out of equilibrium, their evolution can produce non-equilibrium effects in the macroscopic quantities used by the indicator.

\section{Numerical schemes}
\label{secNumSch}
This section describes the discretization in velocity, space and time, and the correction that makes the discrete Maxwellians conservative.
\subsection{Velocity discretization}
We first discretize the velocity space of model \eqref{rescaledBGK}. The velocity domain is truncated to a bounded region. Since the distribution functions do not have compact support, the truncation breaks the conservation of mass, momentum and energy. To recover it, we use the technique of \cite{gam-tha-2009}, described in Section \ref{secConMas}.

The bounded velocity domain $[-v_{\max},v_{\max}]^{d_v}$ is discretized by a uniform Cartesian grid with $N_v$ points per direction and step $\Delta v=2v_{\max}/(N_v-1)$; we denote by $v_m$, $m=1,\dots,N_V$, with $N_V=N_v^{d_v}$, the points of the grid.
Each distribution function $f_i$ is then replaced by the vector $F_i(t, x) \in \mathbb{R}^{N_V}$ of its approximate values at the nodes $v_m$:
\begin{equation}
    F_i(t, x) = \big(f_{i,m}(t, x)\big)_{m=1}^{N_V}, \quad \text{with} \quad f_{i,m}(t, x) \approx f_i(t, x, v_m).
\end{equation}

As a result, the continuous kinetic system transforms into $N_V$ coupled hyperbolic transport equations with stiff relaxation terms for each species $i = 1, \dots, N$. The resulting discrete-velocity evolution equations take the form:
\begin{equation}
    \partial_t f_{i,m} + v_m \cdot \nabla_{x} f_{i,m} = \frac{1}{\varepsilon}\sum_{k=1}^N\nu_{ik}\big(n_i M_{ik,m} - f_{i,m}\big), \quad m=1, \dots, N_V,
\end{equation}
where $M_{ik,m}$ denotes the Maxwellian $M_{ik}$ evaluated at node $v_m$.

\subsection{Spatial discretization}
    The spatial discretization is based on a finite volume framework combined with a third-order Central Weighted Essentially Non-Oscillatory (CWENO3) reconstruction \cite{lev-pup-rus-2000}. Specifically, we consider a two-dimensional uniform Cartesian mesh that partitions the computational domain into $N_x \times N_y$ rectangular cells $I_{j,l} = [x_{j-1/2}, x_{j+1/2}] \times [y_{l-1/2}, y_{l+1/2}]$, for $j=1, \dots, N_x$ and $l=1, \dots, N_y$. We denote by $(x_j, y_l)$ the coordinates of the cell centre, with $\Delta x$ and $\Delta y$ representing the constant spatial grid spacings in the $x$ and $y$ directions, respectively. 

    Within this framework, the cell average of each component of the vectors representing the distribution functions of species $i$ at time $t^n$ is defined as:
    \begin{equation}
        \thickbar{f}_{i,m,j,l}^n = \frac{1}{\Delta x \Delta y} \int_{I_{j,l}} f_{i,m}(t^n, x) \, dx, 
    \end{equation}
    for all $i \in \{1, \dots, N\}$, $m \in \{1, \dots, N_V\}$, $j \in \{1, \dots, N_x\}$ and $l \in \{1, \dots, N_y\}$. 

    Similarly, the cell averages for the macroscopic quantities of the mixture, which are evolved using the single-velocity, single-temperature Euler equations \eqref{eq7}, are given by:
    \begin{equation}
        \thickbar{U}_{j,l}^n = \frac{1}{\Delta x \Delta y} \int_{I_{j,l}} U(t^n, x) \, dx,
    \end{equation}
    where $U(t, x) = [n_1, \dots, n_N, \rho u, E]^\intercal(t, x)$ denotes the continuous state vector of the macroscopic conservative variables for the entire mixture, with $E = \frac{1}{2}\rho |u|^2 + \frac{3}{2}n T$ being the total energy. 

    The numerical fluxes are computed from high-order reconstructions of the distribution functions and of the macroscopic state at the cell interfaces. To avoid spurious oscillations near discontinuities, we use a dimension-by-dimension CWENO3 reconstruction \cite{lev-pup-rus-2000} at suitable quadrature points, which gives the left and right states $f^{L/R}_{i,m}$ and $U^{L/R}$ at each cell edge. The numerical fluxes are then computed with the local Lax--Friedrichs flux \cite{pup-sem-vis-2022}. In the following, we also denote the cells $I_{j,l}$ by $P_c$, $c=1,\dots,N_P$, where $N_P$ is the total number of cells.

\subsection{IMEX asymptotic-preserving scheme}
\label{sec:numerical_properties}

In this section we study the IMEX scheme for the multi-species BGK model \eqref{rescaledBGK}. We first derive the linear systems for the velocities and temperatures of the species in the implicit step, and we show that they have a unique solution with positive temperatures. Then we prove the Asymptotic-Preserving (AP) property: as $\varepsilon \to 0$ the species relax to a common velocity and temperature, and the scheme gives a consistent discretization of the Euler equations \eqref{eq7}.
The analysis of the relaxation step follows the companion paper \cite{PaperA}, and the proofs are reported here for completeness.

\subsubsection{Treatment of the collision parameters}
\label{subsec:linearization}

In a general kinetic framework, the collision frequencies $\nu_{ik}$ may depend on the macroscopic moments. If these dependencies were treated fully implicitly at time $t^{n+1}$, the macroscopic moment systems would become strongly non-linear, requiring computationally expensive iterative solvers.

To preserve the efficiency of the IMEX framework while avoiding heavy notation, all kinetic collision frequencies in what follows are understood to be frozen at the known explicit state. Furthermore, because the density equations decouple and are trivially solved first (yielding $n_k^{n+1} = n_k^*$, as shown below), the interaction parameters can be evaluated using these newly updated densities alongside the frozen kinetic coefficients, using Equation \eqref{eq_aik_bik_gammaik}. This treatment ensures that the macroscopic relaxation steps for velocities and temperatures reduce to linear algebraic systems.

\subsubsection{Macroscopic moment systems}
\label{subsec:solvability}

We consider the first-order IMEX temporal update for a generic species $i \in \{1, \dots, N\}$, applied to the fully rescaled model \eqref{rescaledBGK}:
\begin{equation}
    f_i^{n+1} = f_i^* + \frac{\Delta t}{\varepsilon} \sum_{k=1}^N \nu_{ik} \left( n_i^{n+1} M_{ik}^{n+1} - f_i^{n+1} \right),
    \label{eq:imex_update_formal}
\end{equation}
where $f_i^* = f_i^n - \Delta t \, v \cdot \nabla_x f_i^n$ is the explicit advected distribution. 
By integrating Eq.~\eqref{eq:imex_update_formal} against $\phi_i(v) = \left(1, \, m_i v, \, \frac{1}{2} m_i |v|^2\right)^\intercal$, we derive the macroscopic equations governing the relaxation step. 

\paragraph{\bf Densities.}
Integrating with weight $\phi_i(v) = 1$, we see that the number density of each species is conserved:
\begin{equation}
    n_i^{n+1} = n_i^* + \frac{\Delta t}{\varepsilon} \sum_{k=1}^N \nu_{ik} (n_i^{n+1} - n_i^{n+1}) = n_i^*, \quad \forall i \in \{1, \dots, N\}.
\end{equation}
This fixes the densities and therefore, by \eqref{eq_aik_bik_gammaik}, the parameters $a_{ik}^{n+1}$, $b_{ik}^{n+1}$ and $\gamma_{ik}^{n+1}$ used in the next steps.

\paragraph{\bf Velocities.}
Integrating against $\phi_i(v) = m_i v$, we obtain the momentum balance:
\begin{equation}
    \rho_i^{n+1} u_i^{n+1} = \rho_i^* u_i^* + \frac{\Delta t}{\varepsilon} \sum_{k \neq i} \nu_{ik} \rho_i^{n+1} \left( u_{ik}^{n+1} - u_i^{n+1} \right).
\end{equation}
Substituting the target velocity $u_{ik}^{n+1} = (1 - a_{ik}^{n+1}) u_i^{n+1} + a_{ik}^{n+1} u_k^{n+1}$ and dividing by $\rho_i^{n+1}$ (since $\rho_i^{n+1} = \rho_i^*$), we collect the unknown velocities on the left-hand side to obtain an $N \times N$ linear system for each spatial component:
\begin{equation}
    \left( 1 + \frac{\Delta t}{\varepsilon} \sum_{k \neq i} \nu_{ik} a_{ik}^{n+1} \right) u_i^{n+1} - \frac{\Delta t}{\varepsilon} \sum_{k \neq i} \nu_{ik} a_{ik}^{n+1} u_k^{n+1} = u_i^*, \quad \forall i \in \{1, \dots, N\}.
    \label{eq:ns_velocity_system}
\end{equation}

\paragraph{\bf Temperatures.}
Integrating against $\phi_i(v) = \frac{1}{2} m_i |v|^2$ yields the energy balance. By decomposing the total energies into kinetic ($\mathcal{K}_i \equiv \frac{1}{2} \rho_i |u_i|^2$, and similarly $\mathcal{K}_{ik} \equiv \frac{1}{2} \rho_i |u_{ik}|^2$) and internal thermal components $E_i^{n+1} = \mathcal{K}_i^{n+1} + \frac{3}{2} n_i^{n+1} T_i^{n+1}$, we observe that the velocities $u_i^{n+1}$ are already known from Eq.~\eqref{eq:ns_velocity_system}, allowing all kinetic contributions to act as explicit sources.

Substituting $T_{ik}^{n+1} = (1 - b_{ik}^{n+1}) T_i^{n+1} + b_{ik} ^{n+1}T_k^{n+1} + \gamma_{ik}^{n+1} |u_i^{n+1} - u_k^{n+1}|^2$ and isolating the unknown temperatures, we obtain the $N \times N$ linear system:
\begin{equation}
    \left( 1 + \frac{\Delta t}{\varepsilon} \sum_{k \neq i} \nu_{ik} b_{ik}^{n+1} \right) T_i^{n+1} - \frac{\Delta t}{\varepsilon} \sum_{k \neq i} \nu_{ik} b_{ik}^{n+1} T_k^{n+1} = r_i^*, \quad \forall i \in \{1, \dots, N\},
    \label{eq:ns_temperature_system}
\end{equation}
with the explicitly known right-hand side source term:
\begin{equation}
    r_i^* = \frac{2}{3 n_i^{n+1}} \left[ E_i^* - \mathcal{K}_i^{n+1} + \frac{\Delta t}{\varepsilon} \sum_{k \neq i} \nu_{ik} \left( \mathcal{K}_{ik}^{n+1} - \mathcal{K}_i^{n+1} + \frac{3}{2} n_i^{n+1} \gamma_{ik}^{n+1} |u_i^{n+1} - u_k^{n+1}|^2 \right) \right].
\end{equation}

\subsubsection{Well-posedness of the relaxation systems}
Both systems \eqref{eq:ns_velocity_system} and \eqref{eq:ns_temperature_system} share the same mathematical structure $\mathcal{A} x = y$. The coefficient matrix $\mathcal{A} \in \mathbb{R}^{N \times N}$ has the generic form:
\begin{equation}
    \mathcal{A}_{ii} = 1 + \sum_{k \neq i} c_{ik}, \qquad \mathcal{A}_{ik} = -c_{ik} \text{ for } k \neq i,
\end{equation}
where $c_{ik}$ denotes either $\frac{\Delta t}{\varepsilon} \nu_{ik} a_{ik}^{n+1}$ or $\frac{\Delta t}{\varepsilon} \nu_{ik} b_{ik}^{n+1}$.

We first give bounds on the parameters of the model, which are needed for the positivity of the scheme.

\begin{lemma}[Admissibility of the parameters]
\label{lem:parameters}
Let the collision kernel $g_{i,k}(\mu)$ be non-negative, so that $\lambda_{ik} = 2\pi \int_{-1}^{1}(1-\mu)g_{i,k}(\mu)\, d\mu \ge 0$. Assume that the collision frequencies are positive and satisfy the lower bound
\begin{equation}\label{eq:nu_bound}
\nu_{ik} \ge \frac{\lambda_{ik} n_k^{n+1}}{2}.
\end{equation}
Then, for every $n_k^{n+1} \ge 0$, the mixture parameters defined in \eqref{eq_aik_bik_gammaik} satisfy:
\begin{equation}
    a_{ik}^{n+1} \ge 0, \qquad \gamma_{ik}^{n+1} \ge 0, \quad \text{and} \quad 0 \le b_{ik}^{n+1} \le 1.
\end{equation}
\end{lemma}

\begin{proof}
Since $g_{i,k}(\mu) \ge 0$ and $1-\mu \ge 0$ for all $\mu \in [-1, 1]$, we have $\lambda_{ik} \ge 0$. Given that $\nu_{ik} > 0$ and $n_k^{n+1} \ge 0$, the definition $a_{ik}^{n+1} = \frac{\lambda_{ik} n_k^{n+1} m_k}{\nu_{ik}(m_i + m_k)}$ directly yields $a_{ik}^{n+1} \ge 0$. Consequently, $b_{ik}^{n+1} = \frac{2 m_i}{m_i + m_k} a_{ik}^{n+1} \ge 0$.

Regarding the parameter $\gamma_{ik}^{n+1}$, from definition \eqref{eq_aik_bik_gammaik}, the condition $\gamma_{ik}^{n+1} \ge 0$ holds if and only if:
\[
    a_{ik}^{n+1} \le \frac{2 m_k}{m_i + m_k}.
\]
Substituting the explicit expression for $a_{ik}^{n+1}$, we obtain:
\[
    \frac{\lambda_{ik} n_k^{n+1} m_k}{\nu_{ik}(m_i + m_k)} \le \frac{2 m_k}{m_i + m_k} \iff \frac{\lambda_{ik} n_k^{n+1}}{\nu_{ik}} \le 2 \iff \nu_{ik} \ge \frac{\lambda_{ik} n_k^{n+1}}{2},
\]
which is the lower bound assumed on $\nu_{ik}$.

Finally, substituting the upper bound $a_{ik}^{n+1} \le \frac{2 m_k}{m_i + m_k}$ into the definition of $b_{ik}^{n+1}$ yields:
\[
    b_{ik}^{n+1} \le \frac{2 m_i}{m_i + m_k} \left( \frac{2 m_k}{m_i + m_k} \right) = \frac{4 m_i m_k}{(m_i + m_k)^2}\le 1.
\]
\end{proof}

\begin{theorem}[Unique solvability and positivity]\label{teo1}
Under the assumptions of Lemma~\ref{lem:parameters}, assume that $n_i^*>0$ and $T_i^*>0$ for all $i$. Then, for any $\Delta t > 0$ and any $\varepsilon > 0$, the $N \times N$ linear systems for velocities \eqref{eq:ns_velocity_system} and temperatures \eqref{eq:ns_temperature_system} admit a unique solution. Furthermore, the updated temperatures are strictly positive ($T_i^{n+1} > 0$ for all $i = 1, \dots, N$).
\end{theorem}

\begin{proof}
\emph{Step 1 (sign of the coefficients).} Since $n_i^{n+1}=n_i^*>0$, by Lemma~\ref{lem:parameters} the coefficients $a_{ik}^{n+1}$ and $b_{ik}^{n+1}$ are non-negative, hence $c_{ik}\ge0$ for all $i\ne k$.

\emph{Step 2 (invertibility).}
For every row $i$, the sum of the absolute values of the off-diagonal elements is $\sum_{k \neq i} |\mathcal{A}_{ik}| = \sum_{k \neq i} c_{ik}$. The absolute diagonal element is $|\mathcal{A}_{ii}| = 1 + \sum_{k \neq i} c_{ik}$. Therefore, $|\mathcal{A}_{ii}| > \sum_{k \neq i} |\mathcal{A}_{ik}|$ for all $i$, meaning $\mathcal{A}$ is strictly diagonally dominant. By the L\'evy--Desplanques theorem \cite{tyr-1997}, $\mathcal{A}$ is non-singular and invertible for any $\Delta t > 0$ and $\varepsilon > 0$.
Hence \eqref{eq:ns_velocity_system} has a unique solution $u^{n+1}$ for each spatial component; with $u^{n+1}$ known, the right-hand side $r_i^*$ is determined, and \eqref{eq:ns_temperature_system} has a unique solution as well.

\emph{Step 3 (lower bound for the temperatures).}
Let $T_m^{n+1} = \min_i \{ T_i^{n+1} \}$ be the minimum updated temperature. Evaluating the $m$-th row of system \eqref{eq:ns_temperature_system} and isolating $T_m^{n+1}$ yields:
\[
    T_m^{n+1} = r_m^* + \sum_{k \neq m} c_{mk} \left( T_k^{n+1} - T_m^{n+1} \right).
\]
Since $T_k^{n+1} - T_m^{n+1} \ge 0$ and $c_{mk} \geq 0$, we obtain
\[
    T_i^{n+1} \ge T_m^{n+1} \ge r_m^*, \quad \forall i \in \{1, \dots, N\}.
\]
To prove strict positivity ($T_i^{n+1} > 0$), we must now show that the source term $r_i^*$ is strictly positive.

\emph{Step 4 (positivity of $r_i^*$).}
We analyze the macroscopic source term from Eq.~\eqref{eq:ns_temperature_system}:
\[
r_i^* = \frac{2}{3 n_i^{n+1}} \left[ E_i^* - \mathcal{K}_i^{n+1} + \frac{\Delta t}{\varepsilon} \sum_{k \neq i} \nu_{ik} \left( \mathcal{K}_{ik}^{n+1} - \mathcal{K}_i^{n+1} + \frac{3}{2} n_i^{n+1} \gamma_{ik}^{n+1} |u_i^{n+1} - u_k^{n+1}|^2 \right) \right].
\]
First, we decompose the explicit total energy: $E_i^* = \frac{3}{2} n_i^* T_i^* + \mathcal{K}_i^*$. Let $V_i = \frac{\Delta t}{\varepsilon} \sum_{k \neq i} \nu_{ik} \rho_i^{n+1} (u_{ik}^{n+1} - u_i^{n+1})$. Squaring the implicit velocity update \eqref{eq:ns_velocity_system} yields:
\[
\mathcal{K}_i^* = \mathcal{K}_i^{n+1} - u_i^{n+1} \cdot V_i + \frac{|V_i|^2}{2 \rho_i^{n+1}}.
\]
Expanding the kinetic energy difference $\mathcal{K}_{ik}^{n+1} - \mathcal{K}_i^{n+1}$ using $\frac{1}{2}(|x|^2 - |y|^2) = (x-y)\cdot y + \frac{1}{2}|x-y|^2$, we get:
\[
\frac{\Delta t}{\varepsilon} \sum_{k \neq i} \nu_{ik} \left( \mathcal{K}_{ik}^{n+1} - \mathcal{K}_i^{n+1} \right) = u_i^{n+1} \cdot V_i + \frac{\Delta t}{\varepsilon} \sum_{k \neq i} \nu_{ik} \frac{\rho_i^{n+1}}{2} |u_{ik}^{n+1} - u_i^{n+1}|^2.
\]
Substituting these relations back into $r_i^*$, the terms $u_i^{n+1} \cdot V_i$ cancel and we obtain:
\[
r_i^* = T_i^* + \frac{2}{3 n_i^{n+1}} \left[ \frac{|V_i|^2}{2 \rho_i^{n+1}} + \frac{\Delta t}{\varepsilon} \sum_{k \neq i} \nu_{ik} \left( \frac{\rho_i^{n+1}}{2} |u_{ik}^{n+1} - u_i^{n+1}|^2 + \frac{3}{2} n_i^{n+1} \gamma_{ik}^{n+1} |u_i^{n+1} - u_k^{n+1}|^2 \right) \right].
\]
By Lemma \ref{lem:parameters}, we know that $\gamma_{ik}^{n+1} \ge 0$. Since $\nu_{ik} > 0$ and $\rho_i^{n+1} > 0$, every term inside the square brackets is non-negative. Since $T_i^* > 0$, it follows that $r_i^* > 0$, concluding the proof.
\end{proof}

\begin{rem}[Target Maxwellians]
The pointwise update \eqref{eq:final_pointwise_update} requires the target Maxwellians $M_{ik}^{n+1}$ to be well defined, that is their temperatures to be positive.
The target temperature is defined as:
\[
    T_{ik}^{n+1} = (1 - b_{ik}^{n+1}) T_i^{n+1} + b_{ik}^{n+1} T_k^{n+1} + \gamma_{ik}^{n+1} |u_i^{n+1} - u_k^{n+1}|^2.
\]
Theorem \ref{teo1} proves that the macroscopic temperatures are strictly positive ($T_i^{n+1}, T_k^{n+1} > 0$). Furthermore, Lemma \ref{lem:parameters} establishes that $0 \le b_{ik}^{n+1} \le 1$, so the first two terms form a convex combination of positive values. Since also $\gamma_{ik}^{n+1} \ge 0$ by Lemma \ref{lem:parameters}, we get $T_{ik}^{n+1} > 0$, and all the target Maxwellians are well defined.
\end{rem}

Once the full macroscopic state is determined, the distribution functions are updated point-by-point:
\begin{equation}
    f_i^{n+1} = \frac{f_i^* + \dfrac{\Delta t}{\varepsilon} \sum_{k=1}^N \nu_{ik} \, n_i^{n+1} M_{ik}^{n+1}}{1 + \dfrac{\Delta t}{\varepsilon} \sum_{k=1}^N \nu_{ik}}.
    \label{eq:final_pointwise_update}
\end{equation}

The positivity result concerns the local implicit relaxation step with positive input densities and temperatures. It does not by itself prove positivity of the preceding high-order transport and Runge--Kutta stages.
\subsection{Asymptotic-preserving property}
\label{subsec:AP_Euler}

As $\varepsilon \to 0$, the collisions between different species drive all species of model \eqref{rescaledBGK} to a common velocity $u$ and a common temperature $T$.

\begin{theorem}[Single-fluid Euler limit, first order] 
Assume positive species densities and temperatures, $\lambda_{ik}=\lambda_{ki}>0$ for $i\ne k$, and well-prepared initial data. In the limit $\varepsilon \to 0$, the solutions of the $N$-species moment systems satisfy $u_i \to u$ and $T_i \to T$, and the discrete distribution functions $f_i^{n+1}$ relax to the local Maxwellians $n_i M_i(u, T)$. The scheme then gives a consistent discretization of the single-velocity, single-temperature Euler equations \eqref{eq7} (see also \cite{bis-bos-dim-gro-mar-2022}).
\end{theorem}

\begin{proof}
\emph{Step 1 (synchronization of the macroscopic fields).}
Multiplying the velocity system \eqref{eq:ns_velocity_system} by $\varepsilon$ and taking the limit $\varepsilon \to 0$ gives:
\begin{equation}\sum_{k \neq i} c_{ik} \left( u_i^{n+1} - u_k^{n+1} \right) = 0, \quad \forall i \in \{1, \dots, N\},
\end{equation}
where now $c_{ik} = \Delta t \, \nu_{ik} a_{ik}^{n+1} > 0$, since $\lambda_{ik}>0$ for $i\ne k$ and the densities are positive. We prove synchronization by a discrete maximum principle. For any spatial component, let $u_m^{n+1} = \max_i \{ u_i^{n+1} \}$ be the largest velocity among the species. For species $m$ the limit equation gives:
\begin{equation}
    \sum_{k \neq m} c_{mk} \left( u_m^{n+1} - u_k^{n+1} \right) = 0.
\end{equation}
By definition of the maximum, $u_m^{n+1} - u_k^{n+1} \ge 0$ for all $k$. Since $c_{mk} > 0$, all the terms of the sum are non-negative, so the sum vanishes only if $u_k^{n+1} = u_m^{n+1}$ for all $k$. Hence all the velocities are equal: $u_1^{n+1} = u_2^{n+1} = \dots = u^{n+1}$. For a fixed explicit state, the velocity system also gives $u_i^{n+1}-u_k^{n+1}=\mathcal O(\varepsilon)$, since the collision graph is connected. Thus the right-hand side of \eqref{eq:ns_temperature_system} remains bounded as $\varepsilon\to0$. Multiplying that system by $\varepsilon$ and applying the same maximum principle argument to the temperature system \eqref{eq:ns_temperature_system} in the limit $\varepsilon \to 0$, with positive effective coefficients $c_{ik} = \Delta t \, \nu_{ik} b_{ik} > 0$, shows that all the temperatures are equal: $T_1^{n+1} = \dots = T^{n+1}$. The common values $u^{n+1}$ and $T^{n+1}$ are then fixed by the conservation of the total momentum and energy, $(\rho u)^{n+1} = (\rho u)^*$ and $E^{n+1} = E^*$.

\emph{Step 2 (pointwise relaxation).}
Since $u_i = u_k = u$ and $T_i = T_k = T$, the target velocities and temperatures reduce to $u_{ik} = u$ and $T_{ik} = T$, and every target Maxwellian becomes $M_{ik}^{n+1} = M_i(u^{n+1}, T^{n+1})$.
Taking the limit $\varepsilon \to 0$ of the pointwise update \eqref{eq:final_pointwise_update}:
\begin{equation}
    \lim_{\varepsilon \to 0} f_i^{n+1} = \lim_{\varepsilon \to 0} \frac{\varepsilon f_i^* + \Delta t \sum_{k=1}^N \nu_{ik} \, n_i^{n+1} M_i(u^{n+1}, T^{n+1})}{\varepsilon + \Delta t \sum_{k=1}^N \nu_{ik}} = n_i^{n+1} M_i(u^{n+1}, T^{n+1}).
\end{equation}

\emph{Step 3 (recovery of the time-discrete Euler system).}
We project the time-discrete IMEX kinetic scheme \eqref{eq:imex_update_formal} onto the collisional invariants $\phi_i(v) = (1, m_i v, \frac{1}{2}m_i|v|^2)^\intercal$ and sum over all species $i=1, \dots, N$. By conservation of mass, momentum, and energy during collisions, the stiff interaction terms vanish identically for any $\varepsilon > 0$:
\begin{equation}
    \sum_{i=1}^N \int_{\mathbb{R}^3} \phi_i(v) \left[ \frac{1}{\varepsilon} \sum_{k=1}^N \nu_{ik} \left( n_i^{n+1} M_{ik}^{n+1} - f_i^{n+1} \right) \right] dv = 0.
\end{equation}
In the limit $\varepsilon \to 0$, the distribution functions are the local Maxwellians $f_i^{n+1} = n_i^{n+1} M_i(u^{n+1}, T^{n+1})$ derived in Step 2. For well-prepared data $f_i^n=n_i^nM_i(u^n,T^n)$, or after the initial relaxation step in the limit, integrating the advective fluxes at level $n$ yields the consistent time-discrete single-fluid Euler equations for the mixture:
\begin{equation}\label{eq:semidiscrete_euler}
    \begin{cases}
        \displaystyle \frac{n_i^{n+1} - n_i^n}{\Delta t} + \nabla_x \cdot (n_i^n u^n) = 0, \\[10pt]
        \displaystyle \frac{(\rho u)^{n+1} - (\rho u)^n}{\Delta t} + \nabla_x \cdot \left( \rho^n u^n \otimes u^n \right) + \nabla_x (n^n T^n) = 0, \\[10pt]
        \displaystyle \frac{E^{n+1} - E^n}{\Delta t} + \nabla_x \cdot \left[ \left( E^n + n^n T^n \right) u^n \right] = 0,
    \end{cases}
\end{equation}
where $E = \frac{1}{2}\rho |u|^2 + \frac{3}{2}nT$ is the total mixture energy. Hence the IMEX scheme is AP: in the limit $\varepsilon \to 0$ it gives a consistent time discretization of the Euler system \eqref{eq7}, with no restriction of order $\varepsilon$ on $\Delta t$.
\end{proof}

\subsection{Time integration of the hybrid scheme}
\label{sec:UnifiedFramework}

In this section we describe the time discretization and the coupling between the two models. The same IMEX Runge--Kutta scheme is used in all cells, the kinetic solver is AP, and the coupling is done at each stage, which preserves the order in time when the partition is fixed.

\subsubsection{A unified cross-model IMEX scheme}
\label{sub:UnifiedIMEX}

We use the same time discretization for a hierarchy of models $\{\mathscr{M}_0, \mathscr{M}_1\}$ of increasing accuracy and cost, where $\mathscr{M}_0$ corresponds to the multi-species Euler equations (Layer 0) shown in Equation \eqref{eq7}, and $\mathscr{M}_1$ corresponds to the multi-species BGK kinetic model (Layer 1) in Equation \eqref{rescaledBGK}. The construction can be extended to more than two models \cite{cap-par-rey-2025, cap-par-rey-2026}.

To this end, let us consider the system of equations governing each species $i = 1, \dots, N$, where the right-hand side can be additively split into a non-stiff transport operator $F(U)$, treated explicitly, and a stiff collision/relaxation operator $G(U)$, treated implicitly to avoid severe time-step restrictions. Here, $U$ represents a generalized state vector containing the variables of all species in the cell (either the distribution functions $f_{i}$ or the conservative variables $(n_1,\dots,n_N,\rho u,E)$ of the mixture). At time level $n$, each spatial cell $P_c$ (for $c=1, \dots, N_P$) is associated with explicit and implicit operators $F_{c}^{l_c^n}$ and $G_{c}^{l_c^n}$, corresponding to the locally active model $\mathscr{M}_{l_c^n}$ with $l_c^n \in \{0, 1\}$. Note that the model in each spatial cell is the same for all species. In each cell the solution is advanced with the active model only, that is, with the distribution functions in kinetic cells and with the conservative variables in macroscopic cells; the variables of the other model are recomputed from these when they are needed by the coupling described in Section~\ref{sub:coupling}.

The operators $F_{c}^{l_c^n}$ and $G_{c}^{l_c^n}$ depend on the active model:
\begin{itemize}
    \item in $\mathscr{M}_1$, the transport term $-v \cdot \nabla_x f_i$ is in the explicit part $F_{c}^{1}(U)$, and the collision operator (intra- and inter-species) is in the implicit part $G_{c}^{1}(U)$, solved as in Section \ref{sec:numerical_properties};
    \item in $\mathscr{M}_0$, the fluxes are in the explicit part $F_{c}^{0}(U)$ and $G_{c}^{0}(U) = 0$, so in these cells the scheme reduces to the explicit Runge--Kutta scheme with tableau $\tilde{A}$ given below.
\end{itemize}

To match the third-order spatial accuracy provided by the CWENO3 reconstruction, we employ a third-order, $q$-stage Implicit-Explicit (IMEX) Runge--Kutta scheme. The implicit step is the one of Section \ref{sec:numerical_properties}; here we describe its extension to high order. The stage values in cell $P_c$ are computed in the same form for all models:
\begin{equation}
\label{eq:IMEX_system}
\begin{cases}
U_{c}^{(s)} = \displaystyle U_{c}^n + \Delta t \sum_{j=1}^{s-1} \tilde{a}_{sj} F_{c}^{l_c^n}(t_n + \tilde{c}_j \Delta t, U^{(j)}) + \Delta t \sum_{j=1}^s a_{sj} G_{c}^{l_c^n}(t_n + c_j \Delta t, U^{(j)}), \quad s=1,\dots,q, \\[3mm]
U_{c}^{n+1} = \displaystyle U_{c}^n + \Delta t \sum_{j=1}^q \tilde{\omega}_j F_{c}^{l_c^n}(t_n + \tilde{c}_j \Delta t, U^{(j)}) + \Delta t \sum_{j=1}^q \omega_j G_{c}^{l_c^n}(t_n + c_j \Delta t, U^{(j)}).
\end{cases}
\end{equation}
The scheme is defined by the $q \times q$ matrices $\tilde{A}=(\tilde{a}_{sj})$ (explicit tableau, with $\tilde{a}_{sj}=0$ for $j \ge s$) and $A=(a_{sj})$ (Diagonally Implicit tableau, DIRK, with $a_{sj}=0$ for $j > s$), along with the coefficient vectors $\tilde{c}, c, \tilde{\omega}, \omega$, represented in Butcher notation:
\begin{center}
\begin{tabular}{c|c}
    $\tilde{c}$ & $\tilde{A}$\\      
    \hline
    & $\tilde{\omega}^\intercal$\\
\end{tabular} \qquad \qquad
\begin{tabular}{c|c}
    $c$ & $A$\\                
    \hline
    & $\omega^\intercal$\\
\end{tabular}.
\end{center}
Since $A$ is diagonally implicit, each stage $s = 1, \dots, q$ is computed in two steps:
\begin{enumerate}
    \item explicit step:
    \begin{equation}
        U_{c}^{(s)*} = U_{c}^n + \Delta t \sum_{j=1}^{s-1} \tilde{a}_{sj} F_{c,j} + \Delta t \sum_{j=1}^{s-1} a_{sj} G_{c,j},
    \end{equation}
    where $F_{c,j} = F_{c}^{l_c^n}(t_n + \tilde{c}_j \Delta t, U^{(j)})$ and $G_{c,j} = G_{c}^{l_c^n}(t_n + c_j \Delta t, U^{(j)})$;
    \item implicit step, solved as in Section \ref{sec:numerical_properties} with $\Delta t$ replaced by $a_{ss}\Delta t$:
    \begin{equation}
    \label{eq:ImplicitCorrector}
        U_{c}^{(s)} = U_{c}^{(s)*} + \Delta t\, a_{ss} G_{c}^{l_c^n}(t_n + c_s \Delta t, U^{(s)}).
    \end{equation}
\end{enumerate}
After computing all $q$ stages sequentially, the numerical solution is updated via:
\begin{equation}
    U_{c}^{n+1} = U_{c}^n + \Delta t \sum_{j=1}^q \tilde{\omega}_j F_{c,j} + \Delta t \sum_{j=1}^q \omega_j G_{c,j}.
\end{equation}

\subsubsection{Penalization and GSA property}
\label{sub:AP_Penalization_GSA}
The time discretization has to remain stable and accurate in the stiff limit $\varepsilon\to0$ without restrictions on the time step, that is it has to be AP. For the transport operator $F(U)$ we use the CWENO3 reconstruction and the numerical fluxes described above. For the collision operator $G(U)$, a third-order approximation of the cell average requires a suitable quadrature rule within each spatial cell. 

Denoting the spatial cell as $P_c$, the evolution equation for the cell average $\bar{U}_c(t)$ involves the high-order approximation $\bar{G}^{HO}_c$ of the collision operator. Treating $\bar{G}^{HO}_c$ implicitly would require the solution of a nonlinear system coupling all the quadrature points of the cell. To avoid this, we use a penalization technique \cite{fil-jin-2010, bos-dim-2021}. We introduce the low-order operator $\bar{\mathcal{P}}_{c}^{LO}(U)$, that is the multi-species BGK operator evaluated on the cell average, and we add and subtract it in the semi-discrete system:
\begin{equation}
\label{eq:penalized_system}
    \bar{U}_c' = \bar{F}_c(U) + \left( \bar{G}^{HO}_c(U) - \bar{\mathcal{P}}_{c}^{LO}(U) \right) + \bar{\mathcal{P}}_{c}^{LO}(U).
\end{equation}
In the IMEX scheme the operator $\bar{\mathcal{P}}_{c}^{LO}$ is treated implicitly. Since it involves only the cell average, its inversion only requires the small linear systems of Section \ref{sec:numerical_properties}. The residual $\left( \bar{G}^{HO}_c(U) - \bar{\mathcal{P}}_{c}^{LO}(U) \right)$, which contains the high-order corrections, is treated explicitly. 

The residual is still stiff, and treating it explicitly may destroy the AP property. For this reason we use a \textit{Globally Stiffly Accurate} (GSA) scheme \cite{bos-par-rus-2013}. An IMEX scheme is GSA if its coefficients satisfy
\begin{equation}
\omega_s = a_{q s} \quad \text{and} \quad \tilde{\omega}_s = \tilde{a}_{q s} \qquad \forall s=1,\dots,q.
\end{equation}

With this condition the numerical solution coincides with the last stage, $U_c^{n+1}\equiv U_c^{(q)}$. This identity is useful for preserving the equilibrium constraint at the end of a time step, but GSA alone is not sufficient to control an arbitrary explicit stiff residual: the penalization must also be compatible with the equilibrium limit. The first-order AP argument is given in Section~\ref{subsec:AP_Euler}; the higher-order penalized construction follows \cite{bos-dim-2021}, and its stiff-regime behaviour is assessed numerically in \hyperref[sec:test1]{Test~1}. We use the GSA scheme BPR(3,5,3) \cite{bos-par-rus-2013} (also referred to as BPR(3,4,3) in \cite{bos-dim-2021}), whose explicit and implicit Butcher tableaux are
\begin{center}
\begin{tabular}{c|c c c c c}
$0$   & $0$   & $0$   & $0$   & $0$ & $0$\\
$1$   & $1$   & $0$   & $0$   & $0$ & $0$\\
$2/3$ & $4/9$ & $2/9$ & $0$   & $0$ & $0$\\
$1$   & $1/4$ & $0$   & $3/4$ & $0$ & $0$\\
$1$   & $1/4$ & $0$   & $3/4$ & $0$ & $0$\\
\hline
      & $1/4$ & $0$   & $3/4$ & $0$ & $0$\\
\end{tabular} \qquad \qquad
\begin{tabular}{c|c c c c c}
$0$   & $0$   & $0$    & $0$   & $0$    & $0$\\
$1$   & $1/2$ & $1/2$  & $0$   & $0$    & $0$\\
$2/3$ & $5/18$& $-1/9$ & $1/2$ & $0$    & $0$\\
$1$   & $1/2$ & $0$    & $0$   & $1/2$  & $0$\\
$1$   & $1/4$ & $0$    & $3/4$ & $-1/2$ & $1/2$\\
\hline
      & $1/4$ & $0$    & $3/4$ & $-1/2$ & $1/2$\\
\end{tabular}.
\end{center}

\subsubsection{Interface coupling}\label{sub:coupling}
Once the temporal integrator is defined, the numerical solution must be coupled across neighbouring cells belonging to different physical regimes. Let $\mathcal{T}$ denote the spatial mesh composed of cells $P_c$ ($c=1, \dots, N_P$). Since all species within a given spatial cell share the same physical regime in the present model, we define a single regime indicator for cell $c$:
\begin{equation}
\chi_{c} \in \{\mathrm{K}, \mathrm{E}\}, \quad c=1,\dots,N_P,
\end{equation}
where $\mathrm{K}$ denotes a kinetic description (BGK, Layer 1) and $\mathrm{E}$ denotes a macroscopic description (Euler, Layer 0). Let $\mathcal{N}_c$ denote the spatial numerical stencil associated with cell $c$. Because neighbouring cells can be in different regimes, evaluating the spatial fluxes and the collision terms requires a coupling strategy, which we describe below.

\paragraph{\bf Kinetic cells.}
When a cell $c$ is in the kinetic regime ($\chi_{c}=\mathrm{K}$), all species within it are treated kinetically. The implicit operator $\bar{\mathcal{P}}_{c}^{LO}$ only involves the cell itself. The numerical transport fluxes and the high-order collision term $\bar{G}^{HO}_c$ require the reconstruction of the distribution functions, and hence their values throughout the spatial stencil $\mathcal{N}_c$. Consequently, if a neighbouring cell $j \in \mathcal{N}_c$ belongs to the macroscopic regime ($\chi_{j}=\mathrm{E}$), its kinetic distribution functions must be locally reconstructed as Maxwellians for all species. The effective distribution function $\widetilde{f}_{i,m,j}$ used in these reconstructions for species $i$ is defined as:
\begin{equation}
\label{eq:kinetic_coupling}
\widetilde{f}_{i,m,j} =
\begin{cases}
f_{i,m,j}, & \text{if } \chi_{j}=\mathrm{K}, \\[2mm]
n_{i,j}\,M_{i,m}(u_j, T_j), & \text{if } \chi_{j}=\mathrm{E}.
\end{cases}
\end{equation}

\paragraph{\bf Macroscopic cells.}
When a cell $c$ is in the macroscopic regime ($\chi_{c}=\mathrm{E}$), all species within it are described by the mixture Euler equations \eqref{eq7}. In this regime the solution evolves only by transport, since $G_{c}^{0} = 0$. Computing the numerical fluxes across cell interfaces requires knowing the conservative variables of all species throughout the spatial stencil $\mathcal{N}_c$. Consequently, if a neighbouring cell $j \in \mathcal{N}_c$ belongs to the kinetic regime ($\chi_{j}=\mathrm{K}$), its effective macroscopic state must be locally recovered by integrating its discrete distribution functions against the collision invariants $\phi_i(v) = (1,\, m_i v,\, \frac{1}{2} m_i \vert{}v\vert{}^2)^\intercal$ and summing momentum and energy over the species. The effective macroscopic state $\widetilde{U}_{j}$ used for the Euler fluxes is defined as:
\begin{equation}
\label{eq:macro_coupling}
\widetilde{U}_{j} =
\begin{cases}
U_{j}, & \text{if } \chi_{j}=\mathrm{E}, \\[3mm]
\displaystyle \Big(\int_{\mathbb{R}^{d_v}} f_{1,j}\,dv,\ \dots,\ \int_{\mathbb{R}^{d_v}} f_{N,j}\,dv,\ \sum_{i=1}^N\int_{\mathbb{R}^{d_v}} m_i v\, f_{i,j}\,dv,\ \sum_{i=1}^N\int_{\mathbb{R}^{d_v}} \frac{1}{2}m_i|v|^2 f_{i,j}\,dv\Big)^\intercal, & \text{if } \chi_{j}=\mathrm{K}.
\end{cases}
\end{equation}

\paragraph{\bf Regime update and synchronization of the stages.} The coupling described above passes from kinetic to macroscopic data by computing moments, and from macroscopic to kinetic data by Maxwellian reconstruction. It is applied at two points of each time step:
\begin{enumerate}
\item \emph{Beginning of the time step.} The partition is updated from $\chi_c^n$ to $\chi_c^{n+1}$ with the criterion of Section~\ref{secIndicators}. Cells that change regime update their variables: if $\chi_c^n=\mathrm{K}$ and $\chi_c^{n+1}=\mathrm{E}$, the moments $U_c^n$ are computed by quadrature of $f_{i,c}^n$; if $\chi_c^n=\mathrm{E}$ and $\chi_c^{n+1}=\mathrm{K}$, the distribution functions are set to the Maxwellians $f_{i,m,c}^n=n_{i,c}^n M_{i,m}(u_c^n,T_c^n)$. The first stage of BPR(3,5,3) has no implicit part and coincides with the solution at $t^n$, so its fluxes use the data at $t^n$, taken according to the previous partition. If $\chi_c^{n+1}=\mathrm{K}$, the transport fluxes and the high-order collision term $\bar{G}^{HO}_c$ need the distribution functions in the stencil $\mathcal{N}_c$: for a neighbour $j$ with $\chi_j^n=\mathrm{E}$ the Maxwellians $n_{i,j}^n M_{i,m}(u_j^n,T_j^n)$, built from the moments of cell $j$, are used, and otherwise $f_{i,j}^n$. If $\chi_c^{n+1}=\mathrm{E}$, the Euler fluxes need the moments in the stencil: for a neighbour $j$ with $\chi_j^n=\mathrm{K}$ they are computed by quadrature of $f_{i,j}^n$, and otherwise $U_j^n$ is used. In particular, a neighbour that has just switched from $\mathrm{K}$ to $\mathrm{E}$ still provides its distribution functions at $t^n$, since it was kinetic at that time.
\item \emph{After each stage.} After each implicit stage the coupling is applied again, now with the new partition $\chi_c^{n+1}$: in kinetic cells the moments are computed by quadrature of the current distribution functions, and in macroscopic cells the Maxwellians $n_{i,c}M_{i,m}(u_c,T_c)$ are built from the current moments. In particular, from the second stage on, a cell that has switched from $\mathrm{K}$ to $\mathrm{E}$ no longer provides $f_{i,c}^n$ to its kinetic neighbours, but the Maxwellians built from its current moments. All the other cells provide data according to their current regime. In this way, at each stage all the cells of a stencil provide data at the same stage. Exchanging data only at the end of the time step would introduce a time lag at the interfaces and reduce the accuracy there to first order \cite{cap-ten-2025}.
\end{enumerate}
When the partition is fixed, the stage-by-stage coupling preserves the order of the IMEX scheme, as proved in \cite{cap-ten-2025} and observed in \hyperref[sec:test1]{Test 1}. When the partition changes this is no longer guaranteed, but the first-order time lag at the interfaces is avoided.

\paragraph{\bf Lifting.}
Macroscopic states are lifted with the Maxwellian $n_{i,j}M_{i,m}(u_j,T_j)$. Higher-order reconstructions, for example based on the Chapman--Enskog expansion, do not in general preserve the non-negativity of the distribution function, and negative values of $f_{i,m}$ may cause instabilities. The analytic Maxwellian used for lifting is non-negative, and no instability was observed in these tests. The conservative $L^2$ correction of Section~\ref{secConMas}, however, does not in general preserve pointwise non-negativity.

\subsection{Conservation of macroscopic quantities}\label{secConMas}
It is well known that the discretization of the velocity space on discrete grids introduces a loss of conservation \cite{mie-2000, dim-par-2014}. A natural solution to this problem is to consider a velocity support large enough that the energy outside the domain is negligible. However, the truncation of the domain is not the only issue, as quadrature errors on the grid must also be taken into consideration. It is therefore important that the solver preserves these quantities.

Several approaches have been proposed in the literature, exploiting either entropy minimization (as in \cite{mie-2000}) or $L^2$-minimization (as in \cite{gam-tha-2009, bos-cho-rus-2022}). The \emph{entropy approach} allows the construction of a conservative discrete Maxwellian, while $L^2$-minimization can be applied to more general distribution functions. In this paper, we implement the $L^2$-minimization approach proposed in \cite{gam-tha-2009}. As in \cite{PaperA}, we apply it to each Maxwellian that appears in the scheme, for each species. A weighted version of this projection was recently used for gas mixtures in \cite{cap-par-rey-2025-2}.

\section{Simulations}\label{secSimulations}
In this section we numerically test the proposed scheme. The numerical simulations are executed on a computer equipped with two \emph{Intel(R) Xeon(R) Gold 5418Y} CPUs and 128 GB of RAM, kindly provided by the Mathematics and Computer Science Department of the University of Ferrara. The machine was shared with other users, so the execution times may have been affected by other jobs running at the same time and should be regarded as indicative. The implementation does not require specialized hardware and can also be run on a laptop. All simulations have been executed using the third-order in space and time solvers discussed above.

The time step $\Delta t$ employed in the simulations is generally determined by the relation:
\[
\Delta t = 0.5\,\frac{\Delta h}{v_{\max}},
\]
where $\Delta h=\min(\Delta x, \Delta y)$ and the velocity grid covers $[-v_{\max},v_{\max}]^2$.
This time step is used in all simulations except the convergence test in time (\hyperref[sec:test1]{Test 1}), where the time step is varied.
For the sake of simplicity, in all our simulations we used the following values $\nu_{ik}=\nu_{ik}^{(0)}=1$ and $\lambda_{ik}=1$.
All our simulations are executed on a 2D domain in space and velocity. The 2D velocity domain of the kinetic regime is discretized using a uniform Cartesian grid. All the quantities of Sections~\ref{secModel} and~\ref{secNumSch} are computed accordingly in velocity dimension $d_v=2$. The indicators of Section~\ref{secIndicators} are computed with the three-dimensional formulas of that section, with the velocities and the gradients extended by a zero third component. 

\subsection{Test 1: Time convergence test}
\label{sec:test1}
This test checks the third-order convergence in time of the BPR(3,5,3) scheme for the hybrid solver, that is for the two solvers and for their coupling (Section \ref{sec:UnifiedFramework}). To this end, the regime of each cell is chosen at random, independently of the other cells, and is kept fixed during the simulation. As proved in \cite{cap-ten-2025}, with a fixed partition the coupling of Section \ref{sec:UnifiedFramework} preserves the order of the IMEX scheme.
We consider two species with mass ratio $m_1/m_2=2$, in a physical domain $[0,10]\times[0,0.5]$ discretized with 40 cells along the $x$-direction and 2 cells along the $y$-direction. In the macroscopic cells, the Euler state is initialized with the following macroscopic moments $W_1=(n_1, u^x, u^y, T)$ and $W_2=(n_2, u^x, u^y, T)$ for the first and second species, respectively:
\begin{equation}    
    W_1(x,y) = W_2(x,y) = (1+0.01\exp{(-2 r_x^2)},\,0,\,0,\,1),
\end{equation}
where $r_x^2=(x-x_c)^2$ and $x_c=10/2$. In the kinetic cells, the distribution functions $f_1$ and $f_2$ for both species are initialized as Maxwellian distributions with moments $W_1(x, y)$ and $W_2(x,y)$. The velocity space $[-8,8]^2$ is discretized using 32 Cartesian points in both directions. The initial data depend only on $x$, and zero-gradient boundary conditions are imposed in both directions.

\begin{figure}[htbp]
    \centering
    \includegraphics[width=0.6\linewidth]{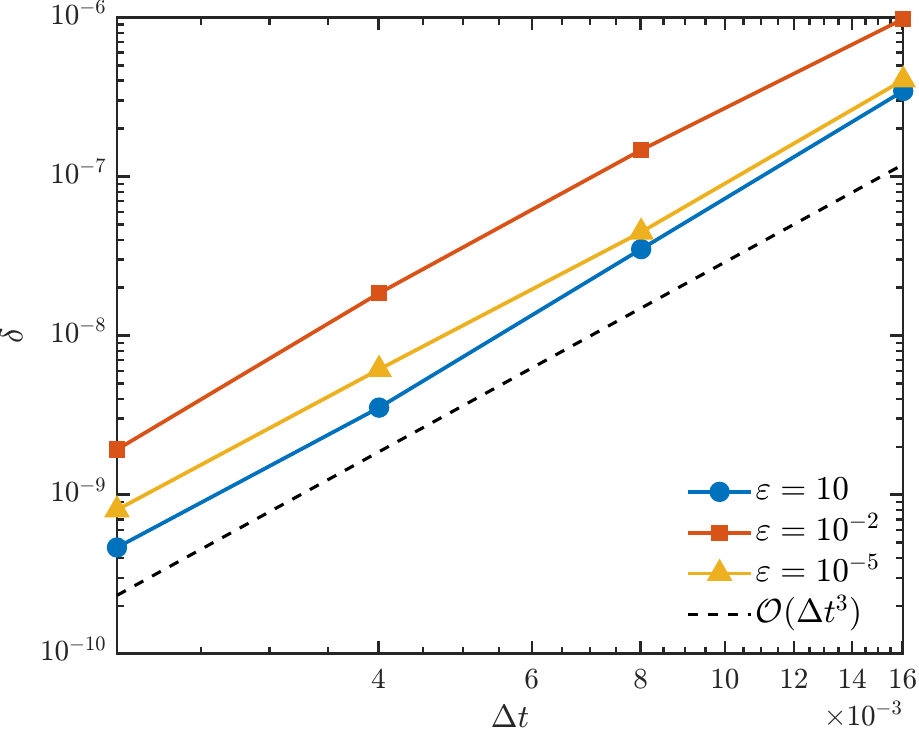}
    \caption{Temporal convergence rates of the BPR(3,5,3) scheme for the benchmark described in \hyperref[sec:test1]{Test 1}. The numerical error $\delta$ ($L^1$ norm) is plotted against the time step $\Delta t$ for $\varepsilon = 10, 10^{-2}, 10^{-5}$, with a fixed random partition. The observed order is close to three for all values of $\varepsilon$.}
    \label{Fig1}
\end{figure}

We take three values of the Knudsen number, $\varepsilon \in \{10^{-5}, 10^{-2}, 10^{1}\}$. For each value of $\varepsilon$ we vary the time step $\Delta t$ and compute the $L^1$ error $\delta$ of the density $\rho$ at the final time, with respect to a reference solution computed with the smallest time step. The corresponding convergence rates are shown in Figure \ref{Fig1}. 

These results show third-order temporal convergence for the three values of $\varepsilon$ tested, including the stiff case $\varepsilon=10^{-5}$. They provide numerical evidence of asymptotic accuracy for this benchmark, with no observed order reduction at the sampled Knudsen numbers.

\subsection{Test 2: Shock tube problem}
\label{sec:test2}
This test is inspired by \cite{bis-bos-dim-gro-mar-2022, bai-rey-2022} and its aim is to verify the ability of the scheme to capture flow discontinuities in a one-dimensional Riemann problem for a gas mixture. In the collision-dominated regime considered here ($\varepsilon=10^{-4}$) the Euler equations already describe the flow well, and the test is used for three checks. First, the kinetic and the Euler solvers for the mixture are compared with the analytic solution of the Riemann problem. Second, we check that the exchange of information across the interfaces between kinetic and Euler cells does not produce oscillations. Finally, the test illustrates the AP property of the scheme and shows where the indicator activates the kinetic solver.

The computational domain is the box $[0,1]\times[0, 0.02]$ discretized using 100 cells along the $x$-direction and 2 cells along the $y$-direction. Zero-gradient boundary conditions are imposed in both directions. We use a two-dimensional mesh for this one-dimensional problem to check that the scheme preserves the 1D structure of the solution, without multidimensional artifacts. The velocity space  $[-15, 15]^2$ is discretized using 64 points per direction. The final time is chosen to be equal to $t_f = 0.15$.

In the macroscopic cells, the Euler solution is initialized with the following macroscopic moments $W_1=(\rho_1, u^x, u^y, T)$ and $W_2=(\rho_2, u^x, u^y, T)$ for the first and second species, respectively:
\begin{equation}
    \begin{cases}
        W_1(x<0.5,y) = ((1-\alpha )\rho_L,\,0,\,0,\,T_L),\\
        W_2(x<0.5,y) = (\alpha \rho_L,\,0,\,0,\,T_L),\\
        W_1(x\geq0.5,y) = (\alpha \rho_R,\,0,\,0,\,T_R),\\
        W_2(x\geq0.5,y) = ((1-\alpha)\rho_R,\,0,\,0,\,T_R),
    \end{cases}    
\end{equation}
where $\alpha=10^{-3}$, $\rho_L=1$, $\rho_R=2^{-3}$, $T_L=P_L/((1-\alpha )\rho_L/m_1)$, $T_R=P_R/((1-\alpha )\rho_R/m_2)$, $P_L=1$ and $P_R=2^{-5}$. In the kinetic cells, the distribution functions $f_1$ and $f_2$ for both species are initialized as Maxwellian distributions with moments $W_1(x, y)$ and $W_2(x,y)$.
The mass ratio of the species is $m_2/m_1=2$, and the threshold of the indicator is $\eta=\varepsilon$.

\begin{figure}
    \centering
    \includegraphics[width=0.96\textwidth]{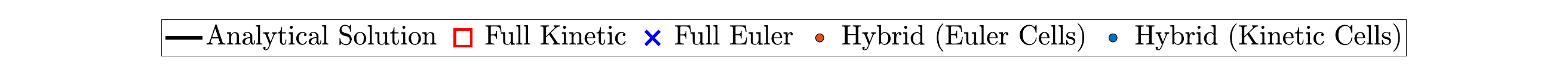}\par\smallskip
    \begin{subfigure}{0.32\textwidth}
        \centering
        \includegraphics[width=\linewidth]{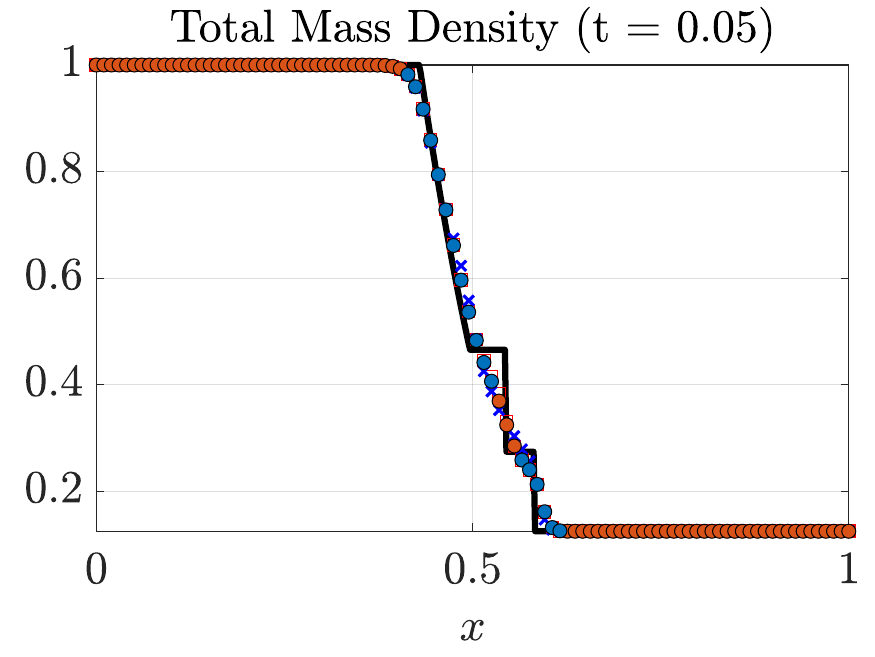}

    \end{subfigure}
    \begin{subfigure}{0.32\textwidth}
        \centering
        \includegraphics[width=\linewidth]{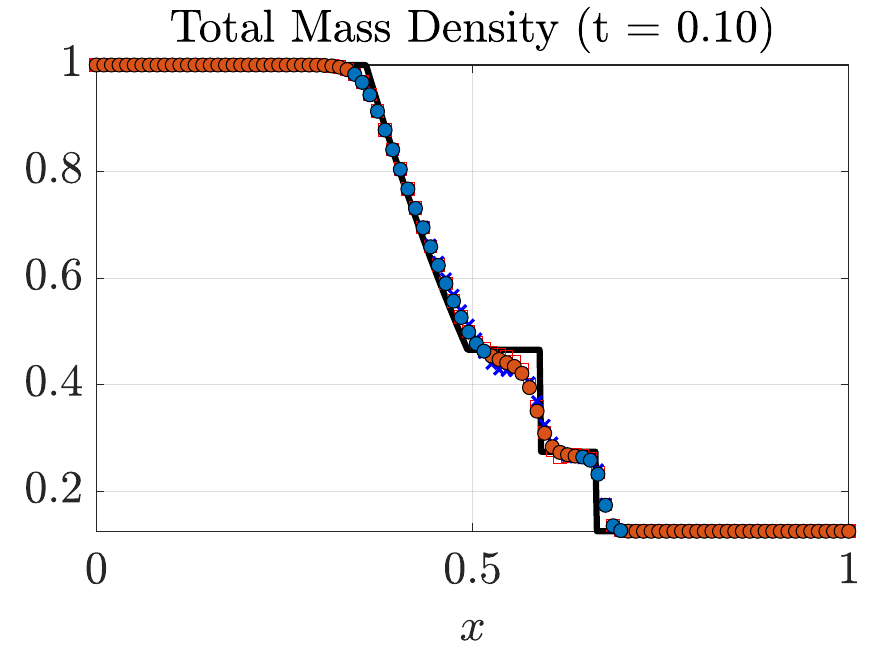}

    \end{subfigure}
    \begin{subfigure}{0.32\textwidth}
        \centering
        \includegraphics[width=\linewidth]{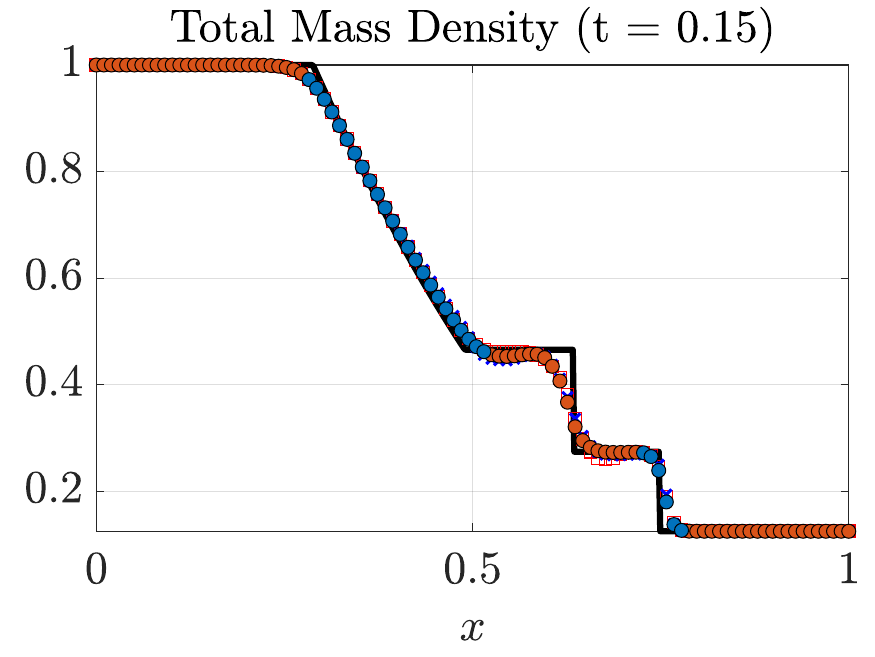}

    \end{subfigure}

    \begin{subfigure}{0.32\textwidth}
        \centering
        \includegraphics[width=\linewidth]{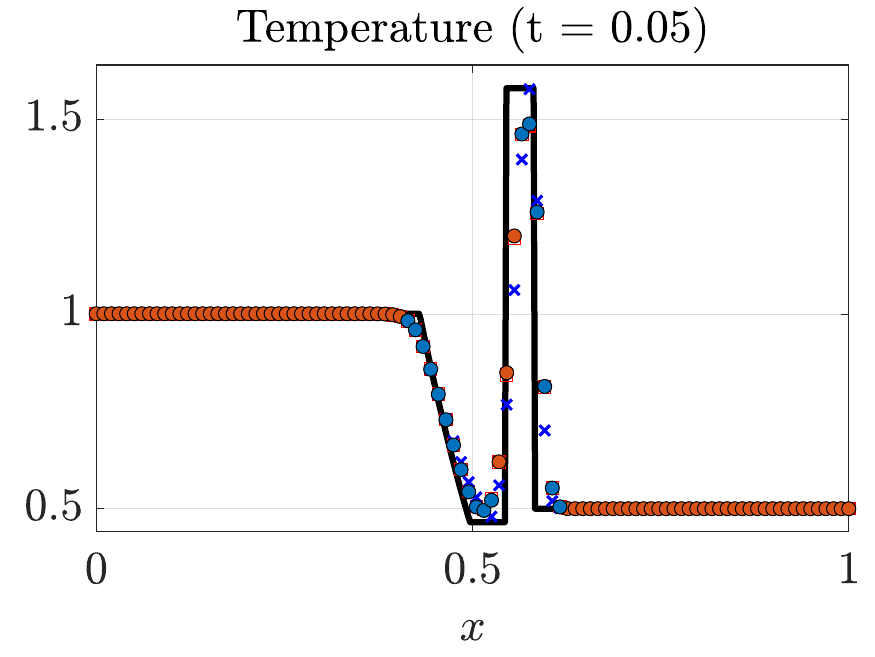}

    \end{subfigure}
    \begin{subfigure}{0.32\textwidth}
        \centering
        \includegraphics[width=\linewidth]{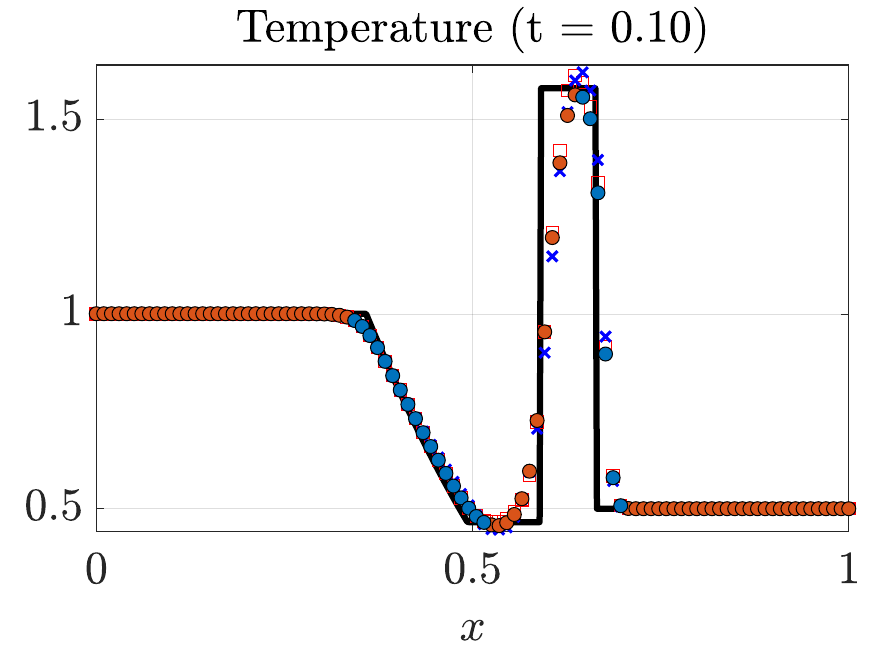}

    \end{subfigure}
    \begin{subfigure}{0.32\textwidth}
        \centering
        \includegraphics[width=\linewidth]{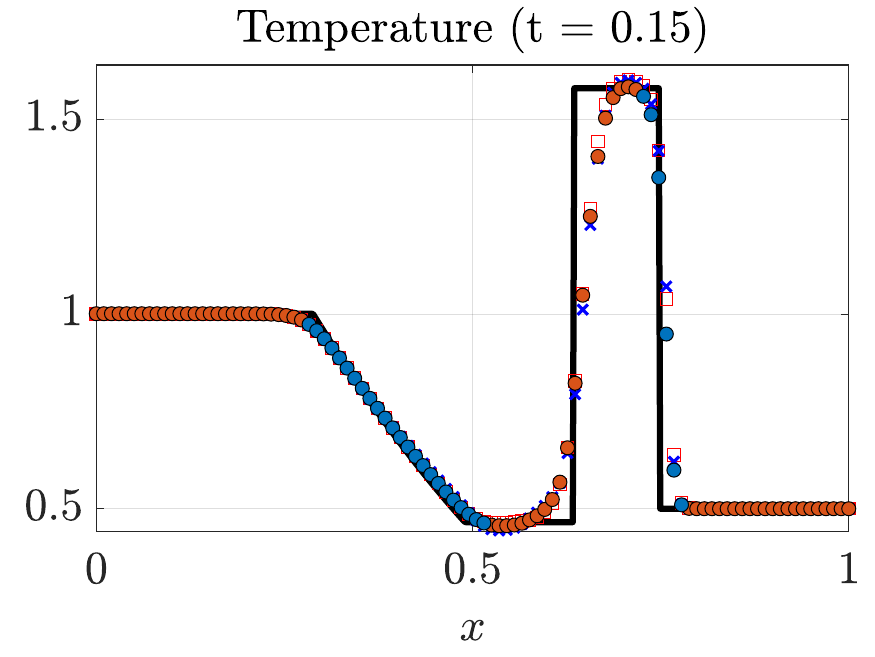}

    \end{subfigure}

    \begin{subfigure}{0.32\textwidth}
        \centering
        \includegraphics[width=\linewidth]{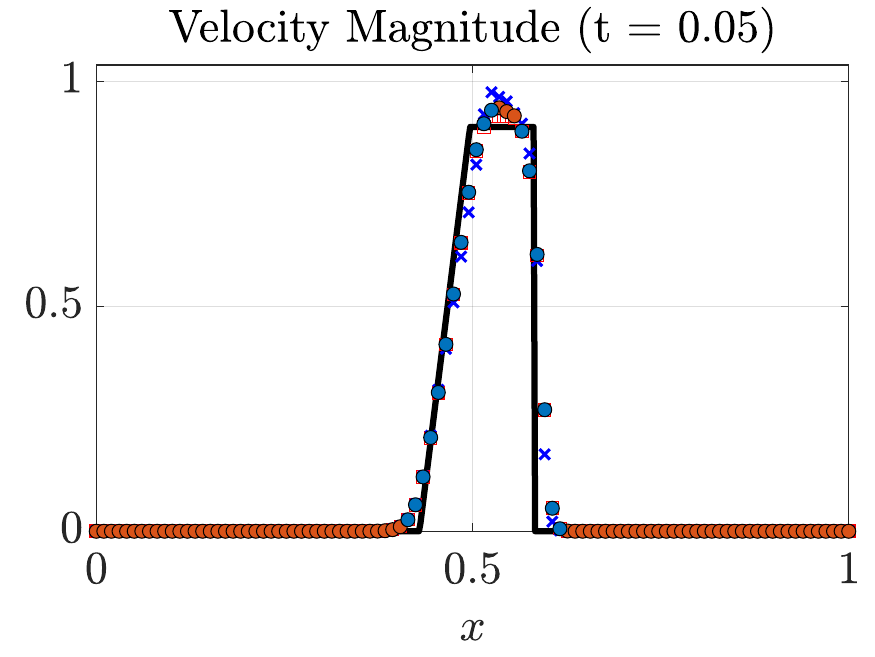}

    \end{subfigure}
    \begin{subfigure}{0.32\textwidth}
        \centering
        \includegraphics[width=\linewidth]{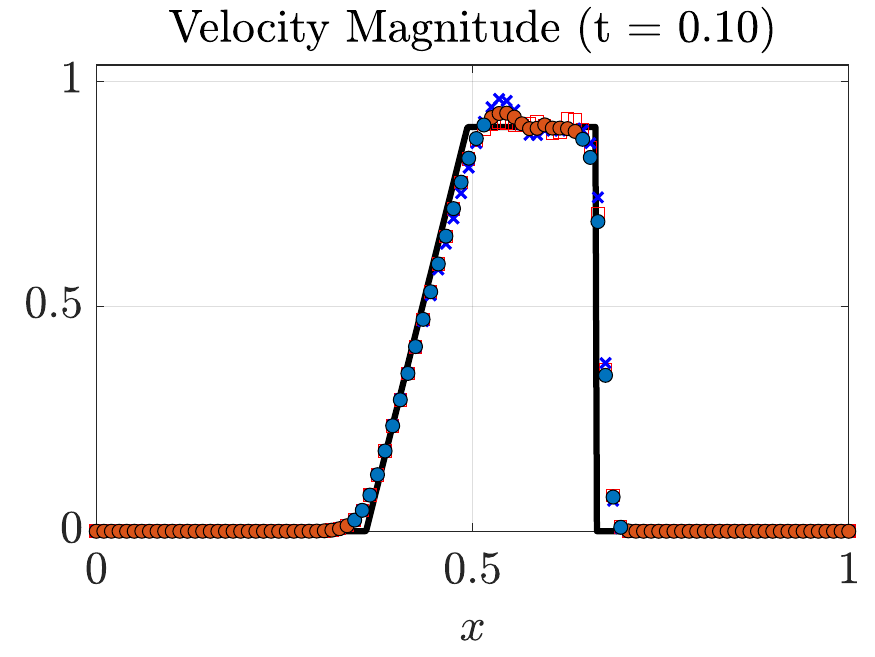}

    \end{subfigure}
    \begin{subfigure}{0.32\textwidth}
        \centering
        \includegraphics[width=\linewidth]{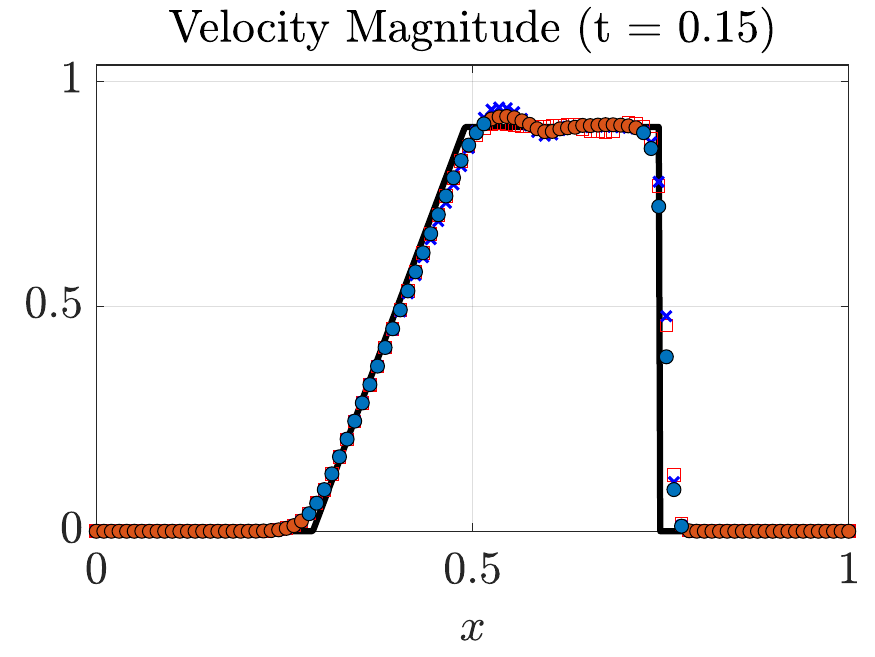}

    \end{subfigure}

    \caption{Total mass density (top), temperature (middle) and velocity magnitude (bottom) of the mixture along $x$ at fixed $y$ for \hyperref[sec:test2]{Test 2}. The profiles report the solutions obtained with the hybrid scheme, the full Euler solver and the full kinetic solver, and the analytic solution. The kinetic solver is active in regions where non-equilibrium effects would appear, and no oscillations are produced at the interfaces between kinetic and Euler cells.}
    \label{Fig2}
\end{figure}

Figure \ref{Fig2} presents the spatial evolution along the $x$-direction of the total mass density, temperature, and velocity magnitude. The profiles obtained with the hybrid scheme, the full Euler solver and the full kinetic solver are reported together with the analytic solution of the Riemann problem for the Euler equations of the mixture with $\gamma=2$, since $d_v=2$. Because of the small Knudsen number, the global solution is physically macroscopic; however, the indicator activates the kinetic solver in regions where non-equilibrium effects would appear.

The three solutions are close to each other and to the analytic one, and the coupling does not produce oscillations at the interfaces between kinetic and Euler cells. The partition changes in time: the cells left behind by the shock go back to the Euler solver, while the kinetic region follows the rarefaction fan as it widens. Because of the small Knudsen number all three solvers approximate the solution of the Euler equations, so the kinetic solution is not a better reference than the Euler one, and we do not report the errors \eqref{eq:errST} for this test.

The computing times are $107.3$ s for the full Euler solver, $2352$ s for the full kinetic solver and $628.7$ s for the hybrid scheme, which is about $3.7$ times faster than the kinetic solver.

\subsection{Test 3: Flow past an inclined square}
\label{sec:test3}
In this test a rarefied gas mixture flows past a square obstacle with a $45^\circ$ angle of attack, and we compare the hybrid, kinetic and Euler solutions around the obstacle.

The physical domain $ [0, 1]\times[0, 1]$ is discretized using 50 cells in both directions. The velocity space $[-8,8]^2$ is discretized using 16 points per direction. We consider a binary gas mixture with a mass ratio $m_2/m_1=2$. The fluid is initialized uniformly across the physical domain in a local thermodynamic equilibrium state. Specifically, the macroscopic moment vectors $W_1 = (n_1, u^x, u^y, T)$ and $W_2 = (n_2, u^x, u^y, T)$ for the first and second species, respectively, are prescribed as:
\begin{equation}
    \begin{cases}        
    W_1(x, y) = (1.0,\, 0.5,\, 0.5,\, 1.0), \\
    W_2(x, y) = (0.8,\, 0.5,\, 0.5,\, 1.0).
    \end{cases}
\end{equation}
The distribution functions $f_1$ and $f_2$ are initialized as local Maxwellian distributions generated by the corresponding macroscopic states $W_1(x, y)$ and $W_2(x, y)$.

The obstacle is the square formed by the $7\times7$ cells around the centre of the domain, with sides parallel to the coordinate axes; since the flow is directed along the diagonal, the angle of attack is $45^\circ$.

Zero-gradient boundary conditions are imposed on the external boundary of the domain, while reflective slip conditions are enforced on the solid surface of the square obstacle. The numerical simulation is executed in the rarefied regime with a Knudsen number $\varepsilon = 0.1$, up to a final time $t_f = 0.15$. The threshold of the indicator is $\eta=\varepsilon$, and the domain partition is dynamically updated at the beginning of each time step.

\begin{figure}[htbp]
    \centering
    \begin{subfigure}{0.32\textwidth}
        \centering
        \includegraphics[width=\linewidth]{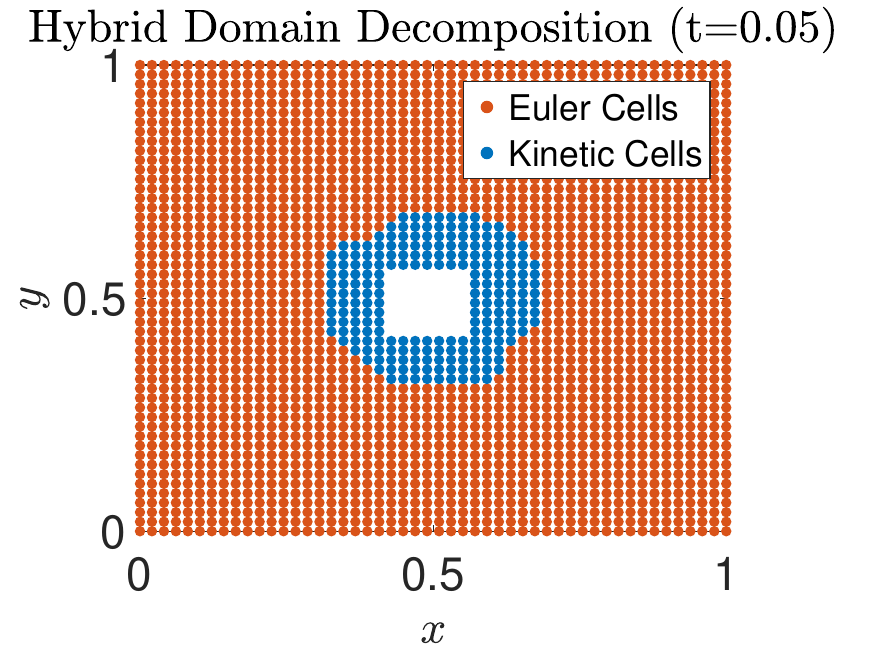}        
    \end{subfigure}
    \hfill
    \begin{subfigure}{0.32\textwidth}
        \centering
        \includegraphics[width=\linewidth]{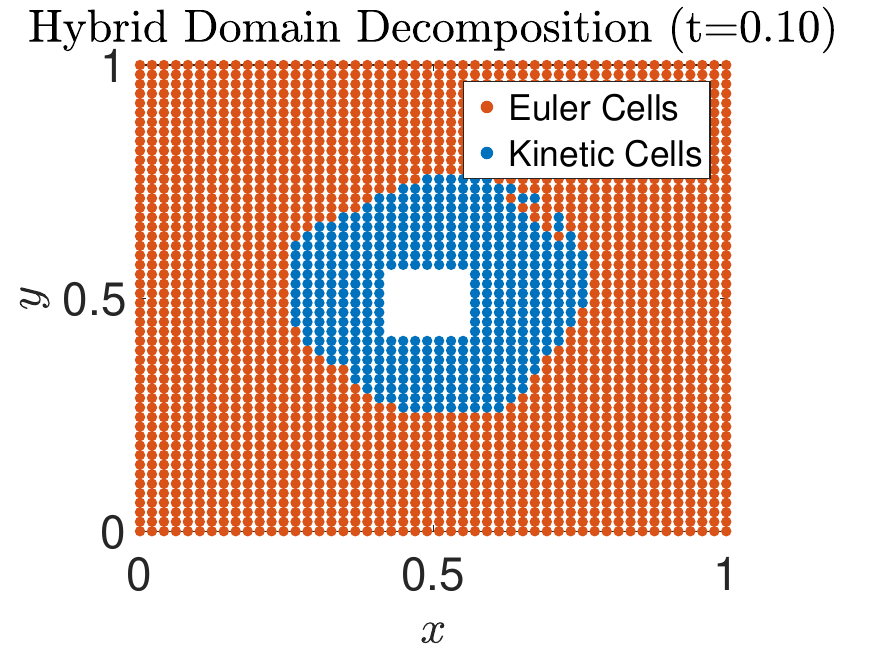}        
    \end{subfigure}
    \hfill
    \begin{subfigure}{0.32\textwidth}
        \centering
        \includegraphics[width=\linewidth]{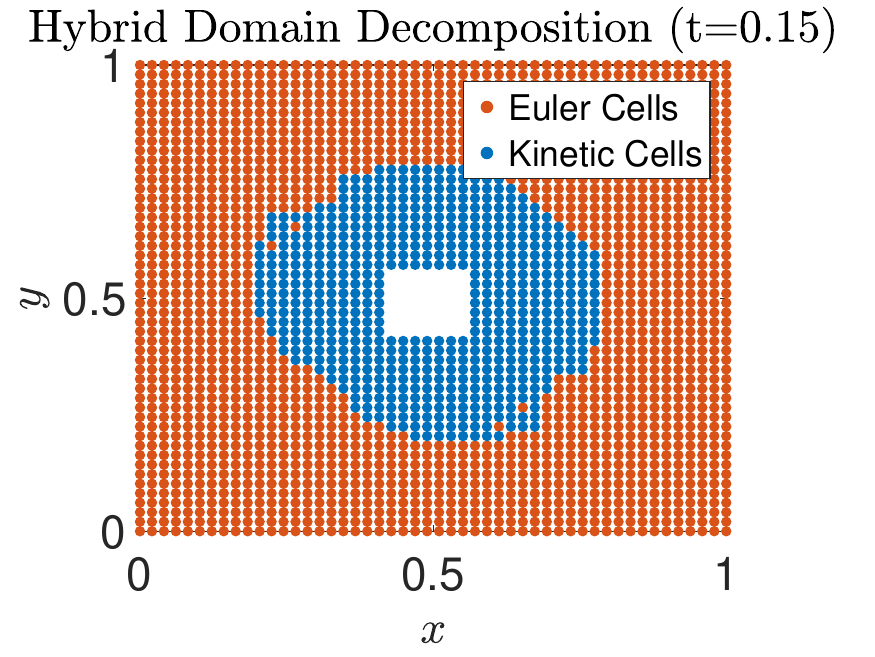}        
    \end{subfigure}

    \begin{subfigure}{0.32\textwidth}
        \centering
        \includegraphics[width=\linewidth]{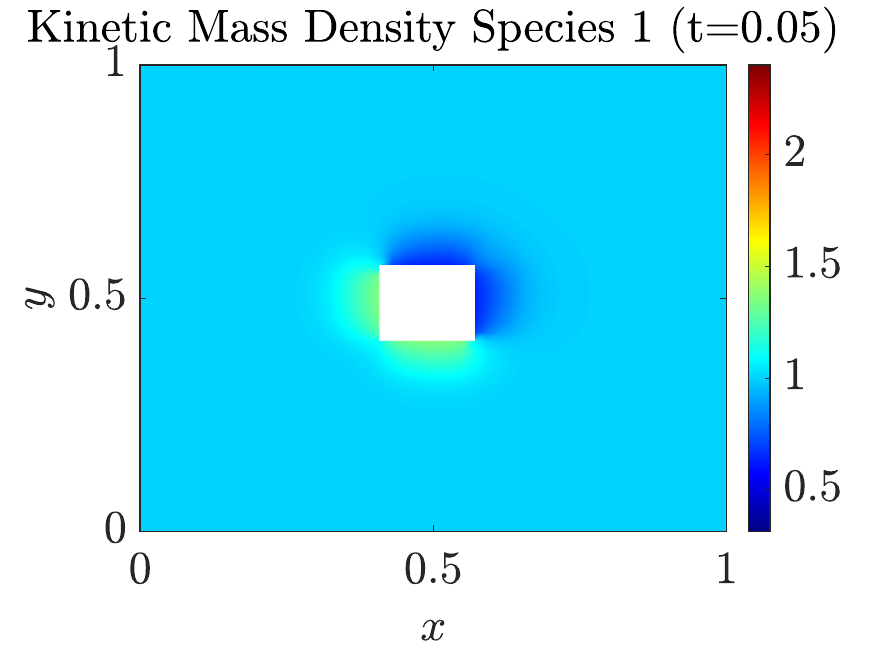}
    \end{subfigure}
    \hfill
    \begin{subfigure}{0.32\textwidth}
        \centering
        \includegraphics[width=\linewidth]{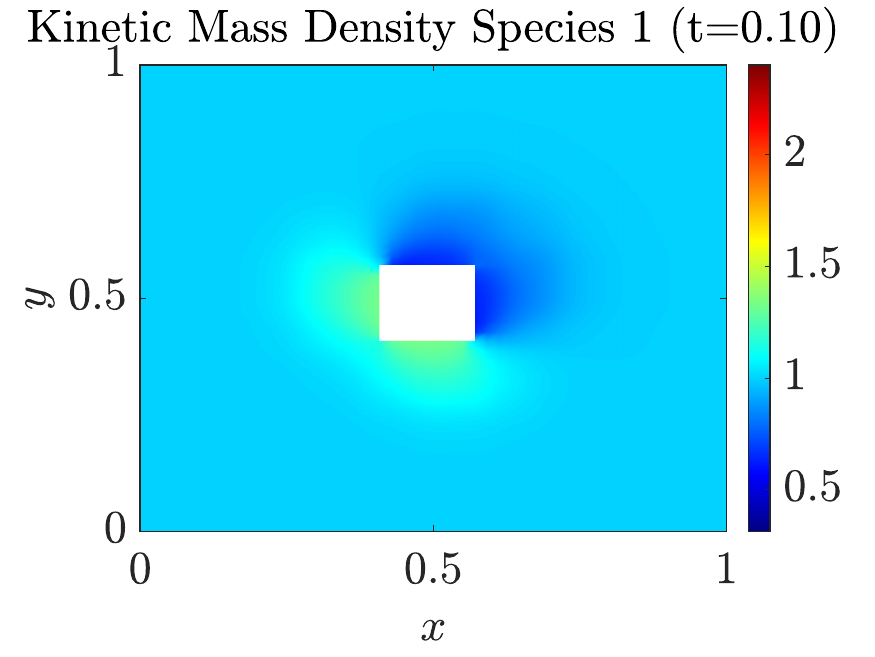}
    \end{subfigure}
    \hfill
    \begin{subfigure}{0.32\textwidth}
        \centering
        \includegraphics[width=\linewidth]{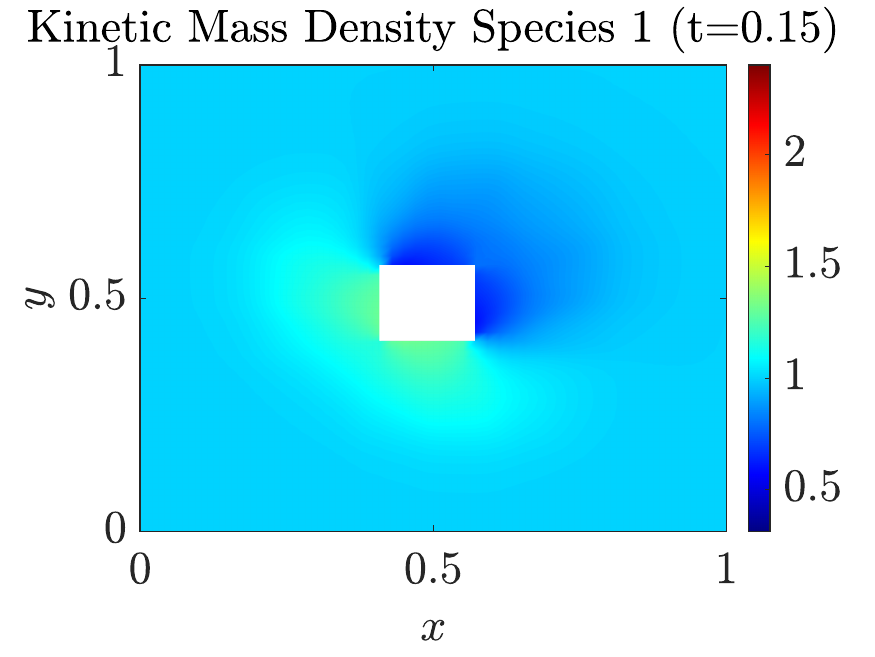}
    \end{subfigure}

    \begin{subfigure}{0.32\textwidth}
        \centering
        \includegraphics[width=\linewidth]{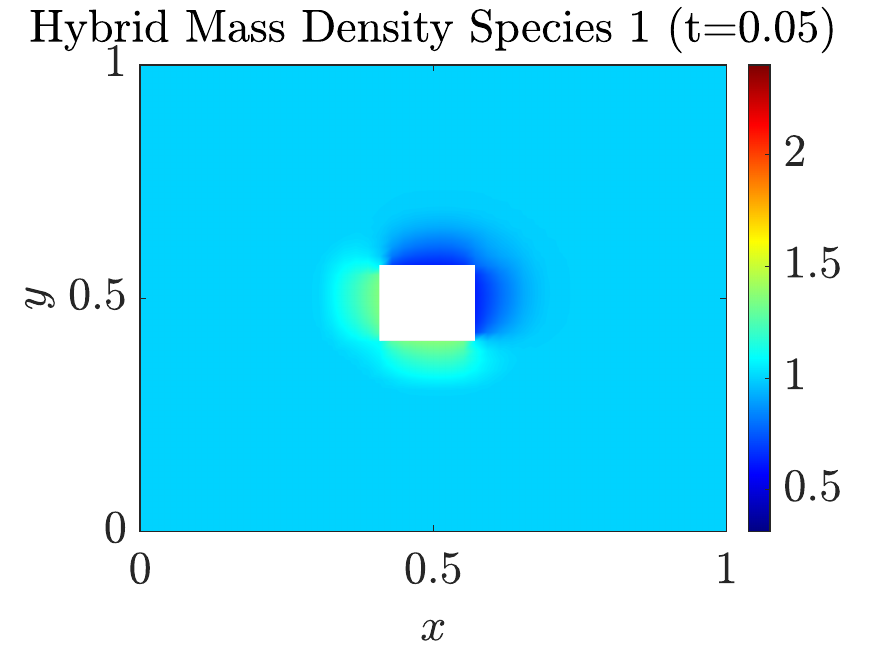}
    \end{subfigure}
    \hfill
    \begin{subfigure}{0.32\textwidth}
        \centering
        \includegraphics[width=\linewidth]{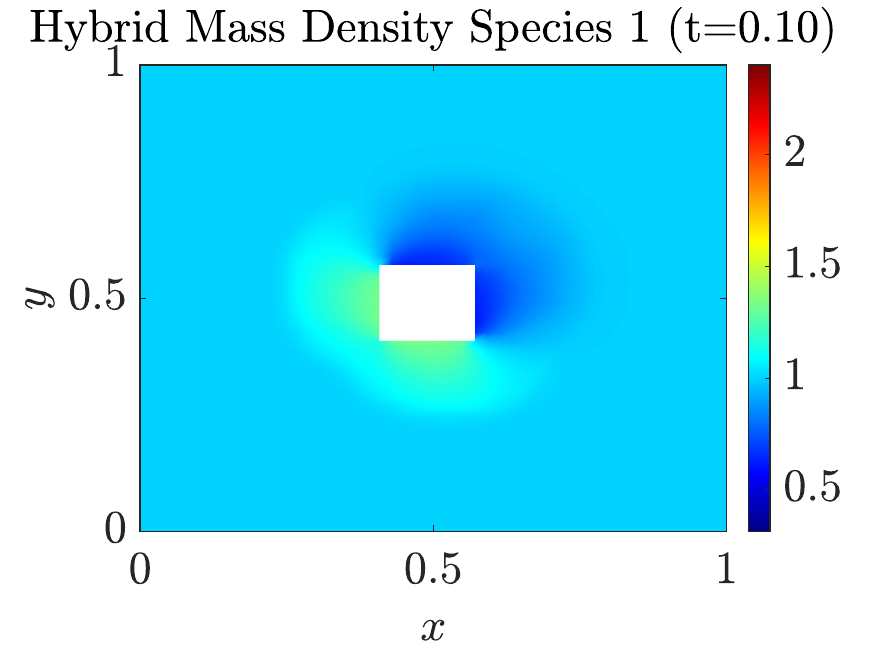}
    \end{subfigure}
    \hfill
    \begin{subfigure}{0.32\textwidth}
        \centering
        \includegraphics[width=\linewidth]{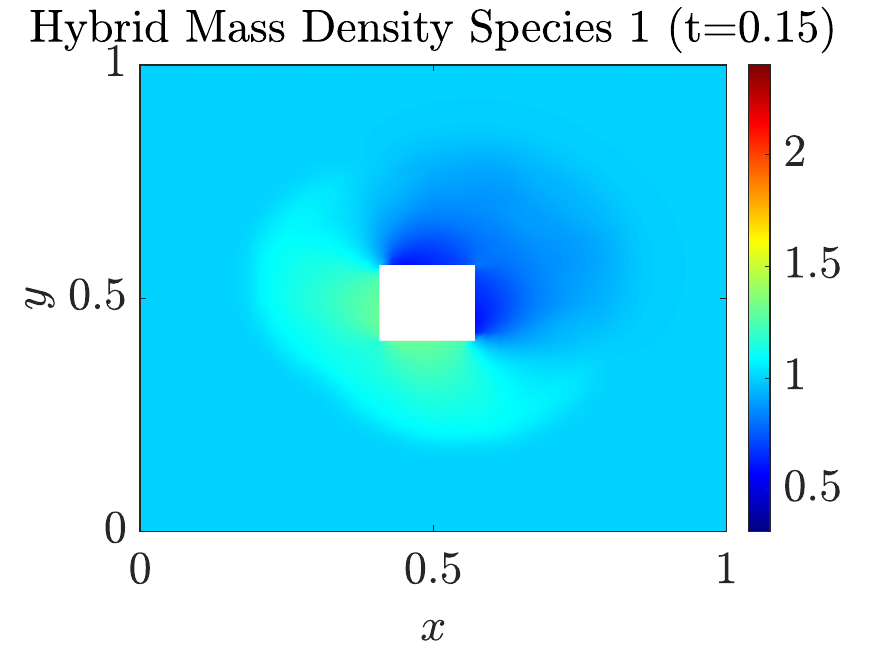}
    \end{subfigure}

    \begin{subfigure}{0.32\textwidth}
        \centering
        \includegraphics[width=\linewidth]{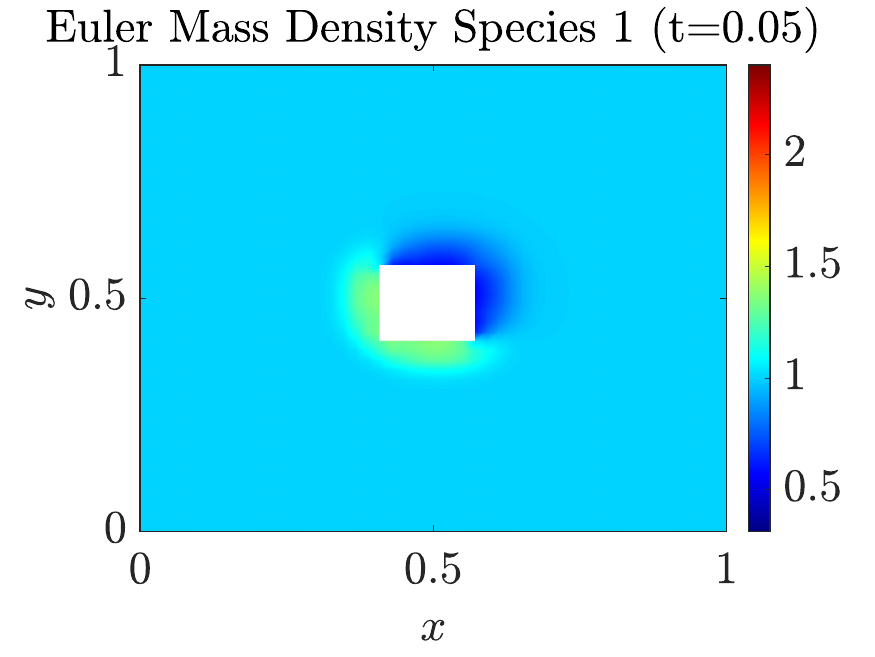}
    \end{subfigure}
    \hfill
    \begin{subfigure}{0.32\textwidth}
        \centering
        \includegraphics[width=\linewidth]{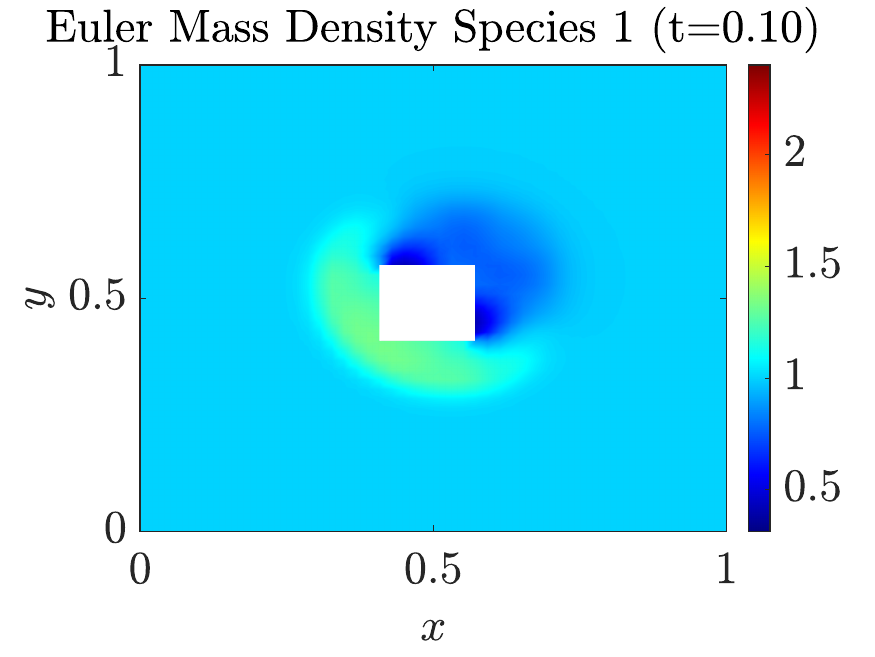}
    \end{subfigure}
    \hfill
    \begin{subfigure}{0.32\textwidth}
        \centering
        \includegraphics[width=\linewidth]{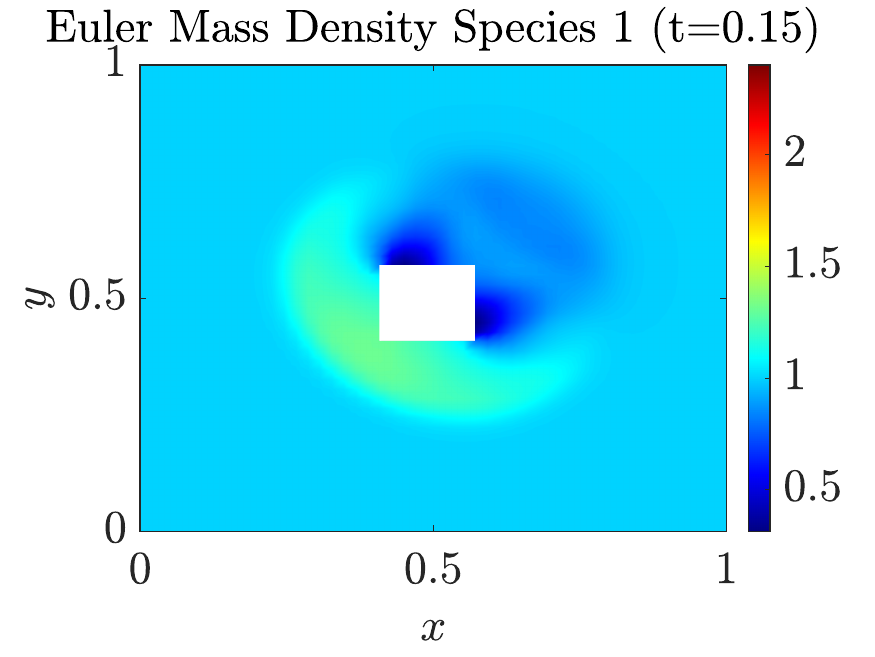}
    \end{subfigure}

    \caption{\hyperref[sec:test3]{Test 3}: Domain adaptation and mass density $\rho_1$ evolution for the first species. The top row illustrates the dynamic domain partition. The second, third, and fourth rows display the mass density distribution obtained using the full kinetic solver, the hybrid solver, and the full Euler solver, respectively.}
    \label{Fig3-Density-Species1}
\end{figure}

\begin{figure}[htbp]
    \centering
    \begin{subfigure}{0.32\textwidth}
        \centering
        \includegraphics[width=\linewidth]{TEST0_HybridDomain_t0.05.pdf}
    \end{subfigure}
    \hfill
    \begin{subfigure}{0.32\textwidth}
        \centering
        \includegraphics[width=\linewidth]{TEST0_HybridDomain_t0.10.pdf}
    \end{subfigure}
    \hfill
    \begin{subfigure}{0.32\textwidth}
        \centering
        \includegraphics[width=\linewidth]{TEST0_HybridDomain_t0.15.pdf}
    \end{subfigure}

    \begin{subfigure}{0.32\textwidth}
        \centering
        \includegraphics[width=\linewidth]{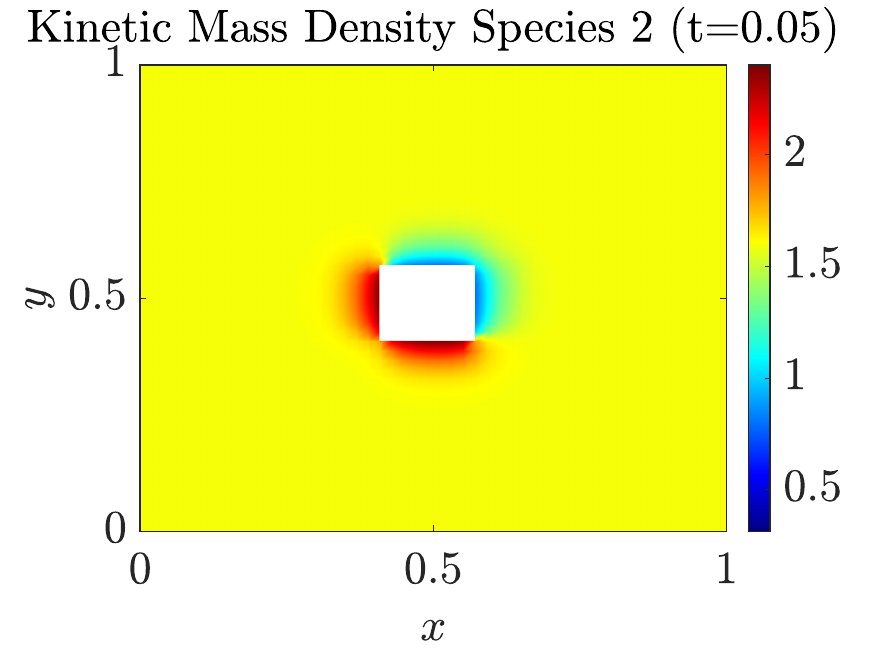}
    \end{subfigure}
    \hfill
    \begin{subfigure}{0.32\textwidth}
        \centering
        \includegraphics[width=\linewidth]{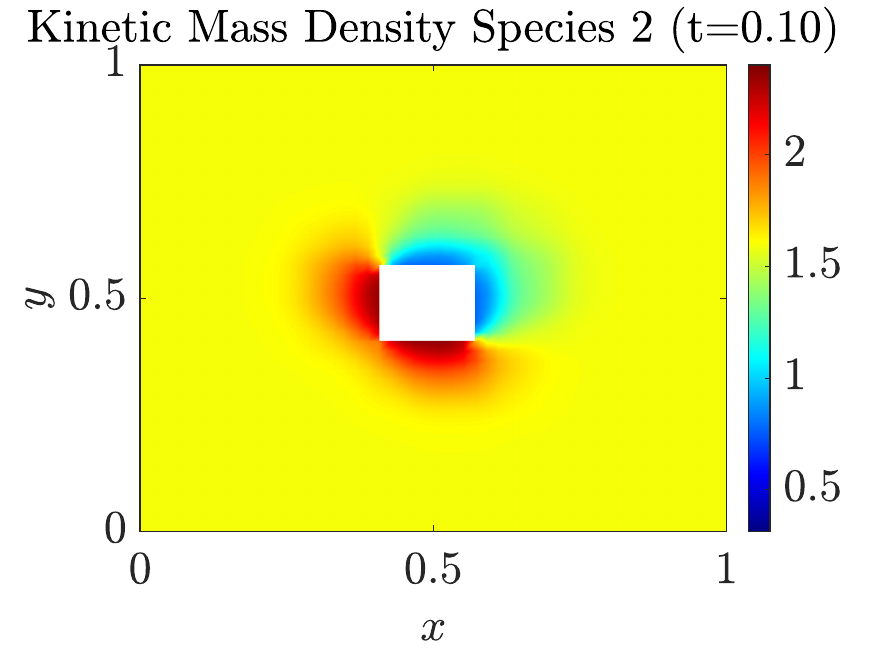}
    \end{subfigure}
    \hfill
    \begin{subfigure}{0.32\textwidth}
        \centering
        \includegraphics[width=\linewidth]{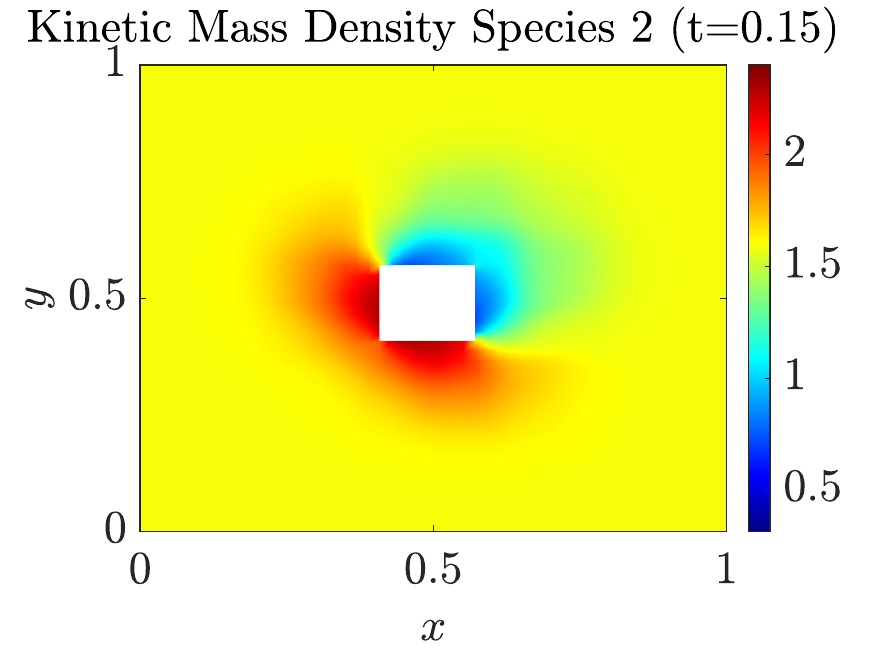}
    \end{subfigure}

    \begin{subfigure}{0.32\textwidth}
        \centering
        \includegraphics[width=\linewidth]{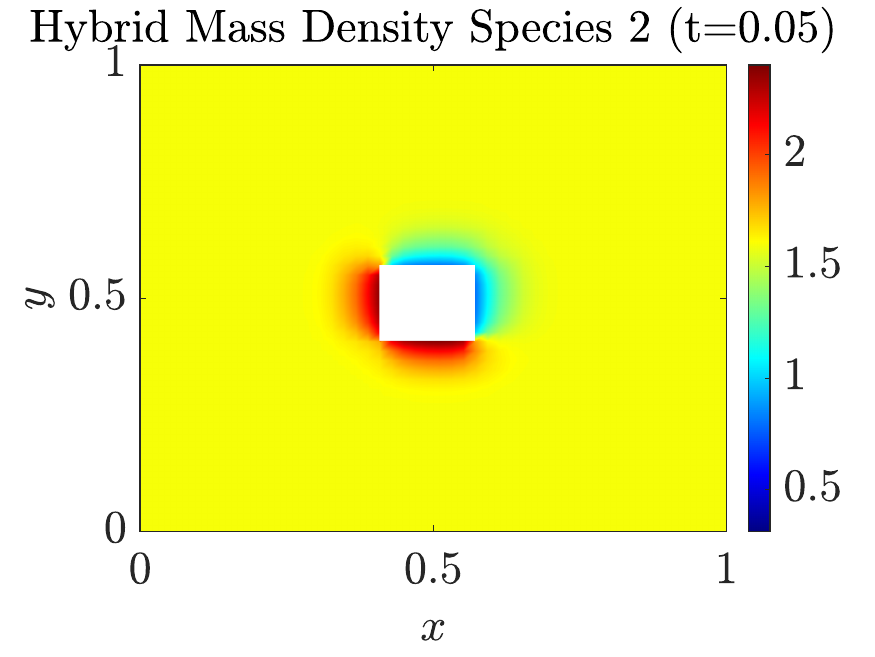}
    \end{subfigure}
    \hfill
    \begin{subfigure}{0.32\textwidth}
        \centering
        \includegraphics[width=\linewidth]{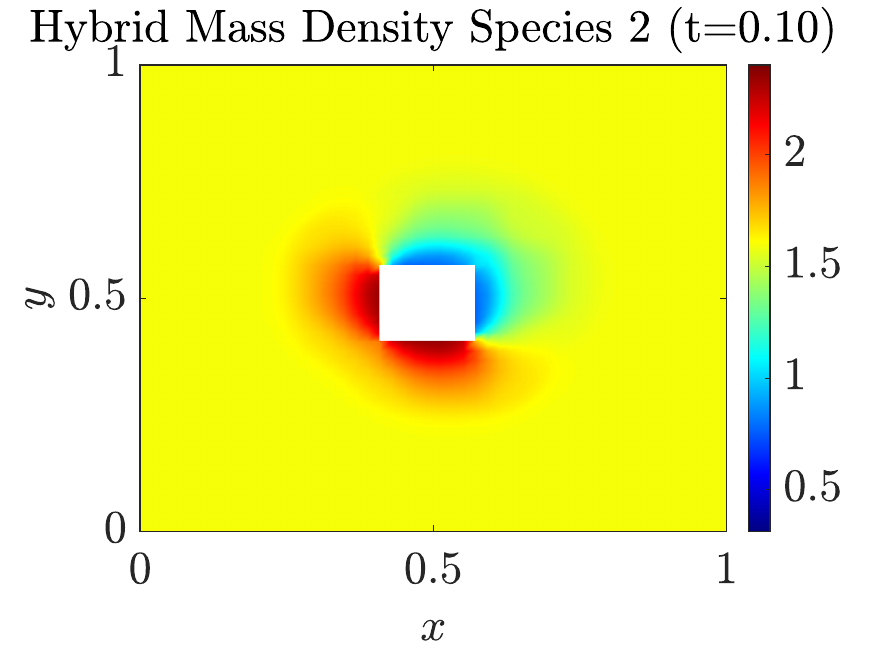}
    \end{subfigure}
    \hfill
    \begin{subfigure}{0.32\textwidth}
        \centering
        \includegraphics[width=\linewidth]{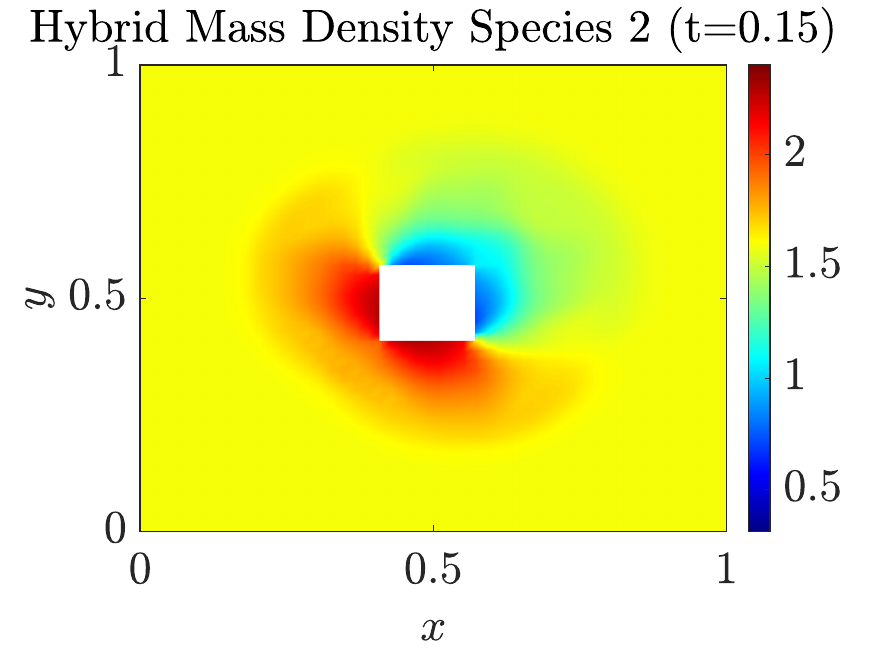}
    \end{subfigure}

    \begin{subfigure}{0.32\textwidth}
        \centering
        \includegraphics[width=\linewidth]{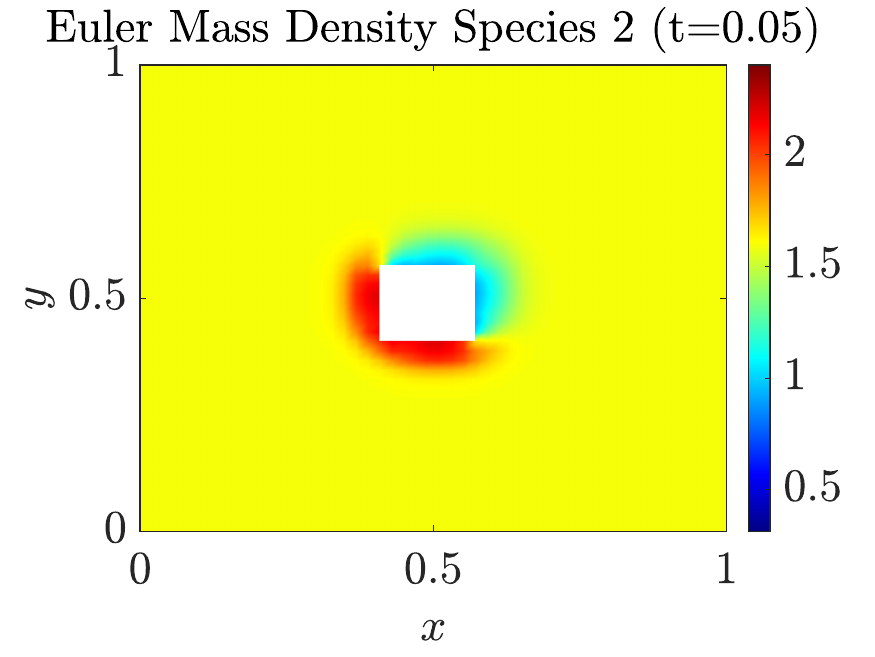}
    \end{subfigure}
    \hfill
    \begin{subfigure}{0.32\textwidth}
        \centering
        \includegraphics[width=\linewidth]{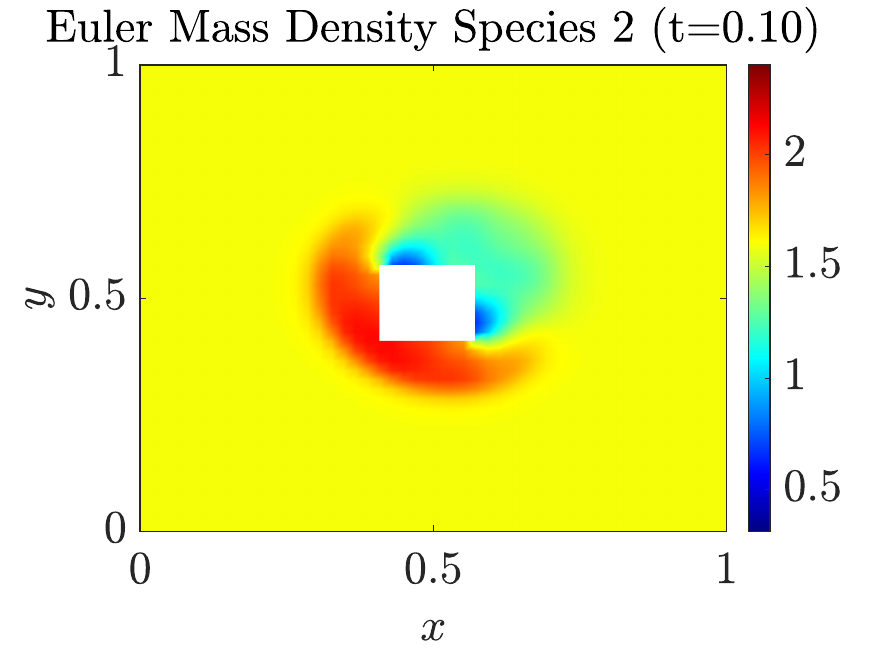}
    \end{subfigure}
    \hfill
    \begin{subfigure}{0.32\textwidth}
        \centering
        \includegraphics[width=\linewidth]{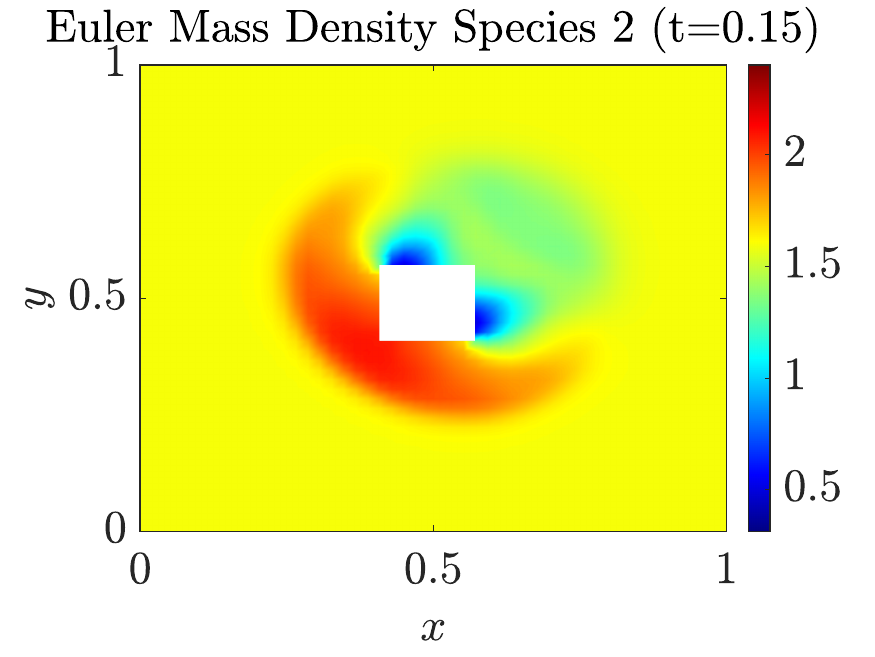}
    \end{subfigure}

    \caption{\hyperref[sec:test3]{Test 3}: Domain adaptation and mass density $\rho_2$ evolution for the second species. The top row illustrates the dynamic domain partition. The second, third, and fourth rows display the mass density distribution obtained using the full kinetic solver, the hybrid solver, and the full Euler solver, respectively.}
    \label{Fig4-Density-Species2}
\end{figure}

\begin{figure}[htbp]
    \centering
    \begin{subfigure}{0.32\textwidth}
        \centering
        \includegraphics[width=\linewidth]{TEST0_HybridDomain_t0.05.pdf}
    \end{subfigure}
    \hfill
    \begin{subfigure}{0.32\textwidth}
        \centering
        \includegraphics[width=\linewidth]{TEST0_HybridDomain_t0.10.pdf}
    \end{subfigure}
    \hfill
    \begin{subfigure}{0.32\textwidth}
        \centering
        \includegraphics[width=\linewidth]{TEST0_HybridDomain_t0.15.pdf}
    \end{subfigure}

    \begin{subfigure}{0.32\textwidth}
        \centering
        \includegraphics[width=\linewidth]{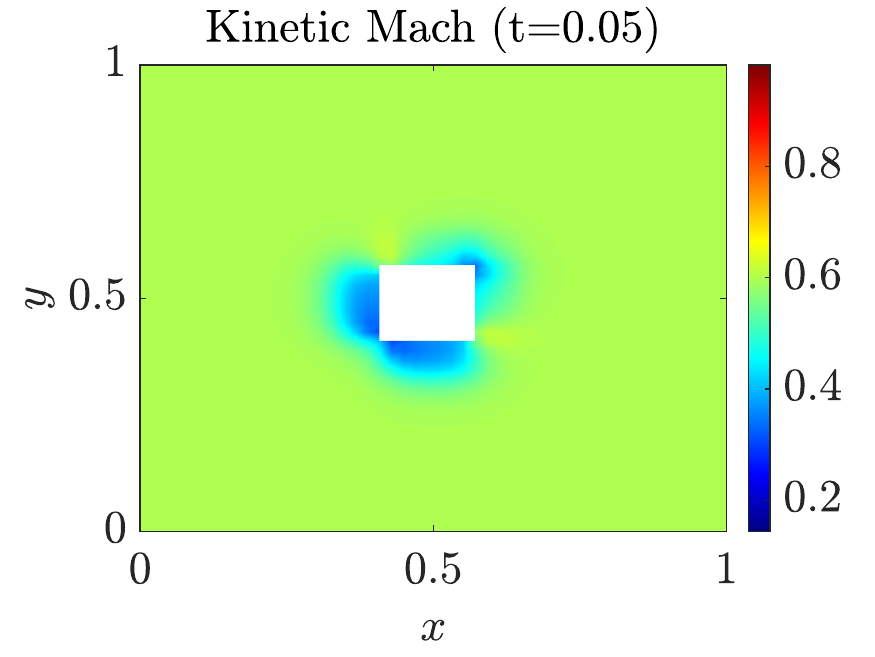}
    \end{subfigure}
    \hfill
    \begin{subfigure}{0.32\textwidth}
        \centering
        \includegraphics[width=\linewidth]{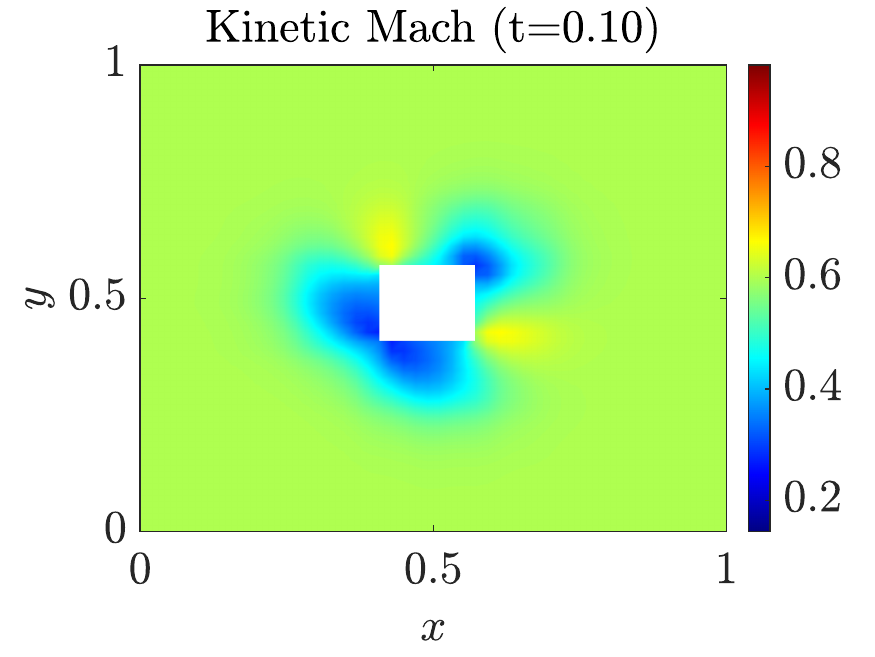}
    \end{subfigure}
    \hfill
    \begin{subfigure}{0.32\textwidth}
        \centering
        \includegraphics[width=\linewidth]{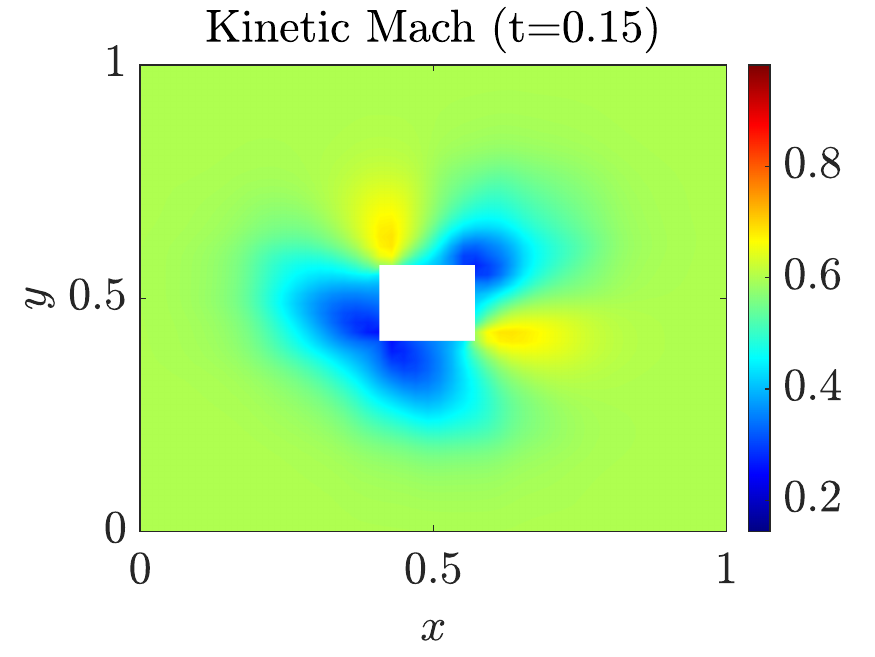}
    \end{subfigure}

    \begin{subfigure}{0.32\textwidth}
        \centering
        \includegraphics[width=\linewidth]{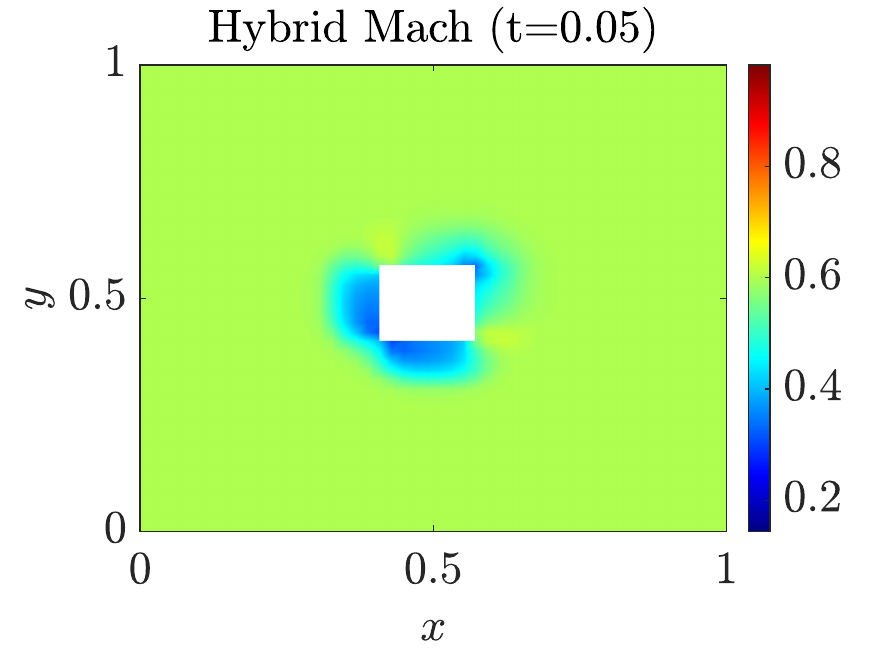}
    \end{subfigure}
    \hfill
    \begin{subfigure}{0.32\textwidth}
        \centering
        \includegraphics[width=\linewidth]{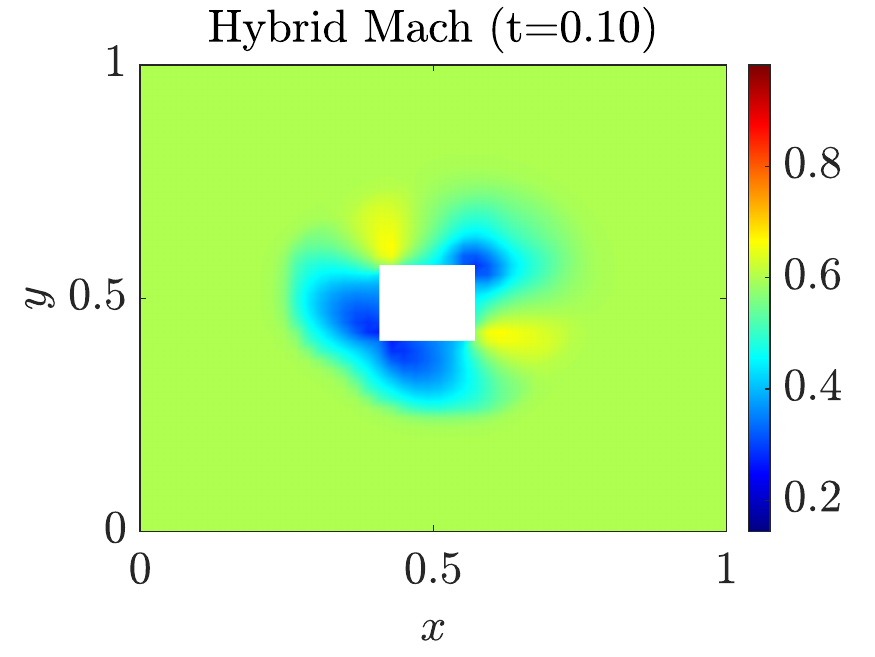}
    \end{subfigure}
    \hfill
    \begin{subfigure}{0.32\textwidth}
        \centering
        \includegraphics[width=\linewidth]{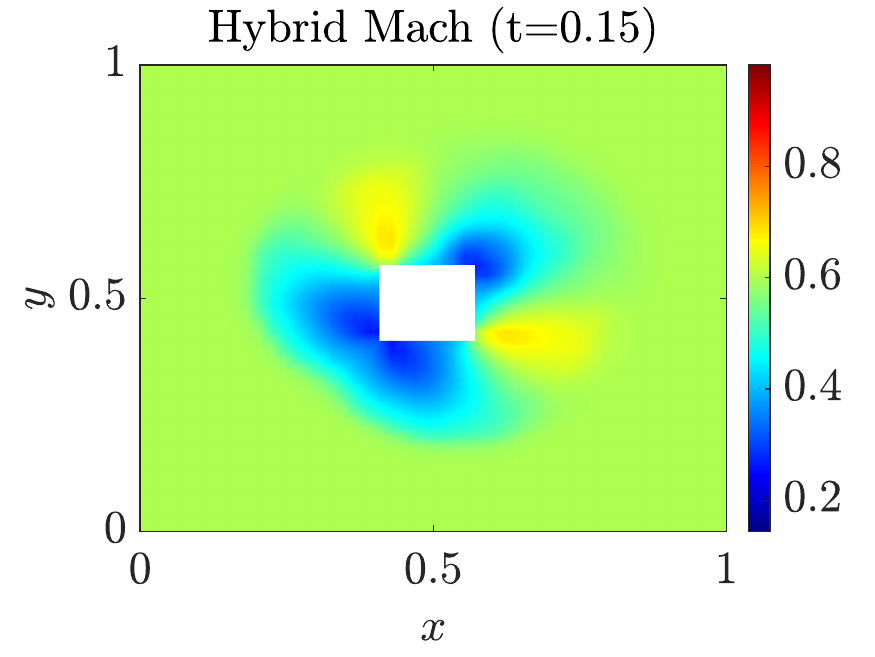}
    \end{subfigure}

    \begin{subfigure}{0.32\textwidth}
        \centering
        \includegraphics[width=\linewidth]{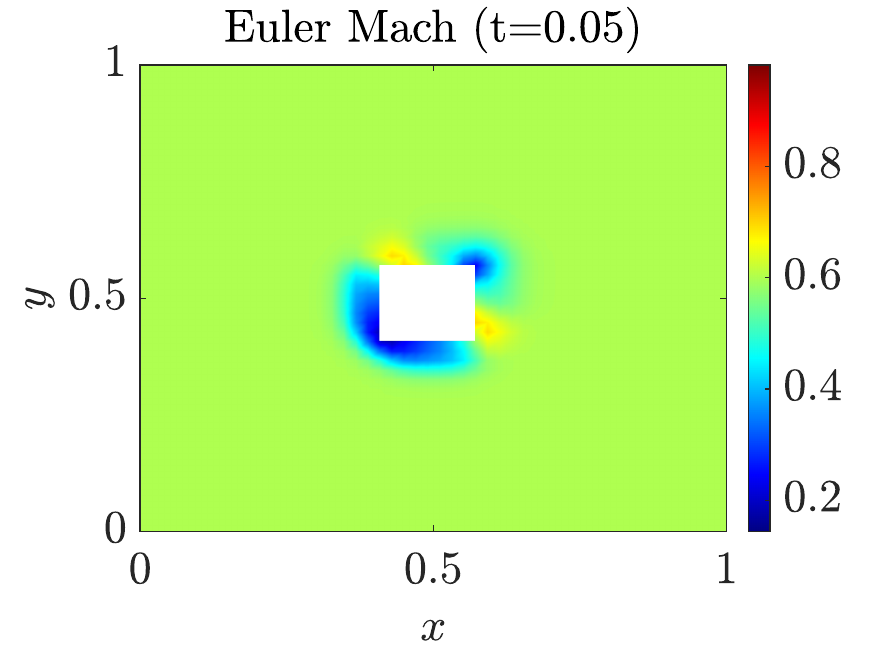}
    \end{subfigure}
    \hfill
    \begin{subfigure}{0.32\textwidth}
        \centering
        \includegraphics[width=\linewidth]{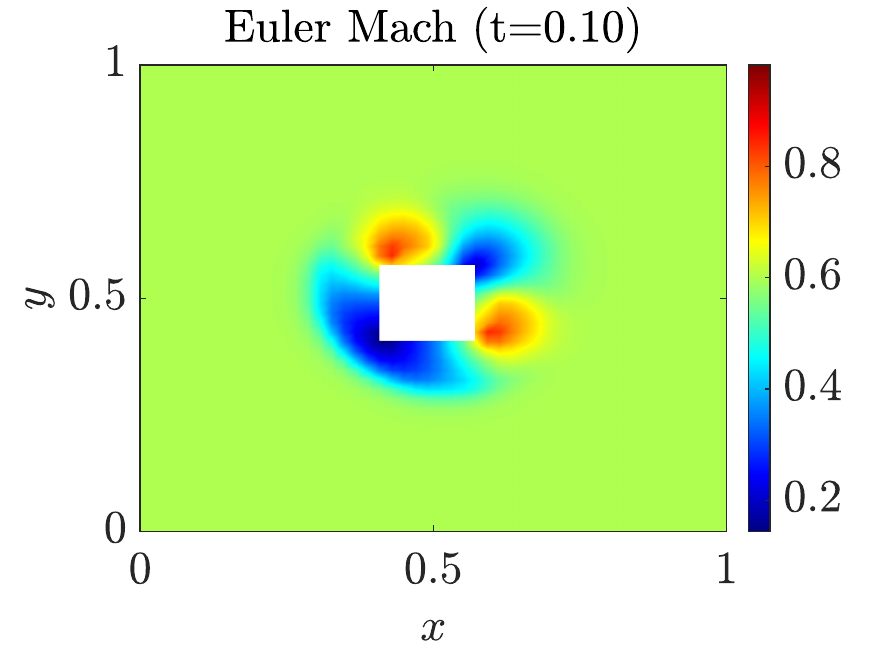}
    \end{subfigure}
    \hfill
    \begin{subfigure}{0.32\textwidth}
        \centering
        \includegraphics[width=\linewidth]{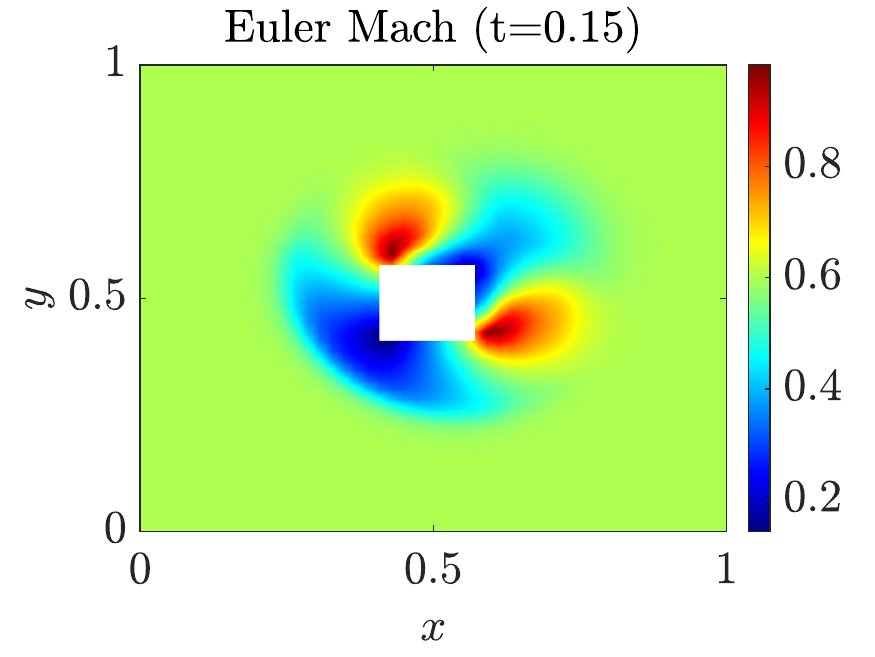}
    \end{subfigure}

    \caption{\hyperref[sec:test3]{Test 3}: Domain adaptation and Mach number evolution for the mixture. The top row illustrates the dynamic domain partition. The second, third, and fourth rows display the Mach number distribution obtained using the full kinetic solver, the hybrid solver, and the full Euler solver, respectively.}
    \label{Fig5-Mach}
\end{figure}

Figures~\ref{Fig3-Density-Species1} and~\ref{Fig4-Density-Species2} show the time evolution of the domain partition and of the mass density of the first and of the second species, respectively, for the hybrid solver, the full Euler solver and the full kinetic solver. Figure~\ref{Fig5-Mach} shows the partition and the Mach number of the mixture, defined as $\text{Ma}=|u|/\sqrt{\gamma P/\rho}$, where $|u|$ is the magnitude of the velocity of the mixture, $P$ its total pressure, $\rho$ its total density, and $\gamma=1+2/d_v=2$ the heat capacity ratio. The same colour scale is used in all the plots.

The hybrid solution is very close to the kinetic one. The kinetic solver is activated in a region around the obstacle, which grows in time as the perturbation produced by the obstacle spreads. The Euler solution differs visibly from the kinetic one, in particular around the corners of the obstacle.

The computing times are $12.37$ s for the full Euler solver, $323.4$ s for the full kinetic solver and $65.32$ s for the hybrid scheme, which is about $5$ times faster than the kinetic solver.

Furthermore, we evaluate, as in \cite{PaperA}, the space-time averaged $L^1$ error with respect to the full kinetic solution,
\begin{equation}\label{eq:errST}
E(q)=\frac{1}{t_f\,|\Omega|}\sum_{n}\Delta t^n\sum_{c}\Delta x\,\Delta y\,\big|q_c^n-q_c^{n,\mathrm{kin}}\big|,
\end{equation}
where the inner sum runs over the cells occupied by the gas, $|\Omega|$ is their total area, and $q$ is the total mass density, the temperature or the velocity magnitude of the mixture. 

The Euler solver gives errors of $1.9\times10^{-2}$ for the mass density, $1.0\times10^{-2}$ for the temperature and $8.4\times10^{-3}$ for the velocity magnitude; the hybrid solver reduces them to $3.1\times10^{-3}$, $3.2\times10^{-3}$ and $1.6\times10^{-3}$, respectively. In this test the hybrid scheme is therefore cheaper than the kinetic solver and more accurate than the Euler solver.

\subsection{Test 4: Two cavities}
\label{sec:test4}
In this test we consider a non-trivial geometry, with two cavities connected by a channel.

The physical domain $[0,1]^2$ is discretized using 50 cells in both directions. Two walls, attached to the upper and lower boundaries, divide the domain into two cavities: each wall is 15 cells wide along $x$, centred in the domain, and 18 cells high along $y$. The system consists of a binary gas mixture with a mass ratio $m_2/m_1 = 2$. The fluid is initially at rest and in local thermodynamic equilibrium, with different states for $x<0.5$ and $x\ge 0.5$. In the macroscopic cells, the Euler solution is initialized with the following macroscopic moments $W_1=(\rho_1, u^x, u^y, T)$ and $W_2=(\rho_2, u^x, u^y, T)$ for the first and second species, respectively:
\begin{equation}
    \begin{cases}
        W_1(x<0.5,y) = ((1-\alpha )\rho_L,\,0,\,0,\,T_L),\\
        W_2(x<0.5,y) = (\alpha \rho_L,\,0,\,0,\,T_L),\\
        W_1(x\geq0.5,y) = (\alpha \rho_R,\,0,\,0,\,T_R),\\
        W_2(x\geq0.5,y) = ((1-\alpha)\rho_R,\,0,\,0,\,T_R),
    \end{cases}    
\end{equation}
where $\alpha=10^{-3}$, $\rho_L=1$, $\rho_R=1$, $T_L=P_L/((1-\alpha )\rho_L/m_1)$, $T_R=P_R/((1-\alpha )\rho_R/m_2)$, $P_L=1$ and $P_R=0.8$. In the kinetic cells, the distribution functions $f_1$ and $f_2$ for both species are initialized as Maxwellian distributions with moments $W_1(x, y)$ and $W_2(x,y)$. The velocity space $[-8,8]^2$ is discretized with $32$ points along both directions.

We set $\varepsilon=0.1$ and $\eta=0.03\,\varepsilon$, and the simulation is run up to $t_f=0.25$. Reflective slip boundary conditions are imposed on all domain edges and on the two walls separating the two cavities.

\begin{figure}[htbp]
    \centering
    \begin{subfigure}{0.32\textwidth}
        \centering
        \includegraphics[width=\linewidth]{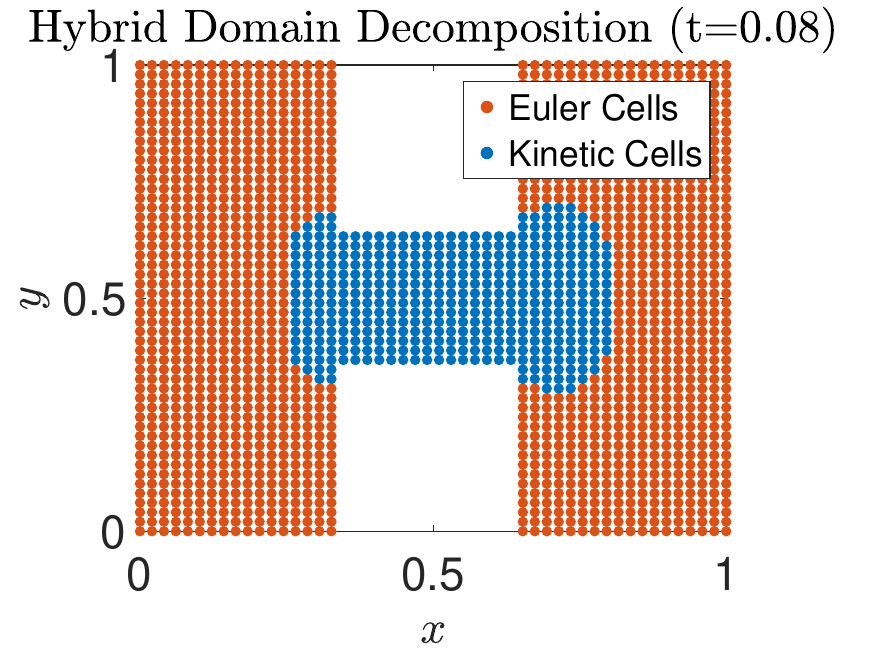}
    \end{subfigure}
    \hfill
    \begin{subfigure}{0.32\textwidth}
        \centering
        \includegraphics[width=\linewidth]{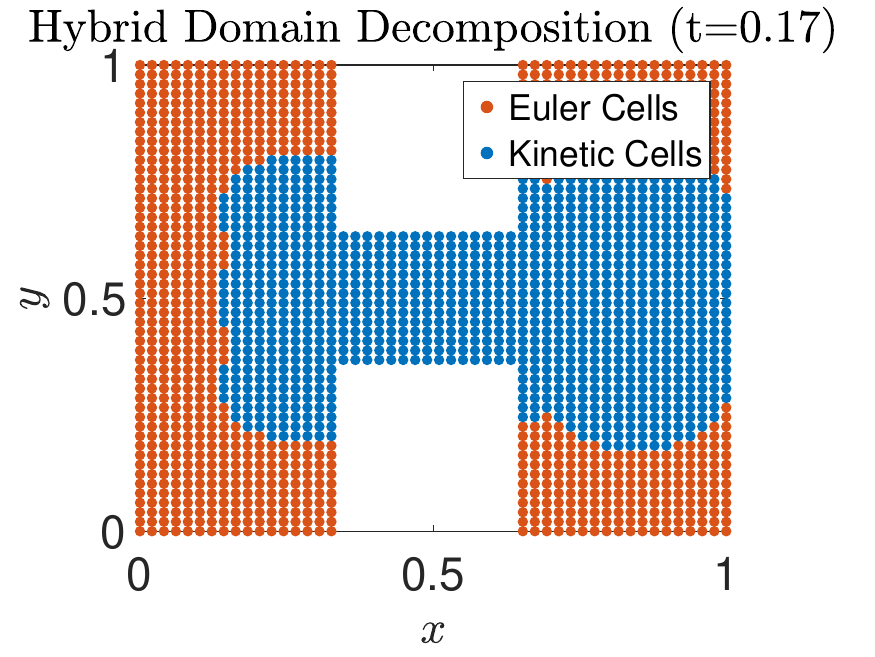}
    \end{subfigure}
    \hfill
    \begin{subfigure}{0.32\textwidth}
        \centering
        \includegraphics[width=\linewidth]{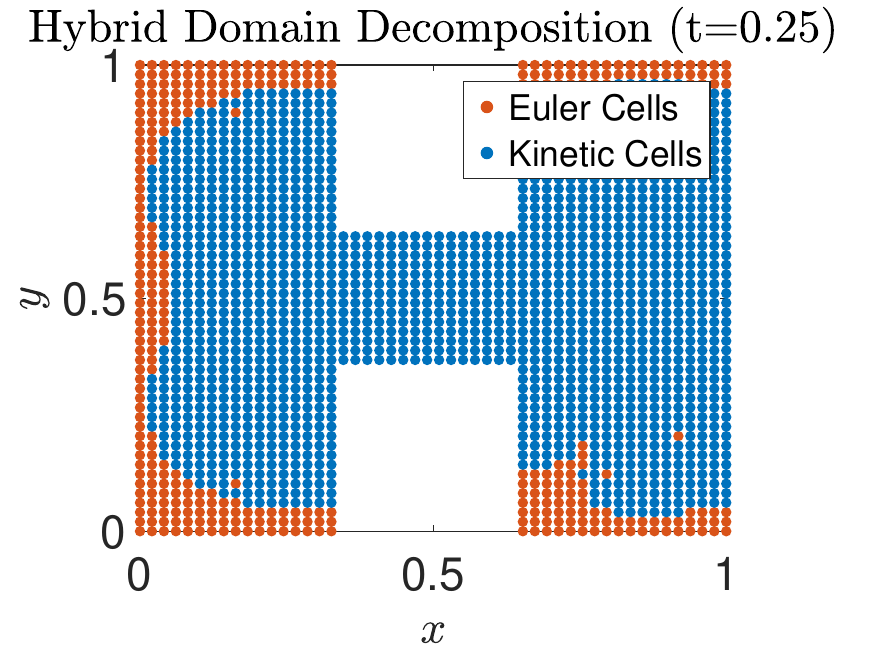}
    \end{subfigure}

    \begin{subfigure}{0.32\textwidth}
        \centering
        \includegraphics[width=\linewidth]{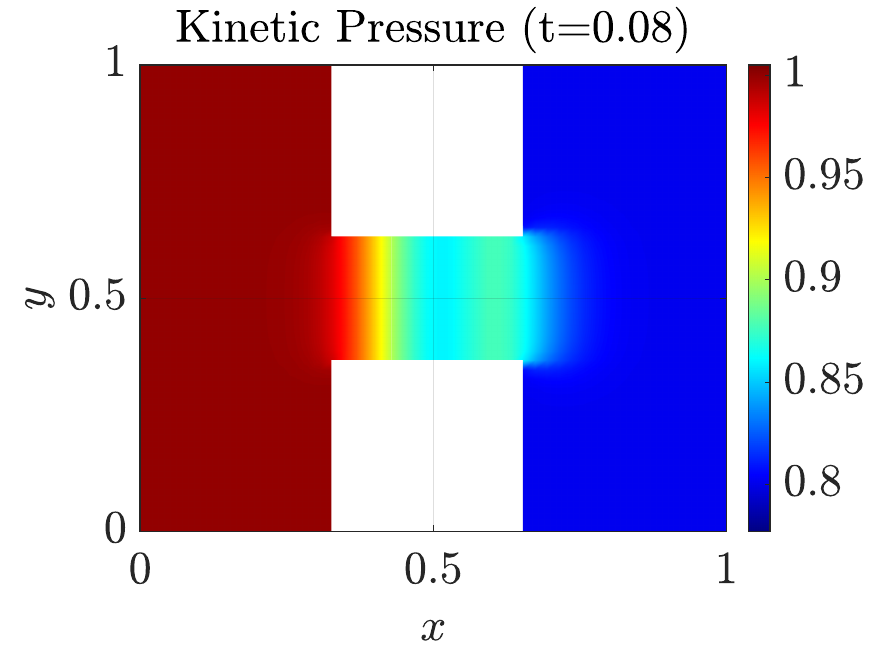}
    \end{subfigure}
    \hfill
    \begin{subfigure}{0.32\textwidth}
        \centering
        \includegraphics[width=\linewidth]{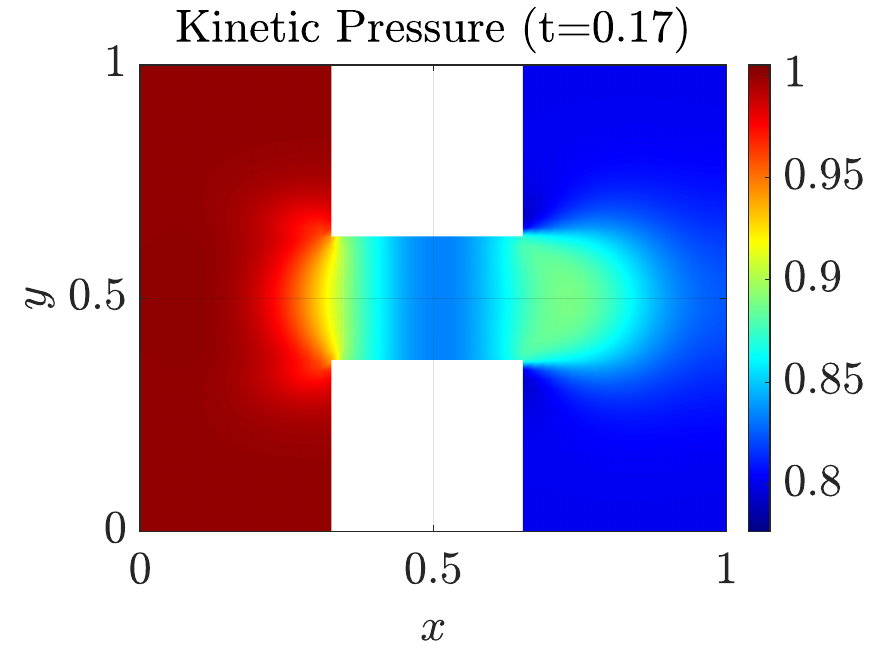}
    \end{subfigure}
    \hfill
    \begin{subfigure}{0.32\textwidth}
        \centering
        \includegraphics[width=\linewidth]{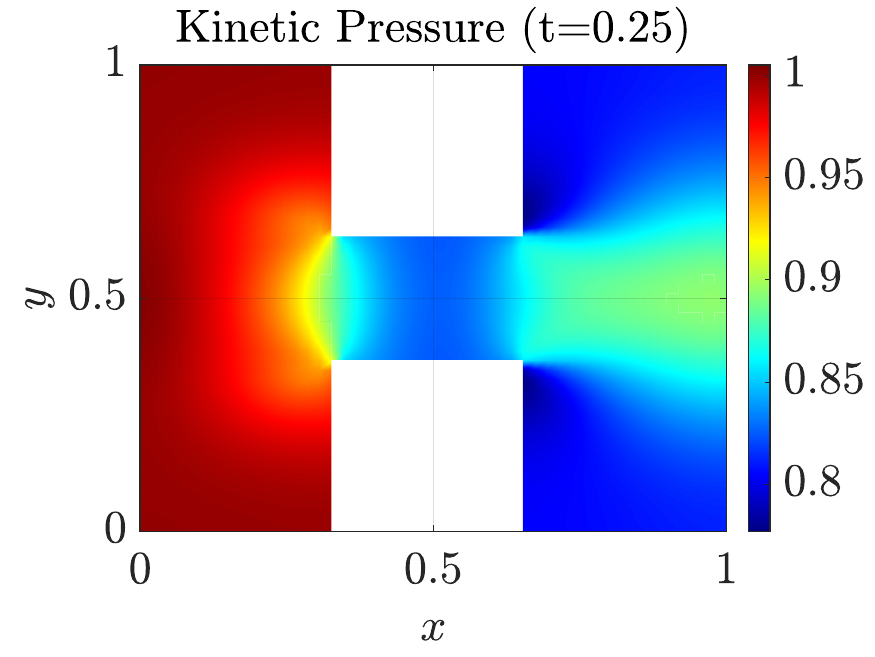}
    \end{subfigure}

    \begin{subfigure}{0.32\textwidth}
        \centering
        \includegraphics[width=\linewidth]{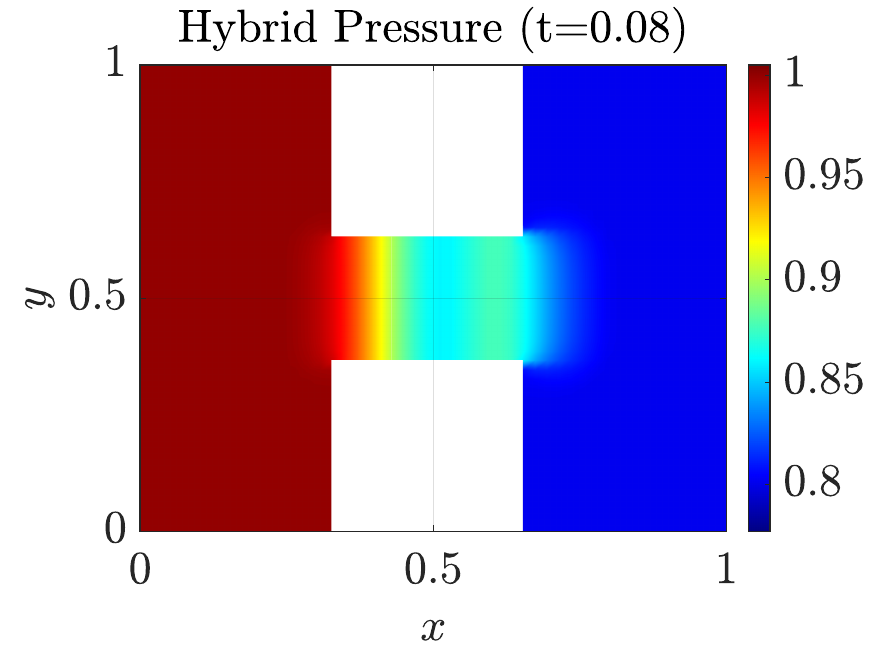}
    \end{subfigure}
    \hfill
    \begin{subfigure}{0.32\textwidth}
        \centering
        \includegraphics[width=\linewidth]{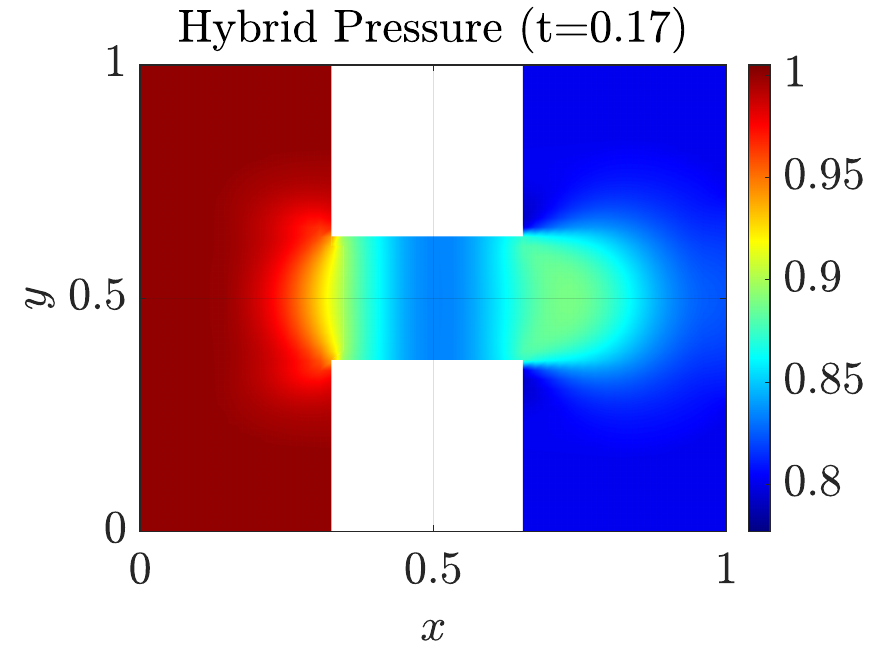}
    \end{subfigure}
    \hfill
    \begin{subfigure}{0.32\textwidth}
        \centering
        \includegraphics[width=\linewidth]{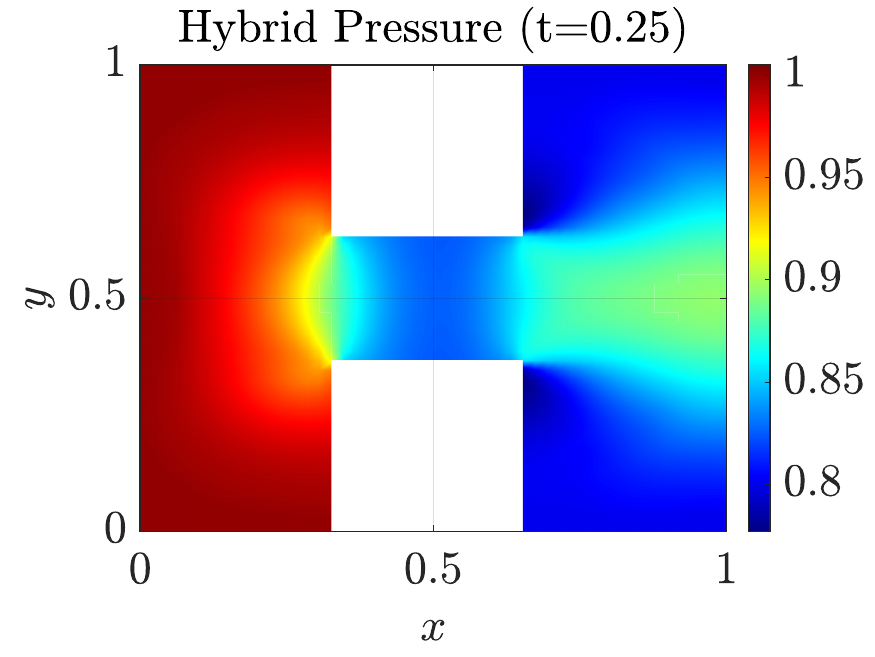}
    \end{subfigure}

    \begin{subfigure}{0.32\textwidth}
        \centering
        \includegraphics[width=\linewidth]{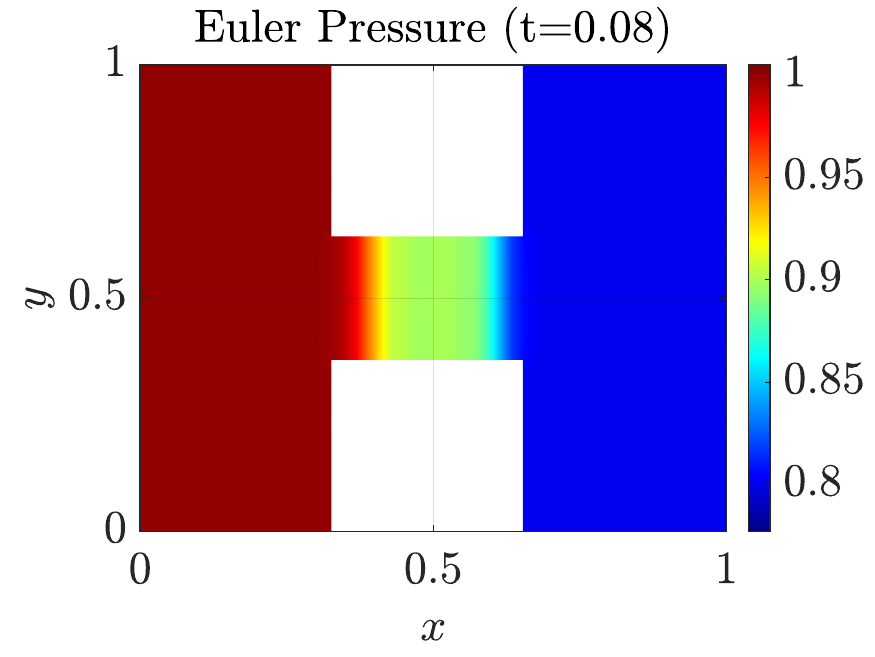}
    \end{subfigure}
    \hfill
    \begin{subfigure}{0.32\textwidth}
        \centering
        \includegraphics[width=\linewidth]{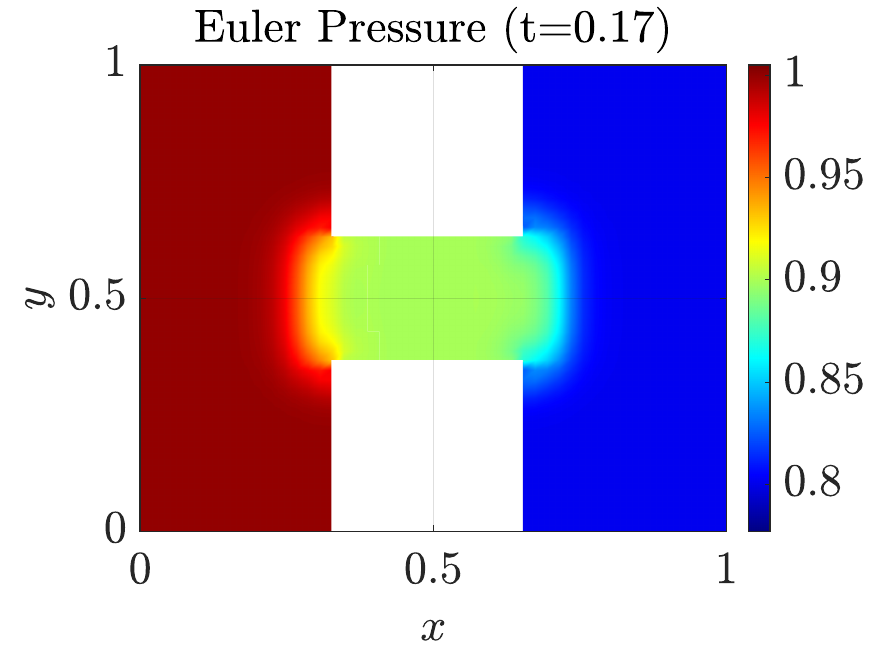}
    \end{subfigure}
    \hfill
    \begin{subfigure}{0.32\textwidth}
        \centering
        \includegraphics[width=\linewidth]{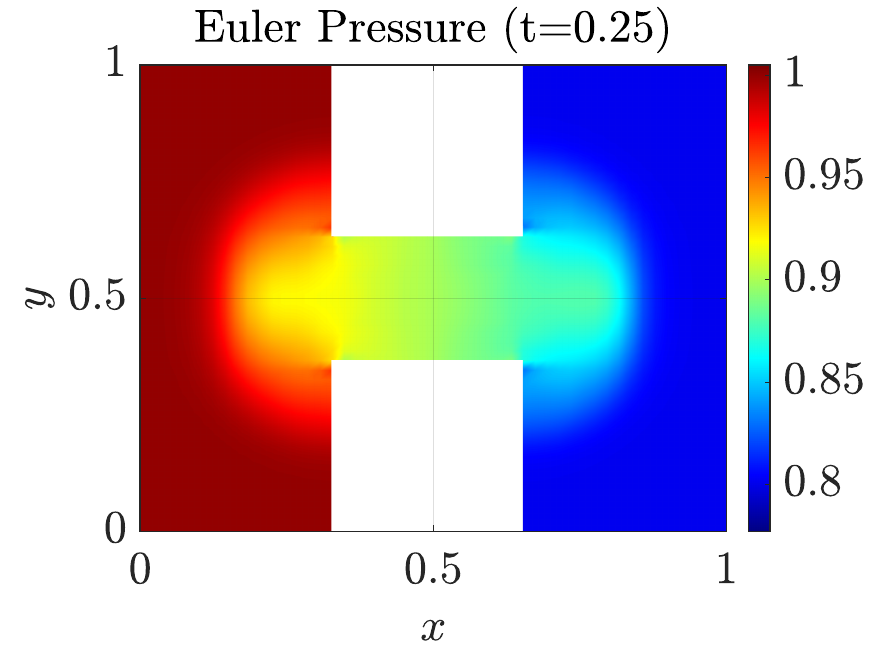}
    \end{subfigure}

    \caption{\hyperref[sec:test4]{Test 4}: Domain adaptation and total pressure evolution for the mixture. The top row illustrates the dynamic domain partition. The second, third, and fourth rows display the pressure distribution obtained using the full kinetic solver, the hybrid solver, and the full Euler solver, respectively.}
    \label{Fig6-P}
\end{figure}

Figure~\ref{Fig6-P} illustrates the time evolution of the dynamic domain partition and of the total pressure of the mixture. The results of the hybrid scheme, the full Euler solver and the full kinetic solver are shown with the same colour scale. The hybrid scheme activates the kinetic solver in regions where non-equilibrium effects appear, and the kinetic region grows gradually in time. The Euler solution differs significantly from the kinetic one, while the hybrid solution is very close to it in the whole domain.

The computing times are $65.75$ s for the full Euler solver, $1279$ s for the full kinetic solver and $576.4$ s for the hybrid scheme, which is about $2.2$ times faster than the kinetic solver. The space-time averaged errors \eqref{eq:errST} of the Euler solver are $3.8\times10^{-3}$ for the mass density, $2.0\times10^{-2}$ for the temperature and $5.6\times10^{-3}$ for the velocity magnitude; the hybrid solver reduces them by more than one order of magnitude, to $1.4\times10^{-4}$, $7.2\times10^{-4}$ and $1.8\times10^{-4}$, respectively.

The switching thresholds are specified separately for each benchmark. The reported errors and speedups therefore describe these choices.
\section{Discussion on numerical efficiency}\label{secNumEff}
In this section, we summarize and discuss the overall numerical efficiency of the proposed dynamic domain decomposition strategy for gas mixtures. Table \ref{table1} reports the execution times, the achieved computational speedups, and the space-time averaged $L^1$ density errors for the shock tube problem (\hyperref[sec:test2]{Test 2}), the flow over a square (\hyperref[sec:test3]{Test 3}), and the two cavities simulation (\hyperref[sec:test4]{Test 4}).

The error metric is the space-time averaged $L^1$ error \eqref{eq:errST} of the total mass density with respect to the full kinetic solution. We do not report it for the shock tube problem (\hyperref[sec:test2]{Test 2}). For $\varepsilon=10^{-4}$ all three solvers approximate the solution of the Euler equations, so the kinetic solution is not a better reference than the Euler one.

In the two rarefied tests the hybrid scheme reduces the density error of the Euler solver by a factor of about $6$ (\hyperref[sec:test3]{Test 3}) and $28$ (\hyperref[sec:test4]{Test 4}), since the kinetic solver is used where non-equilibrium effects are significant. 

For mixtures with large mass ratios the cost of kinetic solvers grows quickly, since the velocity grid has to cover the thermal speed of the lightest species and, at the same time, resolve the Maxwellian of the heaviest one; in this case confining the kinetic solver to a small region is even more important than for single-species frameworks. In this context, the proposed adaptive algorithm reduces this computational cost. In our tests, run on CPU, the hybrid solver is between $2.2$ and $5$ times faster than the kinetic solver.

In the rarefied tests the hybrid solver stays close to the kinetic one, at a fraction of its cost.

\begin{table}[htbp]
    \centering
    \begin{tabular}{|lcc|cc|}
        \toprule
        \textbf{Simulation} & \textbf{Solver} & \textbf{Density Error ($L^1$)} & \textbf{Execution Time} & \textbf{Time / hybrid time} \\
        \midrule
        \multirow{3}{*}{Shock Tube}   & Kinetic & --                    & 2352 s   & 3.7  \\
                                      & Hybrid  & --                    & 628.7 s  & 1.0  \\
                                      & Euler   & --                    & 107.3 s  & 0.17 \\
        \midrule
        \multirow{3}{*}{Square}       & Kinetic & --                    & 323.4 s  & 5.0  \\
                                      & Hybrid  & $3.1\times 10^{-3}$   & 65.32 s  & 1.0  \\
                                      & Euler   & $1.9\times 10^{-2}$   & 12.37 s  & 0.19 \\
        \midrule
        \multirow{3}{*}{Two cavities} & Kinetic & --                    & 1279 s   & 2.2 \\
                                      & Hybrid  & $1.4\times 10^{-4}$   & 576.4 s  & 1.0  \\
                                      & Euler   & $3.8\times 10^{-3}$   & 65.75 s  & 0.11 \\
        \bottomrule
    \end{tabular}
    \caption{CPU execution times, ratio between the time of each solver and that of the hybrid solver, and space-time averaged $L^1$ errors \eqref{eq:errST} of the total mass density with respect to the full kinetic solution, for the shock tube problem (\hyperref[sec:test2]{Test 2}), the flow over a square (\hyperref[sec:test3]{Test 3}) and the two cavities (\hyperref[sec:test4]{Test 4}). Errors are omitted for the shock tube test, where all models converge to the same hydrodynamic solution.}
    \label{table1}
\end{table}

\section{Conclusion}\label{secCon}
In this work we have proposed a dynamic domain decomposition method for multiscale gas mixtures in which all species collide on the same time scale. The switching criterion, inspired by the works of Levermore, Morokoff, and Nadiga \cite{lev-mor-nad-1998}, Filbet and Rey \cite{fil-rey-2014}, Filbet and Xiong \cite{fil-xio-2018}, Tiwari \cite{tiw-2000}, Li, Song, and Wang \cite{li-son-wan-2021}, Xiong and Qiu \cite{xio-qiu-2017}, Caparello, Pareschi, and Rey \cite{cap-par-rey-2025, cap-par-rey-2026}, and Caparello and Tenna \cite{cap-ten-2025}, and applied species by species to mixtures in \cite{PaperA}, uses truncated Chapman--Enskog expansions to decide where the fluid description has to be replaced by the kinetic one. It only uses macroscopic quantities, so the distribution functions are not needed in the macroscopic regions.

With respect to \cite{PaperA}, where collisions between particles of the same species dominate and the partition is computed species by species, the novelty of the present work is the single-velocity, single-temperature regime: the indicator is derived for the whole mixture, and the dynamic spatial partition is common to all species.
The decomposition uses two levels: the Euler equations of the mixture in the regions where the species are close to local thermodynamic equilibrium, and the multi-species BGK model in the rarefied, non-equilibrium regions. These two descriptions are dynamically coupled across the computational domain.

By combining high-order, Asymptotic-Preserving (AP) IMEX time integration with finite volume techniques on Cartesian meshes, the method remains stable over the wide range of Knudsen numbers considered in the tests. When the partition is fixed, the stage-wise coupling preserves the third order in time, and the numerical convergence test shows no temporal order reduction at the sampled Knudsen numbers. In our tests the coupling does not produce oscillations at the interfaces.

In the numerical experiments the kinetic regions are selected automatically and change in time with the flow. The hybrid solver stays close to the kinetic solution while achieving computational speedups between 2.2 and 5, depending on the size of the kinetic region and on the parameters of the simulation.

These results show that the cost of deterministic kinetic simulations of gas mixtures can be reduced by restricting the kinetic solver to the non-equilibrium regions, also in non-trivial geometries. This could make kinetic models easier to use in aerodynamic and aerospace engineering problems involving multicomponent flows. Future work will address the extension of this framework to the full spectral Boltzmann collision operator for gas mixtures, the inclusion of polyatomic gases and chemical reactions, and a parallel implementation for three-dimensional problems.

\section*{Acknowledgments}
    The authors received funding from the European Union's Horizon Europe research and innovation program under the Marie Skłodowska-Curie Doctoral Network DataHyking (Grant No. 101072546). 
    DC and TR are also supported by the French government, through the UniCA$_{JEDI}$ Investments in the Future project managed by the National Research Agency (ANR) with the reference number ANR-15-IDEX-01.
    The research of LP has been supported by the Royal Society under the Wolfson Fellowship “Uncertainty quantification, data-driven simulations and learning of multiscale complex systems governed by PDEs”. LP also acknowledges the  support by the Italian Ministry of University and Research (MUR) through the Fondo Italiano per la Scienza 2022–2023 (Bando FIS 2), under Decree No. 1236 of 01/08/2023, within the FIS2023-01334 Advanced Grant “Tackling Complexity: Advanced Numerical Approaches for Multiscale Systems with Uncertainties” (ADAMUS). This work has been written within the activities of GNCS group of INdAM (Italian National Institute of High Mathematics).

\section*{Declaration of generative AI and AI-assisted technologies}
During the preparation of this work the authors used a generative AI tool to improve the language and readability of the manuscript. After using this tool, the authors reviewed and edited the content as needed and take full responsibility for the content of the publication.
\bibliographystyle{acm}
\bibliography{biblio}
  
\end{document}